\documentclass[12pt,reqno]{amsart}
\usepackage{amsthm}
\usepackage{cases}
\usepackage{amscd}
\usepackage{amsfonts}

\usepackage{amssymb}
\usepackage{amsmath}
\usepackage{booktabs}
\usepackage{hyperref}

\usepackage{enumitem}

\usepackage{subfigure}
\usepackage{algorithm}
\usepackage{algpseudocode}
\renewcommand{\algorithmicrequire}

\usepackage{cleveref}
\usepackage{color}
\usepackage{datetime}

\allowdisplaybreaks

\makeatletter
\newcommand\figcaption{\def\@captype{figure}\caption}
\newcommand\tabcaption{\def\@captype{table}\caption}
\makeatother

\usepackage{geometry}
\usepackage{algorithm}
\usepackage{multirow}
\usepackage{mathrsfs}

\numberwithin{equation}{section}  %

\usepackage{hyperref}
\usepackage{cleveref}

\newcommand{\beq}{\begin{equation}}
\newcommand{\eeq}{\end{equation}}
\newcommand{\bea}{\begin{eqnarray}}
\newcommand{\eea}{\end{eqnarray}}
\newcommand{\beas}{\begin{eqnarray*}}
\newcommand{\eeas}{\end{eqnarray*}}

 \newtheorem{remark}{Remark}[section]
 \newtheorem{theorem}{Theorem}[section]
\newtheorem{lemma}{Lemma}[section]
\newtheorem{proposition}[theorem]{Proposition}
\newtheorem{corollary}[theorem]{Corollary}

\begin{document}
\title[Small inertia limit of the distribution-dependent stochastic  Cucker--Smale equation ]{Small inertia limit of the distribution-dependent  stochastic  Cucker--Smale equation} 

\author[J.Z. Huang, Y. Lv and W. Wang]{Jizhao Huang$^{1}$,  Yan Lv$^{1, *}$ and Wei Wang$^{2}$} 
\address{$^{1}$ School of Mathematics and Statistics, Nanjing University of Science and Technology, Nanjing, China}

\address{$^{2}$ School of Mathematics, Nanjing University, Nanjing, China}

\email{hjzmath21@163.com(J.Z. Huang),  lvyan1998@aliyun.com(Y. Lv), wangweinju@nju.edu.cn(W. Wang)}

\thanks{$^*$Corresponding author.}



\begin{abstract}
This paper investigates the small inertia limit of the distribution-dependent stochastic Cucker--Smale equation with the full velocity alignment source term, which serves as the mean-field equation for the stochastic Cucker--Smale model.
Inspired by the idea of Hottovy et al. \cite{Hottovy2015} to rearrange SDEs via integration by parts,  we repeatedly perform integration by parts on independent copies of our equation.  
 Then we define a global position correction and a quadratic remainder, both with uniformly bounded tails.    Applying the Aldous  criterion, we obtain tightness of the position correction, which in turn yields tightness of the position process.
Subsequently, using the averaging method on the fast time scale, the quadratic functional   can be decomposed into the equilibrium average, local quadratic fluctuations, and exponential transients.
Then, the Poisson--Lyapunov equation is employed to prove the convergence of the time integral.
Finally, under various suitable assumptions, we obtain the weak convergence of the laws of the approximating processes to the law of the solution of the small-inertia limit equation, along with the identification of the limiting process, and also the convergence rate of  strong $L^1$-convergence.  The extra noise-induced drift term can be expressed as the solution of a contraction-type Fredholm equation.
This paper is   complementary to the deterministic particle limit work of  \cite{Carrillo2021,Choi2021} . 
\end{abstract}

\maketitle


\tableofcontents

\section{Introduction}
 
 The Smoluchowski--Kramers approximation characterizes the limiting dynamics of stochastic systems in the small-inertia regime, where  inertial effects are negligible compared to friction. It establishes a rigorous link between second-order SDEs and their effective first-order overdamped counterparts. In the singular limit, the velocity acts as a fast variable and rapidly relaxes to equilibrium, while the position evolves on the slow time scale; this separation of scales gives rise to an effective first-order SDE. This approximation is also known as the small-mass limit; see, e.g., \cite{Cerrai2006,Cerrai2022,Freidlin2004,Shi2021}. 
 This limit has been rigorously characterized from different perspectives in the existing literature: at the trajectory level, the position process is shown to converge, in probability or in the moment sense, to solutions of the overdamped equation \cite{Brzezniak2025}; at the level of probability distribution, the convergence of the Fokker--Planck equation for the joint position-velocity probability distribution to the Smoluchowski equation has been studied \cite{Choi2022}; and at the macroscopic level, aggregation-diffusion equations can also be rigorously derived from compressible Euler-type equations in the large-friction limit \cite{Carrillo2019}.

Over the past few decades, the Smoluchowski--Kramers approximation has been extensively studied. For constant friction coefficients, the classical theory has become relatively mature; nevertheless, recent studies have produced several results in various directions, in particular concerning convergence rates and  the effects of   noise.
For instance, in \cite{Son2020}, by using  techniques of Malliavin calculus,  Son   gave  the rate of convergence for the Smoluchowski--Kramers
approximation for SDEs with fractional Brownian   motion.
In \cite{Liu2026}, Liu  et  al. employed Stein's method to establish the long-time convergence rate of the Smoluchowski--Kramers approximation for   SDE.
Zhao et al. \cite{Zhao2026} used rough path theory to study the Smoluchowski--Kramers approximation of random differential equations  with fast oscillating noise.
Moreover, for state-dependent friction, Hottovy et al. \cite{Hottovy2015}, employing the theory of convergence of stochastic integrals developed by Kurtz and Protter, studied the Smoluchowski--Kramers limit of the following SDEs with state-dependent friction 
  \begin{equation*} 
\left\{
\begin{aligned}
\mathrm{d}   X_t^\epsilon= {}&  V_t^\epsilon   \mathrm{d}t,\\
\epsilon \mathrm{d}V_t^\epsilon
= {} &\mathbf{F}(X_t^\epsilon)\mathrm{d}  t-\gamma(X_t^\epsilon)V_t^\epsilon   \mathrm{d}t
+\sigma(X_t^\epsilon)   \mathrm{d}  B_t,
\end{aligned}
\right.
\end{equation*}
where  $\mathbf{F} $  is  deterministic part of the force,  $\gamma$ is  invertible matrix-valued function, $\sigma $ is the intensity of the noise,  and $\{B_t\}_{t}$ is a Wiener process.
They found that the limiting equation contains an additional noise-induced drift term, which is associated with the coupling between state-dependent friction and noise, and characterized this additional drift as the solution of a Lyapunov matrix equation.  
Specifically, the limiting SDE takes the form
  \begin{equation*} 
\begin{aligned}
  \mathrm{d}X_t 
={}& \gamma^{-1}(X_t )  \mathbf{F}(X_t )\mathrm{d}  t   +\mathbf{S}(X_t)\mathrm{d}  t
+\gamma^{-1}(X_t )   \sigma(X_t )   \mathrm{d}  B_t,
\end{aligned}
\end{equation*}
where $\mathbf{S}(X_t)$ denotes the noise-induced drift, whose $i$-th component is given by
$$
\mathbf{S}_{i}(x)   =   \sum_{j,k=1}^{d}    \frac{\partial}{    \partial x_{k}   }  [   ( \gamma^{-1}  )_{ij}  (x)  ] \mathbf{J}_{jk}(x),
$$
and $\mathbf{J}$ is the matrix solving the Lyapunov equation
$$
\mathbf{J}  (x) \gamma^{\top} (x)  +    \gamma  (x)  \mathbf{J} (x) =  \sigma (x) \sigma^{\top} (x).
$$
 Subsequently, the Smoluchowski--Kramers approximation for SDEs with state-dependent friction has been studied  by various researchers. 
 For instance, in   \cite{Xie2026}, Xie and Zhang
investigated the uniform-in-time diffusion approximation for coupled multiscale stochastic systems, the result that can be applied to derive a uniform Smoluchowski--Kramers approximation for degenerate Langevin systems with state-dependent friction.
Blassel \cite{Blassel2026} studied the Smoluchowski--Kramers approximation for underdamped equilibrium Langevin dynamics  using  $L^2$-hypocoercivity  
 and provided a clear explanation of the noise-induced drift term.
In    \cite{Wang2024},  Wang et al.   investigated the Smoluchowski--Kramers approximation of a Langevin equation with state-dependent friction by  using   the idea of averaging approach \cite[Chapter 5]{Duan2014}  and   coupling the microscopic dynamics to the macroscopic scales.
For related literature on the Smoluchowski--Kramers approximation of SDEs, see e.g.  \cite{Cerrai2020,Cerrai2026,Qian2026,Xie2022} and the references therein.
 
In parallel, the approximation has been established for distribution-dependent SDEs, also known as McKean--Vlasov SDEs, which can be viewed as the mean-field limit of interacting particle systems.
In  \cite{Liu2024,Son2024}, the authors obtained the convergence rates of the Smoluchowski--Kramers approximation for distribution-dependent SDEs driven by Brownian motion and by fractional Brownian motion, respectively, by means of Malliavin calculus.
In  \cite{Li2025},  the authors  studied   the Smoluchowski--Kramers approximation
for a class of McKean--Vlasov SDEs.

In nature, biological groups such as bird flocks  and  fish schools  exhibit rich self-organized collective motions, including flocking, swarming, and schooling. From a biological perspective, individuals have inherent inertia or response delays when adjusting their movement directions. Their perception of and reaction to environmental information and neighboring individuals are not instantaneous but are limited by physiological constraints and nerve conduction velocities. Therefore, the inertia or response time of an individual is typically small. Mathematically, this feature can be described by introducing a small parameter $\epsilon>0$, which represents the inertia coefficient or the response time scale of an individual. When $\epsilon$ is small, the system is in the small inertia or near overdamped regime, where the velocity of each individual rapidly relaxes to an equilibrium state determined by the influences of surrounding individuals. The Cucker--Smale model is a classical mathematical model for describing such self-organized behavior.
The model, proposed by Cucker and Smale \cite{Cucker2007} in 2007, can be described as a second-order Newtonian particle system with a velocity alignment mechanism involving a communication weight function.
In \cite{Wang2022}, due to the difficulties arising from the full velocity alignment term, the authors, under a constant communication weight, used the averaging method to obtain the Smoluchowski--Kramers limit of the stochastic $N$-particle system within the mean-field limit framework.
In  \cite{Wang-Lv-Zhang-Duan-Wang2024},  also under a constant communication weight function, Wang et al. established the commutation of the mean-field limit and the Smoluchowski--Kramers  limit for an interacting particle system with non-Gaussian L\'{e}vy noise. In  \cite{Shi2024}, using the averaging method on a non-autonomous slow-fast system that couples microscopic and macroscopic scales, the authors obtained the Smoluchowski--Kramers approximation for McKean--Vlasov SDEs with a damping term in the velocity alignment.  
 In \cite{Carrillo2021}, Carrillo and Choi considered the monokinetic density
$
\mathbf{f}^\epsilon(x,v)=\varrho^\epsilon (x) \delta_{\mathfrak{u}^\epsilon (x)}(v),
$
and assumed that the communication weight  function $\phi$  is bounded and Lipschitz. They treated both the regular case, in which the interaction potential satisfies
$
\nabla K \in \mathcal{W}^{1,\infty}(\mathbb{R}^d),
$
and the singular repulsive case, in which $K$ is a Coulomb/Riesz potential.
  In both settings,      they  rigorously derive the pressureless Euler-type equations \eqref{eq:pressureless-hydrodynamic-Euler-alignment-models} with nonlocal dissipative terms in velocity from the Cucker--Smale model \eqref{eq:deterministic-Cucker--Smale-model} with confinement and interaction potentials.  These equations read as  
\begin{equation} \label{eq:deterministic-Cucker--Smale-model}
\left\{
\begin{aligned}
\dot{X}_t^{\epsilon,i,N}=   {}  &   V_t^{\epsilon,i,N}    ,  \qquad i=1,\dots, N   , \\
\epsilon_{N}  \dot{V}_t^{\epsilon,i,N}
= {}  &    -   \gamma  V_t^{\epsilon,i,N}  - \nabla \mathcal{V}(X_t^{\epsilon,i,N})
- \frac{1}{N} \sum_{j=1}^{N} \nabla K(X_t^{\epsilon,i,N}   -X_t^{\epsilon,j,N})    
\\
    {} & +\frac{1}{N} \sum_{j=1}^{N}  \phi( X_t^{\epsilon,i,N}  -X_t^{\epsilon,j,N} )   (V_t^{\epsilon,j,N} -V_t^{\epsilon,i,N})     
\end{aligned}
\right.
\end{equation}
and 
\begin{equation} \label{eq:pressureless-hydrodynamic-Euler-alignment-models}
\left\{
\begin{aligned} 
 &\partial _{t}  \varrho  ^\epsilon    +     \nabla \cdot  (\varrho  ^\epsilon    \mathfrak{u}^\epsilon    )  =0 ,\\
 &   \epsilon  [ \partial _{t}  (  \varrho  ^\epsilon    \mathfrak{u}^\epsilon   )    + \nabla \cdot  (\rho  ^\epsilon    \mathfrak{u} ^\epsilon \otimes  \mathfrak{u}^\epsilon      )   ]     = - \gamma \varrho  ^\epsilon  \mathfrak{u} ^\epsilon    -   \varrho  ^\epsilon   \nabla  \mathcal{V}  -    \varrho  ^\epsilon   \nabla K  * \varrho  ^\epsilon   \\
   &    \qquad\qquad   \qquad\qquad  \qquad\qquad \qquad +  \varrho   \int_{\mathbb{R} ^{d}  } \phi(x-  x')   ( \mathfrak{u} ^\epsilon  (x') -\mathfrak{u} ^\epsilon  (x)     )   \varrho  ^\epsilon(  x'  )    \mathrm{d} x'  ,
\end{aligned}
\right.
\end{equation}
where $X_{t}^{\epsilon,i,N}\in\mathbb{R}^{d}$ and $V_{t}^{\epsilon,i,N}\in\mathbb{R}^{d}$ denote the position and velocity, respectively; $\gamma>0$ represents the strength of linear damping on the velocity; $\epsilon_{N}>0$ represents the strength of inertia; $\mathcal{V}$ and $K$ represent the confinement and interaction potentials, respectively; $\phi$ is a communication weight function; and
$
\varrho^\epsilon := \int_{\mathbb{R}^{d}} \mathbf{f}^\epsilon(x,v)\,\mathrm{d}v,
$  $
\varrho^\epsilon \mathfrak{u}^\epsilon := \int_{\mathbb{R}^{d}} v\,\mathbf{f}^\epsilon(x,v)\,\mathrm{d}v
$
represent the local particle density and local momentum, respectively.
Moreover,    under those assumptions,   Carrillo and Choi rigorously derived the small inertia limit from \eqref{eq:pressureless-hydrodynamic-Euler-alignment-models} to the following aggregation equations with nonlocal velocity fields 
\begin{equation*} 
\left\{
\begin{aligned} 
 &\partial_t\varrho +\nabla\cdot(\varrho      \mathfrak{u} )=0,\\
 &  \gamma  \varrho      \mathfrak{u} 
=-\varrho  \bigl(\nabla  \mathcal{V}+\nabla K*  \varrho  \bigr)
+\varrho  \int_{\mathbb{R} ^{d}  } \phi(x- x')\bigl( \mathfrak{u}(x')- \mathfrak{u}(x)\bigr)\varrho (  x' )  \mathrm{d}  x' .
\end{aligned}
\right.
\end{equation*}
 The velocity equation can be written as
$  (\gamma  +   \phi *   \varrho)   \mathfrak{u} 
 - \int_{\mathbb{R} ^{d}  } \phi(x- x')  \mathfrak{u}(x') \varrho (  x' )  \mathrm{d}  x' 
=  -\bigl(\nabla  \mathcal{V}+\nabla K*  \varrho  \bigr).$
Let $\mathfrak{A}_{\varrho}   \mathfrak{u}  :  = (\gamma  +   \phi *   \varrho)   \mathfrak{u} 
 - \int_{\mathbb{R} ^{d}  } \phi(x- x')  \mathfrak{u}(x') \varrho (  x' )  \mathrm{d}  x'  ,$
  then it follows that
$   \mathfrak{u} = - \mathfrak{A}_{\varrho}^{-1}  \bigl(\nabla  \mathcal{V}+\nabla K*  \varrho  \bigr).  $ 
 In \cite{Choi2021}, Choi considered the monokinetic density
$
\mathbf{f}^\epsilon(x,v)=\varrho^\epsilon(x)\delta_{\mathfrak{u}^\epsilon(x)}(v)
$
and a bounded Lipschitz communication weight $\phi$. For both the regular case with $\nabla K\in \mathcal{W}^{1,\infty}(\mathbb{R}^d)$ and the singular   Coulomb case, Choi rigorously derived the small-inertia limit of the hydrodynamic models \eqref{eq:pressureless-hydrodynamic-Euler-alignment-models} with $\mathcal{V}\equiv0$.

The stochastic Cucker--Smale model with interaction potential and  full velocity alignment   term  can be described  by the following many-particle system  
  \begin{equation} \label{eq:Cucker--Smale-model}
\left\{
\begin{aligned}
\mathrm{d}X_t^{\epsilon,i,N}=   {}  &   V_t^{\epsilon,i,N}   \mathrm{d} t,  \qquad i=1,\dots, N   , \\
\epsilon\mathrm{d}V_t^{\epsilon,i,N}
= {}  & \bigg[-V_t^{\epsilon,i,N}
+\frac{1}{N} \sum_{j=1}^{N}  \phi( X_t^{\epsilon,i,N}  -X_t^{\epsilon,j,N} )   (V_t^{\epsilon,j,N} -V_t^{\epsilon,i,N})    \bigg]  \mathrm{d}t  \\
    {} & - \frac{1}{N} \sum_{j=1}^{N} \nabla K(X_t^{\epsilon,i,N}   -X_t^{\epsilon,j,N}) \mathrm{d} t
+      \sigma  \mathrm{d}B_t^{i}   ,
\end{aligned}
\right.
\end{equation}
where $\phi(\cdot):   \mathbb{R}^{d}  \rightarrow  \mathbb{R}_{+}  $ and  $K(\cdot):  \mathbb{R}^{d}  \rightarrow  \mathbb{R} $  respectively represent the     communication weight function   and  interaction potential; $\epsilon > 0$ represents   small inertia/response time of individuals;    $\{B^{i}_{t}\}_{t}$  is  independent Wiener processes   with noise intensity $\sigma > 0$. For the stochastic Cucker--Smale model \eqref{eq:Cucker--Smale-model},  regarded as an SDE with state-dependent friction, the small inertia limit follows from the standard method of Hottovy et al. \cite{Hottovy2015}.
In this paper, we focus on the mean-field equation of the stochastic Cucker--Smale model \eqref{eq:Cucker--Smale-model}, namely the McKean--Vlasov SDE, which is given  
 \begin{equation} \label{eq:McKean-Vlasov-SDE}
\left\{
\begin{aligned}
\mathrm{d}X_t^\epsilon=   {}  &   V_t^\epsilon   \mathrm{d} t,\\
\epsilon\mathrm{d}V_t^\epsilon
= {}  & \left[-V_t^\epsilon
+\int_{\mathbb{R}^{2d}} \phi(X_t^\epsilon   -x')\big(y'-V_t^\epsilon   \big)
\mathfrak{q}_t^\epsilon    (\mathrm{d}x',\mathrm{d}y')\right]  \mathrm{d}t  -  \int_{\mathbb{R}^{2d}} \nabla K  ( X_t^\epsilon-x')  \mathfrak{q}_t^\epsilon    (\mathrm{d}x',\mathrm{d}y')    \mathrm{d}t  \\
&
+      \sigma  \mathrm{d}B_t   ,
\end{aligned}
\right.
\end{equation}
where $\phi(\cdot):   \mathbb{R}^{d}  \rightarrow  \mathbb{R}_{+}  $ and  $K(\cdot):  \mathbb{R}^{d}  \rightarrow  \mathbb{R} $  respectively represent the     communication weight function   and  interaction potential; $\epsilon > 0$ represents  a small inertia time, $\sigma > 0$ is the noise intensity; $\{B_t\}_{t\ge 0}$ is a Wiener process in $\mathbb{R}^d$, and $\mathfrak{q}_t^\epsilon := \mathrm{Law}(X_t^\epsilon, V_t^\epsilon)$.
Define  
$Y_t^\epsilon    :=\epsilon^{1/2} V_t^\epsilon,$
$ q_t^\epsilon      :=\mathrm{Law}(X_t^\epsilon,   Y_t^\epsilon),$
$ \rho_t^\epsilon     :=\mathrm{Law}(X_t^\epsilon)$.
 Then    \eqref{eq:McKean-Vlasov-SDE} can be rewritten as the following McKean--Vlasov
SDE
\begin{equation}\label{eq:scaled-system}
\left\{
\begin{aligned}
\mathrm{d}X_t^\epsilon  &=  \epsilon^{-1/2}  Y_t^\epsilon   \mathrm{d} t,\\
\mathrm{d}Y_t^\epsilon
&=\epsilon^{-1}\left[-Y_t^\epsilon
+\int  _{   \mathbb{R}   ^{2d}}  \phi(X_t^\epsilon   -x')\big(y'-Y_t^\epsilon   \big)
q_t^\epsilon    (\mathrm{d}x',\mathrm{d}y')\right]  \mathrm{d}t  -\epsilon^{-1/2} G_{\rho_t^\epsilon}(X_t^\epsilon  ) \mathrm{d}t
+\epsilon^{-1/2}    \sigma     \mathrm{d}B_t,
\end{aligned}
\right.
\end{equation}
where $
G_\rho  :=(\nabla K*\rho) (\cdot) = \int_{ \mathbb{R}^{d}   } \nabla K(\cdot-x')  \rho(  \mathrm{d}   x')  .
$
 By the disintegration of $q^{ \epsilon}$   with respect to  $\rho^{ \epsilon} $,  we have   $q^{ \epsilon}(\mathrm{d} x,\mathrm{d} y)=\rho^{ \epsilon}(\mathrm{d} x)q^{ \epsilon, x}(\mathrm{d} y)$  and  $ 
u^{ \epsilon }(x):=\int _{\mathbb{R}^d}  y  q^{ \epsilon, x}(\mathrm{d} y) = \mathbb{E}[Y^{ \epsilon }   \mid X^{ \epsilon }  =x]  \quad \text{for} \quad  \rho^{ \epsilon } \text{-a.e.} \ x,$
where $q^{ \epsilon, x}(\mathrm{d}y)$ is a regular conditional distribution of $Y^{ \epsilon }$ given $X^{ \epsilon }=x$.
 We then derive the hydrodynamic
Euler-alignment models    
\begin{equation*} 
\left\{
\begin{aligned}
 &\partial _{t}  \rho  ^\epsilon    +  \epsilon^{-1/2}   \nabla \cdot  (\rho  ^\epsilon    u^\epsilon    )  =0 ,\\
 &     \partial _{t}  (  \rho  ^\epsilon    u^\epsilon   )    +  \epsilon^{-1/2}    \nabla \cdot  (\rho  ^\epsilon    u^\epsilon \otimes   u^\epsilon   + \mathbf{P}^\epsilon  )        = - \epsilon^{-1  }   \rho  ^\epsilon     A_{\rho  ^\epsilon}     u^\epsilon    -  \epsilon^{-1/2}  \rho  ^\epsilon  G _{\rho  ^\epsilon}     ,
\end{aligned}
\right.
\end{equation*}
where  $  \mathbf{P}_{t}^\epsilon  (\mathrm{d} x  )  :=   \rho_{t}^\epsilon(\mathrm{d} x  )   \int_{   \mathbb{R} ^{d}  }    (  u_{t}^\epsilon(x)   -y  )   \otimes  (  u_{t}^\epsilon(x) - y )  q_{t} ^{\epsilon, x }  (      \mathrm{d} y   )$,    $A_{\rho  }$ is defined in  \eqref{eq:A-operator}.

We show that, under the hypotheses stated in Section \ref{Preliminaries-results}, the solution $X_t^\epsilon$ of \eqref{eq:scaled-system} converges to the solution $X_t$ of the limiting equation 
\begin{equation}\label{eq:limit-SDE}
 \mathrm{d} X_t=b_{\rho_t}(X_t) \mathrm{d} t+\lambda_{\rho_t}^{-1}  (X_t)  \sigma \mathrm{d} B_t,
 \qquad \rho_t= \mathrm{Law}(X_t),
\end{equation}
where   
 $b_\rho(x):=-F_\rho(x)+   \overline{\mathcal{Q}}_\rho(x)=  - A_\rho^{-1}G_\rho (x)+   \overline{\mathcal{Q}}_\rho(x)$,  
 $\lambda_\rho(x)  :=1+(\phi   *   \rho)(x)$,  with $F_\rho$ and $\overline{\mathcal{Q}}_\rho$ defined in \eqref{eq:F-series} and \eqref{eq:Qbar}, respectively.
  The noise-induced drift $\overline{\mathcal{Q}}_\rho(x)$ originates from  fast friction  and   noise, and for each $1\le \ell\le d$, its $\ell$-th component can be represented as   the unique   solution to the following contraction-type Fredholm equation
$$ 
 \overline{\mathcal{Q}} _{\rho,\ell}=-\mathfrak{H}_{\rho,\ell}
 +P_\rho   \overline{\mathcal{Q}}  _{\rho,\ell} ,
$$
where
 $\mathfrak{H}_{\rho,\ell} :=\sum_{a=1}^d\left(B_{\rho,a}R_\rho J_{\rho,\ell a}
 +M_{\rho,a}R_\rho J_{\rho,\ell a}\right)$ and 
 $ 2\lambda_\rho(x)J_\rho(x)=   \sigma ^{2}    \mathbf{I}_{d}.$
 Here, $P_{\rho}, R_{\rho}$ and $B_{\rho,a}, M_{\rho,a}$ are defined in Section \ref{global-correction}.
 Then the distribution of the limiting position $X$ no longer merely obeys the aggregation equation driven by the   force balance velocity field, but satisfies a continuity equation with a diffusive flux
 \begin{equation*} 
\left\{
\begin{aligned} 
 &\partial_t\rho  +  \nabla\cdot  \mathscr{S}_\rho=0,\\
 & \mathscr{S}_\rho
=\rho\bigl(-A_\rho^{-1}G_\rho+\overline{\mathcal{Q}}_{\rho}\bigr)
-\frac{1}{2}\nabla\cdot
\bigl(\rho  \lambda_\rho^{-2}\sigma ^{2}\mathbf{I}_{d} \bigr).
\end{aligned}
\right.
\end{equation*}

 We summarize the main contributions of this paper as follows.
 \begin{description}  
  \item[(1)]  To the best of our knowledge, this appears to be the first systematic discussion of the Smoluchowski--Kramers approximation for the distribution-dependent stochastic Cucker--Smale equation with a non-constant communication weight function and a full velocity alignment term. It may also contribute some ideas and references for subsequent work on this topic.
     In the existing literature \cite{Cerrai2020,Hottovy2015,Qian2026,Shi2024,Wang2024,Xie2022}, for SDEs with state-dependent friction or McKean--Vlasov SDEs, the additional noise-induced drift term arising in the Smoluchowski--Kramers approximation can be expressed as the solution of a Lyapunov matrix equation; whereas in our limiting equation  \eqref{eq:scaled-system}, this extra drift term is characterized as the solution of a contraction-type Fredholm equation.
        In addition, the source and the explicit expression of the additional noise-induced drift term are clearly revealed by the derivation procedure.
         \item[(2)]    In  \cite{Carrillo2021,Choi2021},  under the monokinetic ansatz
$
\mathbf{f}^\epsilon(x,v)=\varrho^\epsilon(x)\delta_{\mathfrak{u}^\epsilon(x)}(v),
$
the authors  considered both the regular and singular cases and derived the small-inertia limit of the hydrodynamic Euler-alignment models, yielding the corresponding nonlocal aggregation equations. 
        In contrast, we study the small-inertia limit of a stochastically forced Cucker--Smale equation directly at the McKean--Vlasov level. The random forcing models unresolved fluctuations and, crucially, prevents a monokinetic closure of the velocity distribution. 
        The resulting small-inertia limit  is a distribution-dependent diffusion with, in general, an additional drift generated by the interplay between noise and spatially varying alignment.
        In addition, the limiting position distribution no longer merely satisfies an aggregation equation driven by the force balance velocity field; instead, it satisfies a continuity equation with a diffusive flux.
         This mechanism is absent from the  deterministic models.  This paper is complementary to the deterministic particle  limit work of   \cite{Carrillo2021,Choi2021}. 
        
      \item[(3)]  In contrast to the works \cite{Shi2024,Wang2022,Wang-Lv-Zhang-Duan-Wang2024}, equation \eqref{eq:scaled-system} features a non-constant communication weight function and a full velocity alignment term, which are the main differences and also the main difficulties. For this reason, the methods in \cite{Shi2024,Wang2022,Wang-Lv-Zhang-Duan-Wang2024} are not applicable to   equation  \eqref{eq:scaled-system}. In fact,  equation \eqref{eq:scaled-system}  can be viewed as a SDE with distribution-dependent and state-dependent friction, and it also involves the co-evolution of the source term with the velocity distribution. 
        Inspired by the idea of   Hottovy et al. \cite{Hottovy2015}, we also  rearrange   the equation \eqref{eq:scaled-system}   by  integration by parts.  However, the co-evolution of the source term with the velocity distribution gives rise to a singular scaling $\epsilon^{-1/2}$ in the velocity equation. Subsequently, we repeatedly perform integration by parts on independent copies of the equation, which leads to an infinite singular chain. This naturally motivates the construction of both a truncated position correction $Z_{t} ^{(N)}$ and a truncated quadratic functional $\mathcal{Q}^{(N)}$. Passing to the infinite limit $N \rightarrow  \infty$, we then construct the corresponding global position correction $Z_{t}  $ and the global quadratic functional $\mathcal{Q}$, and in doing so an SDE for the corrected position \eqref{eq:Z-decomp}, in which the singular scaling $\epsilon^{-1/2}$ is eliminated.  We observe that in previous  literature \cite{Hottovy2015,Shi2024,Wang2024}  the fast friction is state-dependent friction, whereas the full fast friction in equation \eqref{eq:scaled-system} is a nonlocal operator $A_\rho$ that depends on the spatial marginal of the distribution.

           \item[(4)] Since the increment estimate $X_t ^\epsilon   -X_s^\epsilon=\epsilon^{-1/2} \int_s^t Y_r^\epsilon  \mathrm{d} r$ contains the singular scaling $\epsilon^{-1/2}$, it is difficult to directly establish path tightness uniformly in $\epsilon$. However, the equation for the corrected position \eqref{eq:Z-decomp} has already eliminated the singular scaling $\epsilon^{-1/2}$, so we first aim to obtain the tightness of $\{Z^\epsilon\}_{0< \epsilon \le1}$. Indeed, using $X_t^\epsilon=X_0^\epsilon+ \epsilon^{-1/2}  \int_0^t Y_s^\epsilon  \mathrm{d }s$ together with the velocity equation \eqref{eq:scaled-system} and differential inequality, we derive uniform moment estimates for the velocity $Y_t^\epsilon$ and the corrected position $Z_{t}^\epsilon$. Then, by applying Aldous criterion, we establish the tightness of $\{Z^\epsilon\}_{0< \epsilon \le1}$, and subsequently obtain the tightness of $\{X^\epsilon\}_{0< \epsilon \le1}$ as well as an estimate for the second moment of the increments. 
                  For the Smoluchowski--Kramers approximation, we adopt the averaging method and freeze the slow variable $X$ and its distribution on the fast time scale. Since the velocity still varies rapidly after freezing, its phase-space law evolves with the fast time; thus we provide a frozen equation that retains the full source term. Subsequently, it is shown to be uniformly exponentially mixing. The quadratic functional $\mathcal{Q}$ is decomposed into the equilibrium average $\overline{\mathcal{Q}}$, local quadratic fluctuations, and exponential transients. Then the Poisson--Lyapunov equation is used to establish an averaging lemma, namely the convergence of $\mathcal{Q}$ to $\overline{\mathcal{Q}}$ in the time integral. Finally, under suitable assumptions, we obtain weak convergence and the convergence rate of strong $L^1$-convergence.
  \end{description}

 	The remainder of this paper is organized as follows.
 In Section \ref{Preliminaries-results}, we introduce the necessary notation and hypotheses, and state our main results, Theorems \ref{thm:sequence-identification} and \ref{thm:strong-convergence}. In Section \ref{correction-tightness}, we prove the tightness of the positions. In Section \ref{Smoluchowski-Kramers-approximation}, we provide the proof of the theorems.
 Throughout this paper, we use $C$, possibly with subscripts, to denote positive constants that may vary from line to line, with $C_p$ indicating dependence on $p$.

\section{Preliminaries and main results}\label{Preliminaries-results}

Let $(\Omega, \mathcal{F}, \mathbb{P})$ be a complete probability space, and denote by $\mathbb{E}$ the expectation with respect to $\mathbb{P}$. 
For a vector $x\in\mathbb{R}^d$, let $|x|$ be the Euclidean norm; for a measurable function $h:\mathbb{R}^d  \to  \mathbb{R}$, define $\|h\|_\infty:=\sup_{x\in\mathbb{R}^d}|h(x)|$.
Given a probability measure $\mu $ on $\mathbb{R}^d$, let $L^p(\mu )$ denote the usual space of $p$-integrable functions with respect to $\mu $. For a linear operator $P: L^p(\mu ) \to L^p(\mu )$, we denote its operator norm by
$$
\|P\|_{L^p(\mu ) \to L^p(\mu )}
:= \sup_{\|h\|_{L^p(\mu )} \le 1} \|P h\|_{L^p(\mu )}.
$$
We assume that the initial data of equation \eqref{eq:scaled-system} satisfy 
$ \sup  \limits_{0<\epsilon  \le  1}  \mathbb{E}\bigl[|X_0^\epsilon|^4+|Y_0^\epsilon |^4\bigr]<\infty.$

For $p>0$, let $\mathcal{P}_p(\mathbb{R}^d)$ denote the space of Borel probability measures on $\mathbb{R}^d$ with finite $p$-th moment:
$$
\mathcal{P}_p(\mathbb{R}^d):=\left\{\mu\in\mathcal{P}(\mathbb{R}^d):\int_{\mathbb{R}^d}|x|^p\,\mu(\mathrm{d}x)<\infty\right\},
$$
where $\mathcal{P}(\mathbb{R}^d)$ is the set of all Borel probability measures on $\mathbb{R}^d$. For $\mu,\nu\in\mathcal{P}_p(\mathbb{R}^d)$, the $p$-th Wasserstein distance is defined by
$$
\mathbb{W}_p(\mu,\nu):=\left(\inf_{\pi\in \mathcal{T}(\mu,\nu)}\int_{\mathbb{R}^d\times\mathbb{R}^d}|x-y|^p\,\pi(\mathrm{d}x,\mathrm{d}y)\right)^{1/p},
$$
where $\mathcal{T}(\mu,\nu)$ denotes the collection of all couplings of $\mu$ and $\nu$, i.e., joint probability measures on $\mathbb{R}^d\times\mathbb{R}^d$ with marginals $\mu$ and $\nu$.

In what follows we  assume the following conditions.
\begin{description}
  \item[(A1)]  $\phi  \in \mathcal{C}_{b}^{1}(\mathbb{R}^d)$ is even and nonnegative. Set $\phi_M   :=\|\phi\|_\infty   $ and $L_\phi     :=\|\nabla\phi\|_\infty$.
   \item[(A2)]  $\nabla K( \cdot ) $ is Lipschitz continuous with Lipschitz constant  $L_{K}$. 
     \item[{(A1\textquotesingle)}]   $\phi  \in \mathcal{C}_{b}^{1,1}(\mathbb{R}^d)$ is even and nonnegative. Set $\phi_M :=\|\phi\|_\infty$ and $L_\phi   :=\|\nabla\phi\|_\infty$. 
   \item[{(A2\textquotesingle)}]  $\nabla K( \cdot ) $ is bounded and Lipschitz continuous with Lipschitz constant  $L_{K}$. 
\end{description}

The two main theorems of this paper are as follows.

 \begin{theorem}\label{thm:sequence-identification}
  Assume that \textup{\textbf{(A1)}-\textbf{(A2)}} hold, that $\rho_0^\epsilon \rightarrow \rho_0$, that $\{|X_0^\epsilon|\}_{\epsilon>0}$ is uniformly integrable, and that the McKean--Vlasov SDE \eqref{eq:limit-SDE} has a unique weak solution. Denote its distribution by $\mathfrak{P}$.  Let $X=\{X_t\}_{t\in[0,T]}$ be the canonical coordinate process on $\mathcal{C}([0,T];\mathbb{R}^d)$ endowed with $\mathfrak{P}$, and set
$
\rho_t := \mathfrak{P}\circ X_t^{-1},$  $ t\in[0,T].
$
Then the whole family $\{\mathrm{Law}(X^\epsilon)\}$ converges weakly in $\mathcal{C}([0,T];\mathbb{R}^d)$ to $\mathfrak{P}$:
$
\mathrm{Law}(X^\varepsilon)  \Longrightarrow  \mathfrak{P}.
$
Consequently, the weak limit of the sequence $X^\epsilon$ is the unique weak solution of the McKean--Vlasov SDE \eqref{eq:limit-SDE}.  In particular, for every $t\in[0,T]$,
$
\rho_t^\epsilon := \mathrm{Law}(X_t^\varepsilon) \Longrightarrow \rho_t.
$
\end{theorem}

\begin{theorem}  \label{thm:strong-convergence}
Assume \textup{\textbf{(A1\textquotesingle)}-\textbf{(A2\textquotesingle)}}   and   
$\mathrm{ Law}(X_0)=\rho_0.
$
Then 
\begin{equation}
\mathbb{E}\Big[\sup_{0\le t\le T} |X_t^\epsilon-X_t|\Big]  \leq   C_{T} \Big( \Big(  \mathbb{E} \Big[  | X _0^\epsilon -X_0  |^{2}  \Big]   \Big)^{1/2}  + \epsilon^{1/5}\Big).
\end{equation}
\end{theorem}

\begin{remark}
  If the initial values satisfy
$
\sup  \limits_{0<\epsilon\le1} \mathbb{E}\Big[|X_0^\epsilon|^8+|Y_0^\epsilon|^8\Big]<\infty$, 
then, under Assumptions \textup{\textbf{(A1\textquotesingle)}-\textbf{(A2\textquotesingle)}}, using Remark \ref{rem:2p-moment-closure}, we obtain
\begin{equation*}
\mathbb{E}\Big[\sup_{0\le t\le T} |X_t^\epsilon-X_t|  \Big] \leq   C_{T} \Big(   \Big(  \mathbb{E} \Big[  | X _0^\epsilon -X_0  |^{2}  \Big]   \Big)^{1/2}   + \epsilon^{1/4}\Big).
\end{equation*}
\end{remark}

\section{  The global   correction and  tightness   }\label{correction-tightness}

  In this section, we prove the tightness of $\{X^\epsilon\}_{0<\epsilon\le 1}$.
  In Section \ref{global-correction}, we introduce the global position correction and the quadratic remainder. In Section \ref{origin-correction-remainder}, we start from the model without velocity alignment and then add the full Cucker--Smale source terms, which naturally yields the truncated position correction and the quadratic remainder, thereby revealing the origin of the global correction and the remainder term. We also derive the SDE for the truncated corrected position. In Section \ref{estimates-truncation-limits}, we   obtain the SDE satisfied by the limiting corrected position. In Section \ref{estimates-tightness}, we establish some uniform moment estimates, which in turn yield tightness of $\{X^\epsilon\}_{0<\epsilon\le 1}$ and increment estimates.

\subsection{     The global correction vector }\label{global-correction}

Introduce the following notation
$$
\lambda_\rho(x)  :=1+(\phi   *   \rho)(x),
\qquad
(L_\rho h)(x):=\lambda_\rho(x)h(x),
$$
and
$$
(K_\rho h)(x):=\int _{\mathbb{R}   ^{d}}   \phi(x-x')h(x')\rho( \mathrm{d} x').
$$
Define the fibered operator
\begin{align}\label{eq:A-operator}
(A_\rho h)(x) := {}  &(L_\rho  h-K_\rho  h)(x)   \notag \\
 ={}& \lambda_\rho(x) h(x) - (K_\rho h)(x)    \notag \\
 ={}  & h(x) + (\phi * \rho)(x) h(x) - \int_{\mathbb{R}^d} \phi(x-x') h(x') \rho(\mathrm{d}  x').
\end{align}
The operator factorises as $A_\rho = L_\rho (I - P_\rho)$ with
\begin{equation*}
P_\rho := \lambda_\rho^{-1} K_\rho   =R_\rho  K_\rho , \qquad
(P_\rho h)(x) = \frac{1}{\lambda_\rho(x)} \int_{\mathbb{R}^d} \phi(x-x') h(x') \rho(\mathrm{d}   x'),
\end{equation*}
where  $R_\rho  : =L_\rho ^{-1}.$
For each coordinate $a=1,\dots,d$, define 
\begin{align}
 (B_{\rho,a}h)(x)
 &:=\frac1{\lambda_\rho(x)}\int    _{  \mathbb{R}  ^{ d}}    \partial_a\phi(x-x')h(x')\rho( \mathrm{d} x'),\label{eq:B-def}\\
 (M_{\rho,a}h)(x)
 &:=\frac{\partial_a\lambda_\rho(x)}{\lambda_\rho(x)}h(x).
 \label{eq:M-def}
\end{align}

\begin{lemma}\label{lem:contraction}
Set $\mu_\rho(\mathrm{d}  x) := \lambda_\rho(x) \rho(\mathrm{d}  x)$. Under \textup{\textbf{(A1)}-\textbf{(A2)}},  for any $ 1\leq   p \leq   \infty$,
$$
\|P_\rho\|_{L^p(\mu_\rho) \to L^p(\mu_\rho)}
\le \frac{\phi_M}{1 + \phi_M} =: r_{*} < 1,
$$
where $p<\infty$ refers to the $L^p(\mu_\rho)$ norm, and $p=\infty$ refers to the $L^\infty(\mu_\rho)$ norm, i.e., the essential supremum with respect to $\mu_\rho$, denoted by $\|h\|_{L^\infty(\mu_\rho)} := \operatorname{ess\,sup}  \limits _{x\sim \mu_\rho} |h(x)|$.
\end{lemma}

\begin{proof}
For $\mu_\rho$, the integral kernel of $P_\rho$ is
$$
p_\rho(x,x')=\frac{\phi(x-x')}{\lambda_\rho(x)\lambda_\rho(x')}.
$$
  By using the evenness of $\phi$,    its row mass satisfies
$$
\int_{\mathbb{R}^d}  p_\rho(x,x')\,\mu_\rho(\mathrm{d} x')
=\frac{(\phi*\rho)(x)}{1+(\phi*\rho)(x)}\le \frac{\phi_M}{1 + \phi_M} = r_{*},
$$
where we have used the monotonicity of $  \frac{x}{1+x}$ and the fact that $(\phi*\rho)(x)\le \phi_M$.
Due  to  symmetry,   we have  
$
\int_{\mathbb{R}^d}   p_\rho(x,x')\,\mu_\rho(\mathrm{d}  x)\le r_{*}.
$
Then  for $h\in L^\infty(\mu_\rho)$,
$$
|P_\rho h(x)|\le r_*\|h\|_{L^\infty(\mu_\rho)},
$$
and for $h\in L^{1}(\mu_\rho)$, Tonelli's theorem together with the column mass bound gives
$$
\|P_\rho h\|_{L^{1}(\mu_\rho)}
\le \int_{\mathbb{R}^d} |h(x')|\left(\int_{\mathbb{R}^d} p_\rho(x,x')\,\mu_\rho(\mathrm{d}  x)\right)\mu_\rho(\mathrm{d}  x')
\le r_*\|h\|_{L^{1}(\mu_\rho)}.
$$
Consequently, by Riesz-Thorin theorem, for every $p\in[1,\infty]$,  we obtain 
$$
\|P_\rho\|_{L^p(\mu_\rho)\to L^p(\mu_\rho)}
\le \frac{\phi_M}{1+\phi_M} =: r_{*} < 1.
$$
\end{proof}

\begin{remark}
In fact, from Lemma \ref{lem:contraction}, we obtain  $\sup_\rho\|P_\rho\|_{L^p(\mu_\rho)\to L^p(\mu_\rho)} \le r_*<1.$  In addition,
for any globally bounded Borel function $h$, we have
$     \sup   \limits   _{ x \in   \mathbb{R}^d}  |P_\rho h(x)|  \le r_*\|h\|_{ \infty } .$
\end{remark}

 \begin{remark}
For each $\epsilon$, the auxiliary space $L^r(\mu_{\rho^\epsilon})$ varies. However, norm estimates in the auxiliary spaces can be converted into estimates for random variables on a common probability space. We can   obtain $\epsilon$-uniform moment estimates, and hence this does not affect our discussion of weak convergence and strong convergence.
\end{remark}

Consequently, $A_\rho^{-1}$ admits the convergent Neumann series
\begin{equation}
A_\rho^{-1} =   (I - P_\rho)^{-1} R_\rho 
=   \sum_{n=0}^{\infty} P_\rho^n  R_\rho ,
\end{equation}
with the operator norm bound $\|A_\rho^{-1}\|_{L^p(\mu_\rho) \to L^p(\mu_\rho)}
\le \frac{1}{   1-r_{*} }$.
Define
\begin{equation}\label{eq:weights}
w_{0,\rho}(x_0)=\lambda_\rho^{-1}  (x_0),
\qquad
w_{n,\rho}(x_0,\ldots,x_n)
=\frac{\prod_{j=0}^{n-1}\phi(x_j-x_{j+1})}
{\prod_{j=0}^{n}\lambda_\rho(x_j)},\quad n\ge1,
\end{equation}
and
$$
g_{n,i,\rho}(x_0,\ldots,x_n):=\nabla_{x_i}w_{n,\rho}(x_0,\ldots,x_n) \quad   \text{ for }  i=0, \dots, n.
$$
Let  $q(\mathrm{d} x,\mathrm{d} y)=\rho(\mathrm{d} x)q^x(\mathrm{d} y),$  
where $q^x(\mathrm{d}y)$ is the conditional distribution of $Y$ given $X=x$.
Set
$$
u(x):=\int _{\mathbb{R}^d}  y  q^x(\mathrm{d} y) = \mathbb{E}[Y   \mid X=x],  \quad  \text{for} \quad  \rho\text{-a.e.} \ x,
$$
$$
\beta_q(x)  :=\int_{\mathbb{R}^{2d}}    \nabla\phi(x-x')\cdot y'\,q(\mathrm{d} x',\mathrm{d} y').
$$
 Next, we define the correction vector field and the correction remainder,    which  are  both   well  defined  by Lemma \ref{lem:kernel-sum}. 
Define the correction vector field
\begin{align*}
 \Xi(q;x_0,y_0)
 :={}&\sum_{n=0}^\infty\int _{  \mathbb{R}  ^{2nd}}
 w_{n,\rho}(x_0,\ldots,x_n)y_n\,
 q^{\otimes n}(\mathrm{d} z_1,\dots  ,\mathrm{d} z_n)\\
={}&   R_\rho(x_{0})y_{0}
 + \sum_{n=1}^\infty   (P_\rho^n    R_\rho u)(x_{0})\\
 ={}& R_\rho(x_0)\bigl(y_0-u (x_0)\bigr)+  (A_\rho^{-1}u )(x_0),
\end{align*}
where $z_i=(x_i,y_i)$, with the convention that the term $z_0=(x_0,y_0)$ is not integrated, and the $n=0$ term is interpreted as $w_{0,\rho}(x_0)y_0$. We slightly abuse notation by writing $R_\rho(x)$ for the value of the multiplication operator with symbol $\lambda_\rho^{-1}$ at $x$.
The correction remainder is defined as
\begin{equation}\label{eq:Q-def}
\begin{aligned}
\mathcal{Q}(q;x_0,y_0)
=\sum_{n=0}^\infty\int  _{  \mathbb{R}  ^{2nd}}    y_n\Bigg[&
\sum_{i=0}^n y_i\cdot g_{n,i,\rho}(x_0,\ldots,x_n)      \\
&+w_{n,\rho}(x_0,\ldots,x_n)
\sum_{i=0}^n\frac{\beta_q(x_i)}{\lambda_\rho(x_i)}
\Bigg]q^{\otimes n}( \mathrm{d} z_1,  \dots ,  \mathrm{d} z_n),
\end{aligned}
\end{equation}
where $z_i=(x_i,y_i)$, with the convention that the term $z_0=(x_0,y_0)$ is not integrated, and the $n=0$ term is interpreted as $y_0\bigl[y_0\cdot \nabla R_\rho(x_0)+  R_\rho(x_0) \frac{\beta_q(x_0)}{\lambda_\rho(x_0)}\bigr]$.

\subsection{The origin of the global correction and the remainder term } \label{origin-correction-remainder}
Now we start with the model without velocity alignment and then add the full Cucker--Smale source terms. In this way, one can directly see that the position corrections and the Neumann series are not guessed, but rather are successively forced when eliminating the fast velocity terms that have not yet been canceled.

\noindent {\textbf{  Case 1:  absence of velocity alignment and constant friction.} }
If $\phi \equiv 0$, then $\lambda \equiv 1$, and the scaled system reduces to
$$
 \mathrm{d} X_t^\epsilon = \epsilon^{-1/2} Y_t^\epsilon\,\mathrm{d} t,
 \qquad
 \mathrm{d} Y_t ^\epsilon= -\epsilon^{-1} Y_t^\epsilon\,\mathrm{d} t
 -  \epsilon^{-1/2}  G_{\rho_t^\epsilon}(X_t^\epsilon)\,\mathrm{d} t
 +\epsilon^{-1/2} \sigma\,\mathrm{d} B_t.
$$
Let the corrected position be $Z_t^\epsilon := X_t^\epsilon + \epsilon^{ 1/2}   Y_t^\epsilon$. Then
\begin{equation}\label{eq:no-alignment-corrected-position}
 \mathrm{d} Z_t^\epsilon = -G_{\rho_t^\epsilon}(X_t^\epsilon)\,\mathrm{d} t + \sigma\,\mathrm{d} B_t.
\end{equation}
The equation \eqref{eq:no-alignment-corrected-position}  does not leave any new $\epsilon^{-1/2}$ singular terms.

\noindent {\textbf{  Case 2:  ordinary position-dependent friction. } }
Consider temporarily the model without the nonlocal velocity source term but with scalar friction $\gamma(x)\ge \gamma_0>0$, and denote $R_\gamma := \gamma^{-1}$. The scaled system reduces to
\begin{equation*} 
\begin{aligned}
\mathrm{d} Y_t^\epsilon
&= -\epsilon^{-1}\gamma(x)Y_t^\epsilon
 \,\mathrm{d} t -\epsilon^{-1/2}G_{\rho_t^\epsilon}(X_t^\epsilon)\,\mathrm{d} t
+\epsilon^{-1/2} \sigma\,\mathrm{d} B_t.
\end{aligned}
\end{equation*}
Multiplying both sides of the above equation by $\epsilon^{ 1/2}  R_\gamma$, we have
\begin{align}
\epsilon^{ 1/2}  \,\mathrm{d}    \bigl(R_\gamma(X_t^\epsilon)Y_t^\epsilon\bigr)
={}&-   \epsilon^{-1/2}Y_t^\epsilon   \mathrm{d}  t
 -R_\gamma(X_t^\epsilon)G_{\rho_t^\epsilon}(X_t^\epsilon)\mathrm{d}  t
 +R_\gamma(X_t^\epsilon)\sigma\mathrm{d}    B_t\notag\\
&+Y_t^\epsilon   (Y_t^\epsilon)^{\top}\nabla R_\gamma(X_t^\epsilon)\mathrm{d}   t,
\label{eq:ordinary-friction-corrector}
\end{align}
where the last term comes from $\epsilon^{ 1/2}    \nabla R_\gamma(X_t^\epsilon  )Y_t^\epsilon  \cdot   \mathrm{d} X_t^\epsilon$ and represents the noise-induced drift under ordinary position-dependent friction.
Let the corrected position be 
$
 Z_t^\epsilon :=X_t^\epsilon+  \epsilon^{ 1/2}  R_\gamma(X_t^\epsilon)Y_t^\epsilon.
$
 Then
\begin{align*}
 \mathrm{d} Z_t^\epsilon
= 
 -R_\gamma(X_t^\epsilon)G_{\rho_t^\epsilon}(X_t^\epsilon)  \mathrm{d} t
 +R_\gamma(X_t^\epsilon)\sigma  \mathrm{d} B_t       +  Y_t^\epsilon   (Y_t   ^\epsilon)^{\top}\nabla R_\gamma(X_t^\epsilon)   \mathrm{d} t.
\end{align*}

\noindent {\textbf{  Case 3:  full Cucker--Smale velocity alignment term included.} }
  Denote the observed copy by $(X_t^{\epsilon,0},Y_t^{\epsilon,0},B_t^{\epsilon,0})=(X_t^{ \epsilon},Y_t^{ \epsilon},B_t^{ \epsilon})$.
The velocity equation  of   \eqref{eq:scaled-system}  is
$$
 \mathrm{d}Y_t^{\epsilon,0}
 =\epsilon^{-1 }\bigl[-\lambda_{\rho_t^\epsilon}(X_t^{\epsilon,0})Y_t^{\epsilon,0}
 +(K_{\rho_t^\epsilon} u_t^\epsilon)(X_t^{\epsilon,0})\bigr]\mathrm{d} t
 -   \epsilon^{-1/2}  G_{\rho_t^\epsilon}(X_t^{\epsilon,0})\mathrm{d} t+   \epsilon^{-1/2}\sigma  \mathrm{d} B_t^{\epsilon,0}.
$$
Multiplying both sides of the above equation by $ \epsilon^{ 1/2} R_{\rho_t ^\epsilon}(X_t^{\epsilon,0})$, we have
\begin{equation}\label{eq:first-CS-rearrangement-actual}
\begin{aligned}
 \epsilon^{-1/2}  Y_t^{\epsilon,0}\mathrm{d} t
={}&-\epsilon^{ 1/2} R_{\rho_t^\epsilon}(X_t^{\epsilon,0})\mathrm{d} Y^{\epsilon,0}
 -R_{\rho_t^\epsilon}(X_t^{\epsilon,0})G_{\rho_t^\epsilon}(X_t^{\epsilon,0})\mathrm{d} t
 +R_{\rho_t^\epsilon}(X_t^{\epsilon,0})\sigma   \mathrm{d} B^{\epsilon,0}\\
 &+   \epsilon^{-1/2}(P_{\rho_t^\epsilon} u_t^\epsilon)(X_t^{\epsilon,0})\mathrm{d} t.
\end{aligned}
\end{equation}
By the product rule,  
$$
 -  \epsilon^{1/2}   R_{\rho_t^\epsilon}(X_t^{\epsilon,0})\mathrm{d} Y_t^{\epsilon,0}
 =-   \mathrm{d}   \bigl[\epsilon^{1/2}   R_{\rho_t ^ \epsilon}(X_t^{\epsilon,0})Y_t^{\epsilon,0}\bigr]
 +\epsilon^{1/2}  Y_t^{\epsilon,0}    \mathrm{d}  \bigl[R_{\rho_t^\epsilon}(X_t^{\epsilon,0})\bigr].
$$
From $\mathrm{d} X_t^{\epsilon,0} = \epsilon^{-1/2}   Y_t^{\epsilon,0}\,\mathrm{d} t$, for any $\varphi\in \mathcal{C}_c^1( \mathbb{R}^d)$,
$$
\frac{\mathrm{d}}{\mathrm{d} t}\int_{   \mathbb{R}^{ d}   }   \varphi(x)\rho_t^{\epsilon}(\mathrm{d} x)
=\epsilon^{-1/2}\int_{   \mathbb{R}^{2d}   }      \nabla\varphi(x)\cdot y\,q_t^{\epsilon}(\mathrm{d} x,\mathrm{d} y)
=\epsilon^{-1/2}   \int_{   \mathbb{R}^{ d}   }    \nabla\varphi(x)\cdot u_t^{\epsilon}(x)\rho_t^{\epsilon}(\mathrm{d} x).
$$
Hence the spatial marginal satisfies the continuity equation in the distributional sense 
$$
\partial_t\rho_t^{ \epsilon}+  \epsilon^{-1/2}   \nabla  \cdot  (\rho_t^{ \epsilon} u_t^{ \epsilon})=0.
$$
 By   same  discussion  for   McKean--Vlasov  SPDEs   \cite{Hong2025},   system \eqref{eq:scaled-system} has a unique strong solution,    and for each fixed $ \epsilon >0$ and $T>0$,
\begin{equation}\label{eq:fixed-eps-moments}
 \mathbb{E}\Bigl[\sup_{t\le T}\bigl(|X_t^{ \epsilon}|^4+|Y_t^\epsilon|^4\bigr)\Bigr]
 \le C_{\epsilon,T}<\infty.
\end{equation}
Taking $\varphi(x')=\phi(x-x')$, first approximate by compactly supported truncations, and then use the boundedness of $\phi,\nabla \phi$ and the first-order velocity moment together with dominated convergence to remove the truncation, we obtain
\begin{equation}\label{eq:lambda-time}
 \partial_t\lambda_{\rho_t^{ \epsilon}}(x)
 =-   \epsilon^{-1/2}  \int_{   \mathbb{R}^{ d}   }   \nabla   \phi(x-x')\cdot u_t^{ \epsilon}(x')\rho_t^{ \epsilon}(\mathrm{d} x')
 =- \epsilon^{-1/2}  \beta_{q_t^{ \epsilon}}(x).
\end{equation}
Then
\begin{align*}
 \mathrm{d}   \bigl[R_{\rho_t^{ \epsilon}}(X_t^{\epsilon,0})\bigr]
 &=\partial_tR_{\rho_t^{ \epsilon}}(X_t^{\epsilon,0})\mathrm{d}    t+
   \nabla R_{\rho_t^{ \epsilon}}(X_t^{\epsilon,0})\cdot   \mathrm{d} X_t^{\epsilon,0} \notag\\
 &=\epsilon^{-1/2}  R_{\rho_t^{ \epsilon}}(X_t^{\epsilon,0})^2
 \bigl[\beta_{q_t^{  \epsilon}}(X_t^{\epsilon,0})-\nabla\lambda_{\rho_t^{  \epsilon}}(X_t^{\epsilon,0})\cdot Y_t^{\epsilon,0}\bigr]  \mathrm{d} t.
\end{align*}
Let  $Z_t^{ \epsilon, (0)  }:=X_t^{\epsilon,0}+   \epsilon^{  1/2}  R_{\rho_t^{ \epsilon}}(X_t^{\epsilon,0})Y_t^{\epsilon,0}$     
\begin{equation}\label{eq:local-CS-corrected-position}
\mathrm{d} Z_t^{\epsilon, (0)  }
 =-(R_{\rho_t ^{ \epsilon}} G_{\rho_t  ^{ \epsilon}})(X_t^{\epsilon,0})\mathrm{d}   t+ \mathcal{Q}_{ 0 } (q_t^{ \epsilon};X_t^{\epsilon,0} ,Y_t^{\epsilon,0})  \mathrm{d} t
 +R_{\rho_t^{ \epsilon}}(X_t^{\epsilon,0})\sigma\mathrm{d}   B_t^{\epsilon,0}
 + \epsilon^{-1/2}  (P_{\rho_t   ^{ \epsilon} } u_t^{ \epsilon})(X_t^{\epsilon,0})   \mathrm{d} t,
\end{equation}
where $\mathcal{Q}_{0}$  is the $0$-th layer of    $\mathcal{Q}$, denoted as
$$
  \mathcal{Q}_{0}(q;x,y)
 :=y\left[y\cdot\nabla R_\rho(x) +R_\rho(x)\frac{\beta_q(x)}{\lambda_\rho(x)}\right],
 \qquad
 \beta_q(x)=\int_{   \mathbb{R}^{2d}   }     \nabla\phi(x-x')\cdot y'\,q(\mathrm{d} x',\mathrm{d} y').
$$
The difference between equation \eqref{eq:local-CS-corrected-position} and the ordinary state-dependent equation is precisely the last singular term.

Next we handle the singular term $\epsilon^{-1/2} (P_{\rho_t ^{ \epsilon}} u_t^{ \epsilon})(X_t^{\epsilon,0})$. Take a copy $(X^{\epsilon,1},Y^{\epsilon,1},B^{\epsilon,1})$ that is independent and identically distributed as the zeroth copy, satisfying
$(X^{\epsilon,1},Y^{\epsilon,1},B^{\epsilon,1})\stackrel{\mathrm{d}}= (X^{\epsilon,0},Y^{\epsilon,0},B^{\epsilon,0})$. We obtain
\begin{equation}\label{eq:Pu-one-copy}
 (P_{\rho_t^{ \epsilon}} u_t^{ \epsilon})(X_t^ {  \epsilon,0})
 =  \mathbb{E}_1      \left[R_{\rho_t^{ \epsilon}}(X_t^{ \epsilon, 0 }   )\phi(X_t^{\epsilon,0}  -X_t^{\epsilon,1} )Y_t^{\epsilon,1} \right],
\end{equation}
where \( \mathbb{E}_1\) denotes expectation with respect to the first copy only, so the right-hand side is still an adapted quantity with respect to the filtration of the zeroth copy.
Repeat \eqref{eq:first-CS-rearrangement-actual} for $Y_t^{\epsilon,1}$
\begin{equation}\label{eq:copy-one-rearrangement}
\begin{aligned}
 \epsilon^{-1/2}   Y_t^{\epsilon,1}  \mathrm{d} t
={}&-  \epsilon^{ 1/2}  R_{\rho_t^{   \epsilon}}(X_t^{\epsilon,1})\mathrm{d} Y_t^{\epsilon,1}
 -R_{\rho_{ t}^{  \epsilon}}(X_t^{\epsilon,1})G_{\rho_t^{  \epsilon}}(X_t^{\epsilon,1})\mathrm{d} t
 +R_{\rho_t^{   \epsilon}}(X_t^{\epsilon,1})\sigma  \mathrm{d} B^{\epsilon,1}\\
 &+  \epsilon^{-1/2}  (P_{\rho_t^{   \epsilon}} u_t^{   \epsilon})(X_t^{\epsilon,1})\mathrm{d} t.
\end{aligned}
\end{equation}
Substituting \eqref{eq:copy-one-rearrangement} into \eqref{eq:Pu-one-copy}, we obtain
\begin{align*} 
 & \epsilon^{-1/2}(P_{\rho_t^{   \epsilon}} u_t^{   \epsilon})(X_t^{\epsilon,0} )  \mathrm{d} t  \\
 = {}& \epsilon^{-1/2}   \mathbb{E}_1      \left[R_{\rho_t^{   \epsilon}}(X_t^{\epsilon,0} )\phi(X_t^{\epsilon,0} -X_t^{\epsilon,1} )Y_t^{\epsilon,1} \right]  \mathrm{d} t\\
 = {} &-  \epsilon^{ 1/2} \mathbb{E}_1      \left[R_{\rho_t^{   \epsilon}}(X_t^{\epsilon,0} )\phi(X_t^{\epsilon,0} -X_t^{\epsilon,1} ) R_{\rho_t^{   \epsilon}}(X_t^{\epsilon,1} )\mathrm{d} Y_t^{\epsilon,1} \right]
 -   \mathbb{E}_1      \left[R_{\rho_t^{   \epsilon}}(X_t^{\epsilon,0} )\phi(X_t^{\epsilon,0} -X_t^{\epsilon,1} )R_{\rho_{ t}^{   \epsilon}}(X_t^{\epsilon,1} )G_{\rho_t^{   \epsilon}}(X_t^{\epsilon,1} )  \right]\mathrm{d} t  \\
 &  + \epsilon^{-1/2}   \mathbb{E}_1      \left[R_{\rho_t^{   \epsilon}}(X_t^{\epsilon,0} )\phi(X_t^{\epsilon,0} -X_t^{\epsilon,1} )(P_{\rho_t^{   \epsilon}} u_t^{   \epsilon})(X_t^{\epsilon,1} )  \right]\mathrm{d} t. 
\end{align*}
Set
$$
 w_{1,\rho_t^{   \epsilon}}(X_t^{\epsilon,0} ,X_t^{\epsilon,1} )
 :=R_{\rho_t^{   \epsilon}}(X_t^{\epsilon,0} )\phi(X_t^{\epsilon,0} -X_t^{\epsilon,1} )R_{\rho_t ^{   \epsilon}}(X_t^{\epsilon,1} ).
$$
Using  the product rule, we get
\begin{align*}
 -   \epsilon^{ 1/2}   \mathbb{E}_1   \left[
 w_{1,\rho_t^{   \epsilon}}(X_t^{\epsilon,0} ,X_t^{\epsilon,1} )\mathrm{d} Y_t^{\epsilon,1} \right]
={}&-\mathrm{d}  \left\{  \epsilon^{ 1/2}   \mathbb{E}_1   \left[
 w_{1,\rho_t^{   \epsilon}}(X_t^{\epsilon,0} ,X_t^{\epsilon,1} )Y_t^{\epsilon,1} \right]\right\}\\
 &+   \epsilon^{ 1/2}   \mathbb{E}_1    \left[Y_t^{\epsilon,1}   \mathrm{d}
 w_{1,\rho_t ^{   \epsilon}}(X_t^{\epsilon,0} ,X_t^{\epsilon,1})\right].
\end{align*}
From $\mathrm{d} X_t^{\epsilon,1}=  \epsilon^{-1/2}   Y_t^{\epsilon,1}  \mathrm{d} t$ and 
$\partial_t\lambda_{\rho_t^{   \epsilon}}(x)=-   \epsilon^{-1/2}\beta_{q_t^{   \epsilon}}(x)$,  we   have
\begin{align*}
\epsilon^{ 1/2}   \,  \mathrm{d} w_{1,\rho_t^{   \epsilon}}(X_t^{\epsilon,0} ,X_t^{\epsilon,1} )
 =\Bigg[ Y_t^{\epsilon,0} \cdot\nabla_{x_0}w_{1,\rho_t^{   \epsilon}}
 +Y_t^{\epsilon,1}      \cdot\nabla_{x_1}w_{1,\rho_t^{   \epsilon}}    +w_{1,\rho_t^{   \epsilon}}\left(
 \frac{\beta_{q_t^{   \epsilon}}(X_t^{\epsilon,0} )}{\lambda_{\rho_t^{   \epsilon}}(X_t^{\epsilon,0}  )}
 +\frac{\beta_{q_t^{   \epsilon}}(X_t^{\epsilon,1} )}{\lambda_{\rho_t^{   \epsilon}}(X_t^{\epsilon,1} )}
 \right)\Bigg]\mathrm{d}   t.
\end{align*}
Then    
\begin{equation}\label{eq:second-CS-rearrangement-actual}
\begin{aligned}
 \epsilon^{-1/2}    (P_{\rho_t^{   \epsilon}} u_t^{   \epsilon})(X_t^{\epsilon,0} )\mathrm{d}  t
={}&-   \mathrm{d} \bigl[   \epsilon^{ 1/2}  (P_{\rho _t^{   \epsilon}}  R_{\rho_t^{   \epsilon} }u_t^{   \epsilon})(X_t^{\epsilon,0} )\bigr]
 -(P_{\rho_t^{   \epsilon} }   R_{\rho_t^{   \epsilon}} G_{\rho_t^{   \epsilon}})(X_t^{\epsilon,0} )\mathrm{d}  t\\
 &+ \mathcal{Q}_{1 } (q_t^{   \epsilon};X_t^{\epsilon,0}  ,Y_t^{\epsilon,0} )  \mathrm{d}  t
 +   \epsilon^{-1/2}    (P_{\rho_t^{   \epsilon}}^2  u_t^{   \epsilon})(X_t^{\epsilon,0}  )\mathrm{d}  t,
\end{aligned}
\end{equation}
where
\begin{align*}
 (P_{\rho_t   ^{   \epsilon}} R_{\rho_t ^{   \epsilon}} u_t^{   \epsilon})(X_t^{\epsilon,0} )
 ={}& \mathbb{E}_1   \left[w_{1,\rho_t^{   \epsilon}}(X_t ^{\epsilon,0} ,X_t ^{\epsilon,1}  )Y_t ^{\epsilon,1} \right],\\
(P_{\rho_t ^{   \epsilon} }  R_{\rho_t ^{   \epsilon} }  G_{\rho_t ^{   \epsilon}})(X_t ^{\epsilon,0} )
 = {}&\mathbb{E}_1   \left[w_{1,\rho_t ^{   \epsilon}}(X_t ^{\epsilon,0} ,X_t ^{\epsilon,1} )G_{\rho_t ^{   \epsilon} }(X _t ^{\epsilon,1} ) \right],\\
 (P_{\rho_t ^{   \epsilon}}^2 u_t ^{   \epsilon})(X_t ^{\epsilon,0} )   ={} &    \mathbb{E}_1  \left[
 R_{\rho_t ^{   \epsilon}}(X_t ^{\epsilon,0} )\phi(X_t ^{\epsilon,0} -X_t ^{\epsilon,1} )(P_{\rho_t ^{   \epsilon}} u_t ^{   \epsilon})(X_t ^{\epsilon,1} )\right] 
\end{align*}
and $\mathcal{Q}_{1}$ is  the $1$st layer  of $\mathcal{Q}$, denoted as
\begin{align*}
  \mathcal{Q}_{1} (q;x_{0},y_{0})  := {}&  \mathbb{E}_1    \Bigg[      Y ^1     \Bigg( y_{0}  \cdot\nabla_{x_0}w_{1,\rho}
 +Y^1\cdot\nabla_{x_1}w_{1,\rho_t}    +w_{1,\rho}\left(
 \frac{\beta_{q}(x_{0})}{\lambda_{\rho}(x_{0})}
 +\frac{\beta_{q}(X^1)}{\lambda_{\rho}(X^1)}
 \right)\Bigg)  \Bigg]\\
={}&    \int_{\mathbb{R}   ^{2d}} y_1\Bigg[ 
\sum_{i=0}^1 y_i\cdot g_{1,i,\rho}(x_0,  x_1)   +w_{1,\rho}(x_0, x_1)
\sum_{i=0}^1\frac{\beta_q(x_i)}{\lambda_\rho(x_i)}    \Bigg]q ( \mathrm{d} x_1 ,  \mathrm{d} y_1 ).
\end{align*}
Substituting  \eqref{eq:second-CS-rearrangement-actual} into
\eqref{eq:local-CS-corrected-position} and define
$
 Z_t^{\epsilon, (1) }:=X_t^{\epsilon,0} +  \epsilon^{ 1/2}  R_{\rho_t  ^{   \epsilon}}(X_t^{\epsilon,0}  )Y_t^{\epsilon,0} 
 + \epsilon^{ 1/2}     (P_{\rho_t^{   \epsilon}} R_{\rho_t ^{   \epsilon}}  u_t^{   \epsilon})(X_t^{\epsilon,0} ),
$
we have 
\begin{equation}\label{eq:first-nonlocal-CS-corrected-position}
\begin{aligned}
 \mathrm{d} Z_t^{ \epsilon ,  (1)  }
={}&-\bigl(R_{\rho_t ^{   \epsilon}} G_{\rho_t  ^{   \epsilon}}  +P_{\rho_t  ^{   \epsilon}} R_{\rho_t   ^{   \epsilon}} G_{\rho_t   ^{   \epsilon}}\bigr)(X_t^{   \epsilon, 0} )\mathrm{d}   t
 +( \mathcal{Q}_0(q_t^{   \epsilon};X_t^{\epsilon,0} ,Y_t^{\epsilon,0} ) +  \mathcal{Q}_1(q_t^{   \epsilon};X_t^{\epsilon,0} ,Y_t^{\epsilon,0}) )\mathrm{d} t\\
 &+R_{\rho_t^{   \epsilon}}(X_t^{\epsilon,0})\sigma  \mathrm{d} B_t^{\epsilon,0} 
 +   \epsilon^{-1/2}        (P_{\rho_t  ^{   \epsilon}}^2   u_t^{   \epsilon})(X_t^{\epsilon,0} )\mathrm{d}  t.
\end{aligned}
\end{equation}
We find that we do indeed cancel \(P_\rho u\), but leave \(P_\rho^2u\). Repeating the same procedure on the new terminal copy yields a singular chain that cannot be truncated at finite layers 
$$
 P_\rho u\longmapsto P_\rho^2u
 \longmapsto P_\rho^3u\longmapsto\cdots .
$$

In general, take mutually independent copies $(X^{\epsilon, j},Y^{\epsilon, j},B^{\epsilon, j})_{1\le j\le n}$, and set
\begin{align*}
w_{0,\rho}(x_0)=\lambda_\rho^{-1} (x_0),
\qquad   w_{n,\rho}(x_0,\ldots,x_n)
:=\frac{\prod_{j=0}^{n-1}\phi(x_j-x_{j+1})}
{\prod_{j=0}^{n}\lambda_\rho(x_j)},\quad n\ge1,
\end{align*}
and
$
g_{n,i,\rho}:=\nabla_{x_i}w_{n,\rho}.
$
Define
$$
 \mathcal{H}_n(q_t^{   \epsilon};X_t^{\epsilon,0} )
 :=   \mathbb{E}_{1:n}   \left[
 w_{n,\rho_t^{   \epsilon}}(X_t^{\epsilon,0} ,\ldots,X_t^{\epsilon,n} )Y_t^{\epsilon,  n } \right],
$$
where $ \mathbb{E}_{1:n}$ denotes the expectation with respect to the joint distribution of the independent copies indexed by $1,\ldots,n$   and 
we adopt the convention that $\mathcal{H}_0(q_t^{   \epsilon};X_t^{\epsilon,0} ) :=R_{\rho_t  ^{   \epsilon}}(X_t^{\epsilon,0} )Y_t^{\epsilon,0} $.
It follows from   $\mathrm{d} X_t^{  \epsilon, i } =  \epsilon^{-1/2}  Y_t^{ \epsilon, i }   \mathrm{d} t$ and 
$\partial_t\lambda_{\rho_t  ^{   \epsilon}}(x)=-    \epsilon^{-1/2}  \beta_{q_t^{   \epsilon}}(x)$   that
\begin{align*}
 \epsilon^{ 1/2}  \,  \mathrm{d} w_{1,\rho_t^{   \epsilon}}(X_t^{\epsilon,0}  ,\ldots,X_t ^{\epsilon,n} )
 =\Bigg[   \sum_{i=0}^n   Y_t^{ \epsilon, i  } \cdot
 g_{n,i,\rho_t^{   \epsilon}}(X_t^{\epsilon,0} ,\ldots,X_t^{\epsilon,n} ) + w_{n,\rho_t^{   \epsilon}}
   \sum_{i=0}^n\frac{\beta_{q_t^{   \epsilon}}(X_t^{ \epsilon, i } )}{\lambda_{\rho_t  ^{   \epsilon}}(X_t^{ \epsilon, i }  )}      \Bigg]\mathrm{d}   t.
\end{align*}
Therefore, combining with the velocity equation   of   \eqref{eq:scaled-system}, we obtain     for each $n \geq 0$
\begin{equation}\label{eq:single-chain-semi}
\begin{aligned}
   \epsilon^{ 1/2} \,\mathrm{d}\mathcal{H}_n(q_t^{   \epsilon};X_t^{\epsilon,0} )
={}& \mathcal{Q}_n(q_t ^{   \epsilon};X_t^{\epsilon,0} ,Y_t^{\epsilon,0}   )\mathrm{d} t
 +   \epsilon^{-1/2}  (- \mathcal{D}_n(q_t  ^{   \epsilon};X_t^{\epsilon,0} ,Y_t^{\epsilon,0} )+\mathcal{S}_n(q_t ^{   \epsilon};X_t^{\epsilon,0}))\mathrm{d}   t\\
&-  \digamma_n(q_t^{   \epsilon};X_t^{\epsilon,0})  \mathrm{d}  t
+\boldsymbol{1}_{\{n=0\}}R_{\rho_t^{   \epsilon}}(X_t^{\epsilon,0} )\sigma  \mathrm{d} B_t^{\epsilon,0} ,
\end{aligned}
\end{equation}
where we adopt the convention that   $\mathcal{D}_{0 }(q_t^{   \epsilon};X_t^{\epsilon,0} ,  Y_t^{\epsilon,0}  ) :=Y_t^{\epsilon,0} $ and  $\digamma_0(q_t^{   \epsilon};X_t^{\epsilon,0} ):=   R_{\rho_t^{   \epsilon}} G_{\rho_t^{   \epsilon}} (X_t^{\epsilon,0} )   $ ,   
\begin{align*}
 \mathcal{D}_n(q_t^{   \epsilon};X_t^{\epsilon,0} ,  Y_t^{\epsilon,0}  )
 &:=\mathbb{E}_{1:n}   \left[
 w_{n,\rho_t^{   \epsilon}}(X_t^{\epsilon,0}  ,\ldots,X_t^{\epsilon,n} )
 \lambda_{\rho_t^{   \epsilon}}(X_t^{\epsilon,n} )Y_t^{\epsilon,n} \right],\\
 \mathcal{S}_n(q_t^{   \epsilon};X_t^{\epsilon,0} )
 &:= \mathbb{E}_{1:n+1}   \left[
 w_{n,\rho_t^{   \epsilon}}(X_t^{\epsilon,0} ,\ldots,X_t^{\epsilon,n} )
 \phi(X_t^{\epsilon,n} -X_t^{ \epsilon,  n+1 } )Y_t^{ \epsilon ,  n+1   } \right],\\
 \digamma_n(q_t^{   \epsilon};X_t^{\epsilon,0} ) &:=   \mathbb{E}_{1:n}   \left[
 w_{n,\rho_t ^{   \epsilon}  }(X_t^{\epsilon,0}   ,\ldots,X_t^{\epsilon,n} )G_{\rho_t  ^{   \epsilon}}(X_t^{\epsilon,n})\right]
\end{align*} 
and $\mathcal{Q}_{n}$ is the $n$-th layer of $\mathcal{Q}$, denoted as
\begin{align*}
 &  \mathcal{Q}_{n} (q;x_{0},y_{0}) \\
   :={}&  \mathbb{E}_{1:n}    \Bigg[      Y ^n     \Bigg( y_{0}\cdot
 g_{n,0,\rho }    +  \sum_{i=1}^n   Y^i\cdot
 g_{n,i,\rho}  + w_{n,\rho}
    \frac{\beta_{q}(x_{0})}{\lambda_{\rho}(x_{0} )}    + w_{n,\rho}
   \sum_{i=1}^n\frac{\beta_{q}(X^i)}{\lambda_{\rho}(X^i )}    \Bigg)  \Bigg]\\
={}&    \int_{\mathbb{R}   ^{2nd}}     y_n\Bigg[
\sum_{i=0}^n y_i\cdot g_{n,i,\rho}(x_0,\ldots,x_n)+w_{n,\rho}(x_0,\ldots,x_n)
\sum_{i=0}^n\frac{\beta_q(x_i)}{\lambda_\rho(x_i)}
\Bigg]q^{\otimes n}( \mathrm{d} z_1,    \dots ,    \mathrm{d} z_n).
\end{align*}
Note that  \begin{align}
 w_{n+1,\rho}(x_0,\ldots,x_{n+1})
 \lambda_\rho(x_{n+1})
 =w_{n,\rho}(x_0,\ldots,x_n)
 \phi(x_n-x_{n+1}),\label{eq:weight-telescope}
\end{align}
which implies 
$
\mathcal{S}_n(q_t^{   \epsilon};X_t^{\epsilon,0} )
 =   \mathcal{D}_{n+1}(q_t ^{   \epsilon} ;X_t ^{\epsilon,0}  ,Y_t^{\epsilon,0} )
$
for each $n \geq 0$.   Then   we have
\begin{align}
 &\sum_{n=0}^N(-  \mathcal{D}_{n }(q_t ^{   \epsilon}  ;X_t^{\epsilon,0} ,  Y_t^{\epsilon,0}  )    +  \mathcal{S}_n(q_t ^{   \epsilon};X_t^{\epsilon,0} )   )\notag\\
=&-   \mathcal{D}_{0 }(q_t^{   \epsilon};X_t^{\epsilon,0} ,  Y_t^{\epsilon,0}  )       +   \mathcal{D}_{N+1 }(q_t^{   \epsilon};X_t^{\epsilon,0} ,  Y_t^{\epsilon,0}  )      \notag\\
=&  -Y_t^{\epsilon,0} +
 \mathscr{R}_{\rho_t  ^{   \epsilon}}^{(N)}(q_t ^{   \epsilon};X_t^{\epsilon,0} ),\label{eq:fast-telescope}
\end{align}
where  
\begin{align*}
\mathscr{R}_\rho^{(N)}(q;x_0)
  :={}&\int _{  \mathbb{R}  ^{2(N+1)d}}  \lambda_\rho(x_{N+1})w_{N+1,\rho}(x_0,\ldots,x_{N+1})
 y_{N+1}q^{\otimes(N+1)}(\mathrm{d} z_1,   \dots,   \mathrm{d}    z_{N+1})\\
= {}& (P_\rho^{N+1}u)(x_0). 
\end{align*}
Define the truncated correction
\begin{align}\label{eq:XiN-def}
 \Xi^{(N)}(q;x_{0},y_{0})
 :={} &  R_\rho(x_{0})y_{0}   +  \sum_{n=1}^N   \mathbb{E}_{1:n}   \left[
 w_{n,\rho_t }(x_{0}, X_t^1,\ldots,X_t^n)Y_t^n \right] \notag\\
   = {}& \sum_{n=0}^N   \int _{  \mathbb{R}  ^{2nd}} w_{n,\rho}(x_0,\ldots,x_n)y_n\,
 q^{\otimes n}(\mathrm{d} z_1,   \dots ,   \mathrm{d} z_n) \notag   \\
  = {}&    R_\rho(x_{0})y_{0} + \sum_{n=1}^N(P_\rho^n    R_\rho u)(x_{0})
 \end{align}
 and  the truncated corrected position
 $$
 Z_t^{\epsilon,  (N) }:  =X_t^{\epsilon,0}    + \epsilon^{ 1/2}    \Xi ^{(N)} (q_{t}^{   \epsilon};X_{t}^{\epsilon,0} ,Y_{t}^{\epsilon,0} ) .  
$$
 Set 
\begin{align*}
 F_{\rho }^{(N)} (x_0)   := {}  &  R_\rho G_\rho (x_0)  +  \sum_{n=1}^N    \mathbb{E}_{1:n}   \left[
 w_{n,\rho }(x_0,\ldots,X ^n)G_{\rho }(X ^n)\right]\\
 ={}  & \sum_{n=0}^N   \int _{   \mathbb{R} ^{nd}} w_{n,\rho}(x_0,\ldots,x_n)
 G_\rho(x_n)\rho^{\otimes n}(\mathrm{d}   x_1,   \dots,   \mathrm{d} x_n)\\
 ={}&\sum_{n=0}^N    P_\rho^n   R_\rho G_\rho (x_0)   
 \end{align*}
 and
 $$
  \mathcal{Q}^{(N)} (q;x_{0},y_{0})  :=  \sum_{n=0}^N      \mathcal{Q}_n   (q;x_{0},y_{0}) .
$$
Taking the sum of \eqref{eq:single-chain-semi} from $n=0$ to $N$, we have
\begin{equation}\label{eq:N-rearrangement-motivation}
\begin{aligned}
 \mathrm{d}   Z_t^{ \epsilon,  (N) } =  {}&  \mathrm{d} \Big(X_t^{\epsilon,0} 
 +    \epsilon^{  1/2}   R_{\rho_t  ^{   \epsilon}}(X_t^{\epsilon,0} )Y_t^{\epsilon,0} 
 +    \epsilon^{  1/2}     \sum_{n=1}^N(P_{\rho_t  ^{   \epsilon}  }^n  R_{\rho_t  ^{   \epsilon}} u_t^{   \epsilon})(X_t^{\epsilon,0} )   \Big) \\
 ={}&-\sum_{n=0}^N(P_{\rho_t  ^{   \epsilon}}^n  R_{\rho_t^{   \epsilon}} G_{\rho_t^{   \epsilon}})(X_t^{\epsilon,0} )\mathrm{d} t
 +\sum_{n=0}^N  \mathcal{Q}_n (q_{t}^{   \epsilon};X_t^{\epsilon,0} ,Y_t^{\epsilon,0} ) \mathrm{d} t
 +R_{\rho_t ^{   \epsilon}}(X_t^{\epsilon,0} )\sigma  \mathrm{d} B_t^{\epsilon,0}  \\
   &+   \epsilon^{-1/2}   (P_{\rho_t ^{   \epsilon}}^{N+1}   u_t^{   \epsilon})(X_t^{\epsilon,0} )\mathrm{d}  t\\
 ={}& - F_{\rho_t^{   \epsilon} }^{(N)} (X_t^{\epsilon,0}  )     \mathrm{d} t
 +   \mathcal{Q}^{(N)} (q_{t} ^{   \epsilon};X_t^{\epsilon,0} ,Y_t^{\epsilon,0} ) \mathrm{d} t
 +R_{\rho_t^{   \epsilon}}(X_t^{\epsilon,0} )\sigma  \mathrm{d} B_t^{\epsilon,0}    + \epsilon^{-1/2}       (P_{\rho_t^{   \epsilon}}^{N+1}   u_t^{   \epsilon})(X_t^{\epsilon,0} )  \mathrm{d}  t  .
\end{aligned}
\end{equation}

\subsection{    Finite truncation limits } \label{estimates-truncation-limits}

From  equation  \eqref{eq:N-rearrangement-motivation},
the truncated corrected position $Z_t^{\epsilon,(N)}  :=X_t^ \epsilon   +  \epsilon^{ 1/2}   \Xi ^{ (N)} (q_{t}^ \epsilon;X_{t}^ \epsilon,Y_{t}^ \epsilon) $ satisfies the adapted SDE
\begin{equation}\label{eq:ZN-decomp}
 \mathrm{d} Z_t^{\epsilon,(N)}
 =
 -F_{\rho_t^\epsilon}^{(N)}(X_t^\epsilon)\mathrm{d} t
 + \mathcal{Q}^{(N)}(q_t^\epsilon;X_t^\epsilon,Y_t^\epsilon)\mathrm{d} t
 +\lambda_{\rho_t^\epsilon}^{-1}  (X_t^\epsilon)\sigma   \mathrm{d} B_t   
   +  \epsilon^{-1/2} \mathscr{R}_{\rho_t^\epsilon}^{(N)}       (q_t^\epsilon;X_t^\epsilon)\mathrm{d} t,
\end{equation}
where  $ \mathscr{R}_{\rho_t^\epsilon}^{(N)}
       (q_t^\epsilon;X_t^\epsilon)   =   (P_{\rho_t^\epsilon}^{N+1}   u_t^\epsilon)(X_t^\epsilon).$

\begin{lemma}\label{lem:finite-corrector}
Under \textup{\textbf{(A1)}-\textbf{(A2)}},   for each fixed $\epsilon > 0$,    the   continuous adapted process $Z_t^{\epsilon}   :  =X _t^{\epsilon}   +   \epsilon^{ 1/2}   \Xi   (q _t^{\epsilon};X _t^{\epsilon},Y_t ^{\epsilon})  $ satisfies 
 \begin{equation}
 \mathbb{E}   \Big[\sup_{t\le T}|Z_t^{\epsilon, (N)}-Z_t^\epsilon|^{2} \Big] \rightarrow    0,    \qquad   N  \rightarrow  \infty .
\end{equation} 
Moreover  $Z_t^{\epsilon}$ satisfies the following equation
\begin{equation}\label{eq:Z-decomp}
 \mathrm{d} Z_t^\epsilon
 =-F_{\rho_t^\epsilon}(X_t^\epsilon)\mathrm{d} t
 + \mathcal{Q}(q_t^\epsilon;X_t^\epsilon,Y_t^\epsilon)\mathrm{d} t
 +\lambda_{\rho_t^\epsilon}^{-1}   (X_t^\epsilon)   \sigma  \mathrm{d} B_t,
\end{equation}
where
\begin{equation}\label{eq:F-series}
 F_\rho :=\sum_{n=0}^\infty P_\rho^n   R_\rho G_\rho=A_\rho^{-1}G_\rho.
\end{equation}
\end{lemma}

\proof

For any $N'> N \geq 1,$
from   Lemmas  \ref{lem:contraction},      \ref{lem:operator-insertion}   and   \ref{lem:kernel-sum},     we have  
\begin{align*}
& \|\mathscr{R}_\rho^{(N')}(q;\cdot)   -\mathscr{R}_\rho^{(N)}(q;\cdot)  \|_{L^2(\rho)}   =  \|  (P_\rho^{N'+1}  - P_\rho^{N+1})  u \|_{L^2(\rho)}
\le C r_*^{N+1}\|u   \|_{L^2(\rho)},\\
& \|F_\rho ^{(N')}    -F_\rho^{(N)}\|_{L^2(\rho)}   =     \| \sum_{n=N+1}^{N'} P_\rho^n   R_\rho G_\rho\|_{L^2(\rho)}
\le \frac{C r_*^{N+1}}{1-r_*}
       \|G_\rho\|_{L^2(\rho)},\\
& \|\mathcal{Q}^{(N')}-\mathcal{Q}^{(N)}\|_{L^2(q)}
\le C(1+m_4(q)^{1/2})
   \sum_{n>N}(n+1)^2r_*^{(n-2)_+}.
\end{align*}
For fixed $\epsilon>0$, by Assumption \textbf{(A2)} and  \eqref{eq:fixed-eps-moments}, there exists $C_{\epsilon,T}<\infty$ such that for $t\le T$,
$$
\|u_t  ^\epsilon   \|_{L^2(q_t  ^\epsilon )}+\|G_{\rho_t  ^\epsilon }\|_{L^2(q_t  ^\epsilon )}  +   m_4(q_t  ^\epsilon )^{1/2}     \le C_{\epsilon,T}.
$$
Furthermore, we  obtain
\begin{align}
 \mathbb{E}  \bigg[\sup_{t\le T}\left|
 \epsilon^{-1/2} \int_0^t  (\mathscr{R}_{\rho_s^\epsilon}^{(N')}(q_s^\epsilon;X_s^\epsilon)   - \mathscr{R}_{\rho_s^\epsilon}^{(N)}(q_s^\epsilon;X_s^\epsilon)    )\mathrm{d} s  \right|^{2}  \bigg]
 &\le C_{\epsilon,T}r_*^{N+1},\label{eq:R-tail-time}\\
 \mathbb{E}   \bigg[ \sup_{t\le T}\left|
 \int_0^t(F_{\rho_s^\epsilon}^{(N')}   -F_{\rho_s^\epsilon}^{(N)})(X_s^\epsilon)\mathrm{d} s\right|^{2}  \bigg]
 &\le C_{\epsilon,T}r_*^{N+1},\label{eq:F-tail-time}\\
 \mathbb{E}\bigg[\sup_{t\le T}\left|
 \int_0^t(\mathcal{Q}^{(N')}(q_s^\epsilon;X_s^\epsilon,Y_s^\epsilon)   -\mathcal{Q}^{(N)}(q_s^\epsilon;X_s^\epsilon,Y_s^\epsilon))\mathrm{d} s\right|^{2}   \bigg]
 &\le C_{\epsilon,T}
 \sum_{n>N}(n+1)^2r_*^{(n-2)_+} \label{eq:Q-tail-time}
\end{align}
and 
\begin{equation}\label{eq:Xi-initial-tail}
  \mathbb{E}[|\epsilon^{1/2}  (\Xi^{(N')}(q_0^\epsilon;X_0^\epsilon,Y_0^\epsilon)   -\Xi ^{(N)}(q_0^\epsilon;X_0^\epsilon,Y_0^\epsilon))|  ^{2}]
 \le C_{\epsilon,T}   r_{*}^{N+1}.
\end{equation}
Consequently, \eqref{eq:R-tail-time}-\eqref{eq:Xi-initial-tail} imply that $\{Z^{\epsilon,N}\}_N$ is a Cauchy sequence in $L^2(\Omega;\mathcal{C}([0,T];\mathbb{R}^d))$. By completeness of this space, there exists a  limit $Z^\epsilon$ such that
$$
\mathbb{E}   \Big[\sup_{t\le T}|Z_t^{\epsilon,N}-Z_t^\epsilon|^{2} \Big]  \rightarrow  0,    \qquad   N  \rightarrow  \infty.
$$
Moreover, for fixed $\epsilon$, since   $
\mathbb{E} \Big[      \sup_{t\le T}|\Xi^{(N)}(q_t^\epsilon;X_t^\epsilon,Y_t^\epsilon)  -\Xi    (q_t^\epsilon;X_t^\epsilon,Y_t^\epsilon)|^{2}  \Big]   \rightarrow  0,      N  \rightarrow  \infty 
$, $Z_t^\epsilon$ is denoted by    $ Z_t^\epsilon   =  X_t^\epsilon   + \epsilon^{ 1/2}    \Xi   (q_t^\epsilon;X_t^\epsilon,Y_t^\epsilon)$      and  a continuous adapted process.
On the other hand,  by a similar discussion as above for the estimates, we may take the limit as $N  \rightarrow   \infty$ in equation \eqref{eq:ZN-decomp} and obtain that $Z_t^\epsilon$ satisfies the following equation
\begin{equation*} 
 \mathrm{d} Z_t^\epsilon
 =-F_{\rho_t^\epsilon}(X_t^\epsilon)\mathrm{d} t
 + \mathcal{Q}(q_t^\epsilon;X_t^\epsilon,Y_t^\epsilon)\mathrm{d} t
 +\lambda_{\rho_t^\epsilon}^{-1}    (X_t^\epsilon)  \sigma  \mathrm{d} B_t,
\end{equation*}
where
$ 
 F_\rho =\sum  \limits_{n=0}^\infty P_\rho^n   R_\rho G_\rho=A_\rho^{-1}G_\rho.
$
\endproof

\subsection{Moment estimates and  tightness}\label{estimates-tightness}
In the following, we provide several uniform  estimates that are crucial for proving the tightness of $\{X^\epsilon\}_{0<\epsilon\le 1}$.

\begin{proposition}\label{prop:moment-closure}
Under \textup{\textbf{(A1)}-\textbf{(A2)}}, for each $T>0$,
\begin{equation}\label{eq:uniform-Y4-X2}
 \sup_{0<  \epsilon  \le1}\sup_{t\le T} \mathbb{E}[|Y_t^\epsilon |^4]
 +\sup_{0<\epsilon \le1}  \mathbb{ E}\Big[\sup_{t\le T}|Z_t^\epsilon |^2   \Big]
 +\sup_{0<\epsilon \le1}\sup_{t   \le T}   \mathbb{E}[|X_t^\epsilon |^2]
 \le C_T,
\end{equation}
and
\begin{equation}\label{eq:Q-uniform-derived}
 \sup_{0<\epsilon \le  1}\sup_{t\le T}  \mathbb{E}[| \mathcal{Q}(q_t^\epsilon;X_t^\epsilon,Y_t^\epsilon) |^2]  \le C_T.
\end{equation}
\end{proposition}

\proof

Applying It\^{o}'s formula, we obtain
$$
\begin{aligned}
\mathrm{d}   |Y_t^\epsilon|^4
=&\Bigg[4\epsilon^{-1}|Y_t^\epsilon|^2   \langle Y_t^{\epsilon },   \int _{   \mathbb{R}^{2d} } \phi(X_t^\epsilon-x')(y'-Y_t^\epsilon)\,q_t^\epsilon(\mathrm{d} x',\mathrm{d}  y')-|Y_t^\epsilon|^2    \rangle
\\
& -4  \epsilon^{-1/2} |Y_t^\epsilon|^2 \langle  Y_t^{\epsilon },  G_{\rho_t^\epsilon}(X_t^\epsilon)\rangle 
+2(d+2)\sigma^{2} \epsilon^{-1 } |Y_t^\epsilon|^2 \Bigg]\mathrm{d}     t  
+4\sigma  \epsilon^{-1/2} |Y_t^\epsilon|^2  \langle Y_t^{\epsilon },  \mathrm{d}    B_t\rangle.
\end{aligned}
$$
By symmetry, we have
\begin{align*}
&4\epsilon^{-1}   \mathbb{E}   \left[ |Y_t^\epsilon|^2   Y_t^{\epsilon }   \cdot  \int  _{    \mathbb{R} ^{2d}} \phi(X_t^\epsilon-x')(y'-Y_t^\epsilon)\,q_t^\epsilon(\mathrm{d} x',\mathrm{d}  y') \right]  \\
={} &-2\epsilon^{-1}\mathbb{ E}   \left[ \mathbb{\tilde{E}}   \left[
 \phi(X_t^\epsilon-{\tilde{X} }_t^\epsilon)\bigl(|Y_t^\epsilon|^2Y_t^\epsilon-|{\tilde{Y} }_t^\epsilon|^2{\tilde{Y} }_t^\epsilon\bigr) \cdot(Y_t^\epsilon-{\tilde{Y} }_t^\epsilon) \right]  \right]\\
 \le {} &0,
\end{align*}
where we used that $z\mapsto |z|^2z$ is monotone.
By Young's inequality 
\begin{equation}\label{eq:force-young-fourth}
\epsilon^{-1/2} |Y_t^\epsilon|^3|G_{\rho  _t^\epsilon}(X_t^\epsilon)|
 \le \xi  \epsilon^{-1}|Y_t^\epsilon|^4+C_\xi      \epsilon|G_{\rho_t^\epsilon}(X_t^\epsilon)|^4,
\end{equation}
where $\xi  > 0$ is chosen later.
Under Assumption \textbf{(A2)}   
$$
 |G_{\rho_t^\epsilon}(X_t^\epsilon)|^4
 \le \mathbb{\tilde{E}}  \Big[  \bigl|\nabla K(X_t^\epsilon-\tilde{X}_t^\epsilon )\bigr|^4    \Big]
 \le C \left(1+\mathbb{\tilde{E}} \Big[|X_t^\epsilon-\tilde{X}_t^\epsilon |^4   \Big]   \right),
$$
and consequently,
$$
 \mathbb{E}\Big[|G_{\rho_t^\epsilon}(X_t^\epsilon)|^4 \Big]   \le C\Bigl(1+ \mathbb{E} \Big[|X_t^\epsilon|^4 \Big] \Bigr).
$$
Then  we have 
\begin{equation*} 
 \epsilon^{-1/2}    \mathbb{E}\Big[    |Y_t^\epsilon|^3|G_{\rho  _t^\epsilon}(X_t^\epsilon)|   \Big] 
 \le \xi  \epsilon^{-1}    \mathbb{E}\Big[  |Y_t^\epsilon|^4  \Big] +  C_\xi      \epsilon  \Bigl(1+ \mathbb{E} \Big[|X_t^\epsilon|^4 \Big] \Bigr).
\end{equation*}
Using  Young's inequality,  we obtain  
$$2(d+2)\sigma^{2} \epsilon ^{-1}  \mathbb{E}\Big[    |Y_t^\epsilon|^2  \Big]    \leq    2(d+2)\sigma^{2} \epsilon ^{-1}  \xi
\mathbb{E}\Big[    |Y_t^\epsilon|^4  \Big]     +  C_{\xi  }\epsilon ^{-1}  . $$
Taking $\xi$ sufficiently small   yields
\begin{equation}\label{eq:M4-differential}
 \mathrm{d}     \mathbb{E} \Big[  |Y_t^\epsilon   |^4\Big]   \le    -\frac{c}{\epsilon}  \mathbb{E} \Big[  |Y_t^\epsilon   |^4\Big]  \mathrm{d}     t  + \frac{C}{\epsilon}  \mathrm{d}     t  +  C\epsilon \mathbb{E} \Big[|X_t^\epsilon|^4 \Big] \mathrm{d}     t.
\end{equation}

On the other hand, from
$
 X_t^\epsilon=X_0^\epsilon+  \epsilon^{-1/2}  \int_0^t   Y_s^\epsilon\,\mathrm{d}  s
$
and H\"older's inequality, we obtain
\begin{equation}\label{eq:N4-crude}
\mathbb{E} \Big[|X_t^\epsilon|^4 \Big]    \le C \mathbb{E} \Big[|X_0^\epsilon|^4 \Big]   +C_T \epsilon^{-2}  \int_0^t \mathbb{E} \Big[  |Y_s^\epsilon   |^4\Big]   \,\mathrm d s.
\end{equation}
Substituting \eqref{eq:N4-crude} into \eqref{eq:M4-differential} and using  differential inequality, we have
\begin{align*}
 \mathbb{E} \Big[  |Y_t^\epsilon   |^4\Big]
 &\le  \mathbb{E} \Big[  |Y_0^\epsilon   |^4\Big] e^{-ct/\epsilon}
 +\int_0^t e^{-c(t-s)/\epsilon}
 \left[ \frac{C}{\epsilon}+C\epsilon \mathbb{E} \Big[|X_0^\epsilon|^4 \Big]  +   C_{T}   \epsilon^{-1}\int_0^s \mathbb{E} \Big[  |Y_r^\epsilon   |^4\Big] \,\mathrm dr \right]\mathrm ds \notag\\
 &=: \mathbb{E} \Big[  |Y_0^\epsilon   |^4\Big] e^{-ct/\epsilon}
 +I_1+I_2+I_3.
\end{align*}
 By  standard   calculation   we have 
\begin{align*}
 I_1
 &=\frac{C}{\epsilon}\int_0^t e^{-c(t-s)/\epsilon}\,\mathrm{d}   s\\
 &=\frac{C}{\epsilon}\cdot\frac{\epsilon}{c} \left(1-e^{-ct/\epsilon}\right)\\
 &\le C ,
\end{align*}
\begin{align*}
 I_2
 &=C\epsilon \mathbb{E} \Big[|X_0^\epsilon|^4 \Big]\int_0^t e^{-c(t-s)/\epsilon}\,\mathrm{d} s\\
 &=C\epsilon \mathbb{E} \Big[|X_0^\epsilon|^4 \Big]    \frac{\epsilon  \left(1-e^{-ct/\epsilon}\right)}{c}
\\
 &\le  C \epsilon^2 \mathbb{E} \Big[|X_0^\epsilon|^4 \Big] 
\end{align*}
and 
\begin{align*}
 I_3 ={}&\frac{C}{\epsilon}  \int_0^t   e^{-c(t-s)/\epsilon}
   \left( \int_0^s    \mathbb{E} \Big[  |Y_r^\epsilon   |^4\Big]   \mathrm{d}   r\right)\mathrm{d}  s\\
= {}&\frac{C}{\epsilon}\int_0^t     \mathbb{E} \Big[  |Y_r^\epsilon   |^4\Big]
   \left(\int_r^te^{-c(t-s)/\epsilon}\mathrm{d}  s\right)\mathrm{d}  r  \\
 \le{} & C\int_0^t   \mathbb{E} \Big[  |Y_r^\epsilon   |^4\Big]   \mathrm{d}  r,
\end{align*}
  where we used Fubini's theorem.
Therefore
$$
 \mathbb{E} \Big[  |Y_t^\epsilon   |^4\Big]  \le C_T+C_T\int_0^t    \mathbb{E} \Big[  |Y_r^\epsilon   |^4\Big] \,\mathrm{d} r,
$$
and by Gronwall's lemma, for $0\le t\le T$,
$$
 \mathbb{E} \Big[  |Y_t^\epsilon   |^4\Big]   \le C_T e^{C_T t} \le C_T,
$$
which  implies    
$$ \sup_{0<  \epsilon  \le1}\sup_{t\le T}
 \mathbb{E}\Big[|Y_t^\epsilon |^4   \Big]
 \le C_T.$$
Combining with Lemma  \ref{lem:kernel-sum} gives
\begin{equation}\label{eq:Xi-pointwise-fourth}
 \sup_{0<  \epsilon  \le1}\sup_{t\le T}   \mathbb{E}\Big[|\Xi(q_t^\epsilon;X_t^\epsilon,Y_t^\epsilon)|^4 \Big]    \le C_T       \qquad    \sup_{0<  \epsilon  \le1}\sup_{t\le T}    \mathbb{E}\Big[   \mathcal{Q}(q_t^\epsilon;X_t^\epsilon,Y_t^\epsilon)|^2 \Big]    \le C_T .
\end{equation}

Using  $ 
 F_\rho =\sum  \limits_{n=0}^\infty P_\rho^n   R_\rho G_\rho=A_\rho^{-1}G_\rho 
$   and    $\|A_\rho^{-1}\|_{L^2(\mu_\rho) \to L^2(\mu_\rho)}
\le \frac{1}{   1-r_{*} }$,   we have 
\begin{align*}
  \mathbb{E}[|F_{\rho_t^ \epsilon    }(X_t^\epsilon)|^2  ] ={} &  \|F_{\rho_t^ \epsilon    } \|^2  _{L^2( \rho_t^ \epsilon)} \\ 
  = {}&  \|  A_{\rho_t^ \epsilon }^{-1}G_{\rho_t^ \epsilon }   \|^2  _{L^2( \rho_t^ \epsilon)}  \\ 
 \le {}&  C\|   G_{\rho_t^ \epsilon }   \|^2  _{L^2( \rho_t^ \epsilon)}  \\  
 \le {}&C  \mathbb{E}[|G_{\rho_t^\epsilon }(X_t^\epsilon )|^2 ]\\
 \le {}&  C\bigl(1+ \mathbb{E} [  |X_t^\epsilon|^2]  \bigr).
\end{align*}
From  It\^{o}'s formula,  we  get 
\begin{align*} 
\mathrm{d} |Z_t^\epsilon|^2 &= \Big[ 2\big\langle Z_t^\epsilon,\,-F_{\rho_t^\epsilon}(X_t^\epsilon)+\mathcal{Q}(q_t^\epsilon;X_t^\epsilon,Y_t^\epsilon)\big\rangle + \sigma ^2 \lambda_{\rho_t^\epsilon}^{-2}   (X_t^\epsilon) \Big]\mathrm{d} t \\ &\quad + 2\lambda_{\rho_t^\epsilon}^{-1} (X_t^\epsilon) \big\langle Z_t^\epsilon,\,\sigma\,\mathrm{d} B_t  \big\rangle, 
\end{align*} 
By Burkholder-Davis-Gundy's inequality, the Cauchy--Schwarz inequality, and   $ X_t^\epsilon =  Z_t^\epsilon - \epsilon^{ 1/2}   \Xi   (q_t^\epsilon;X_t^\epsilon,Y_t^\epsilon)   $, we obtain
\begin{align*}
 \mathbb{E}  \Big [ \sup_{r\le t}|Z_r^\epsilon|^2  \Big]
\le {}&  \mathbb{E}\Big [|Z_0^\epsilon|^2 \Big]  +C_T\int_0^t   \mathbb{E}\Big [\sup_{r\le s}|Z_r^\epsilon|^2  \Big] \mathrm{d} s
 +C_T\int_0^t   \mathbb{E}\Big [|F_{\rho_s^\epsilon}(X_s^\epsilon)|^2 \Big] \mathrm{d} s\notag\\
 {}&  +C_T\int_0^t\mathbb{E}   \Big  [| \mathcal{Q}(q_s^\epsilon;X_s^\epsilon,Y_s^\epsilon)|^2  \Big]  \mathrm{d}  s+\frac{1}{2}  \mathbb{E}\Big [\sup_{r\le t}|Z_r^\epsilon|^2 \Big  ] +C_T\\
\le {}&   \mathbb{E}\Big [|Z_0^\epsilon|^2 \Big]  +C_T\int_0^t   \mathbb{E}\Big [\sup_{r\le s}|Z_r^\epsilon|^2  \Big] \mathrm{d} s
 +C_T\int_0^t    \bigl(1+ \mathbb{E} [  |X_s^\epsilon|^2]  \bigr) \mathrm{d} s\notag\\
 {}& +C_T\int_0^t\mathbb{E}   \Big  [| \mathcal{Q}(q_s^\epsilon;X_s^\epsilon,Y_s^\epsilon)|^2  \Big]  \mathrm{d}  s+\frac{1}{2}  \mathbb{E}\Big [\sup_{r\le t}|Z_r^\epsilon|^2 \Big  ] +C_T\\
 \le {}& C_T+C_T\int_0^t\mathbb{E}  \Big[  \sup_{r\le s}|Z_r^\epsilon|^2  \Big]  \mathrm{d} s
 +C_{T}   \epsilon\sup_{s\le T}\mathbb{E}\Big  [|\Xi(q_s^\epsilon;X_s^\epsilon,Y_s^\epsilon)|^2  \Big]+\frac{1}{2} \mathbb{E}\Big  [\sup_{r\le t}|Z_r^\epsilon|^2  \Big  ]\\
   \le {} & C_T+C_T\int_0^t  \mathbb{ E}  \Big [\sup_{r\le s}|Z_r^\epsilon|^2   \Big]\mathrm{d}s
 + \frac{1}{2} \mathbb{E}\Big  [\sup_{r\le t}|Z_r^\epsilon|^2\Big ],
\end{align*}
where we have used \eqref{eq:Xi-pointwise-fourth} and the inequality $|X_t^ \epsilon  |^2   \le2    |Z_t^ \epsilon  |^2+ 2\epsilon|\Xi (q_t^\epsilon;X_t^\epsilon,Y_t^\epsilon)  |^2$.
Applying Gronwall's lemma yields
$$\sup_{0<\epsilon \le1}  \mathbb{ E}\Big[\sup_{t\le T}|Z_t^\epsilon |^2   \Big]
 \le C_T,$$
which  implies 
$$ 
  \sup_{0<\epsilon \le1}\sup_{t   \le T}   \mathbb{E}[|X_t^\epsilon |^2]
 \le C_T.
 $$
\endproof

With the uniform moment estimate at hand, we establish the tightness of the family $\{X^\epsilon\}_{0<\epsilon\le 1}$.

\begin{proposition}\label{prop:tightness-complete}
Under \textup{\textbf{(A1)}-\textbf{(A2)}},  we have 
\begin{equation}\label{eq:Xi-ucp-L2}
 \mathbb{E}  \Big[    \sup_{t\le T}|\epsilon^{ 1/2}  \Xi(q_t^\epsilon;X_t^\epsilon,Y_t^\epsilon) |^2   \Big] \leq  C_T \epsilon^{1/2}   ,
\end{equation}
and the family \(\{X^ \epsilon   \}_{0< \epsilon   \le1}\) is  tight in \(\mathcal{C}([0,T];\mathbb{R}^d)\). Moreover, we also obtain
$$
\sup_{0< \epsilon  \le1}\mathbb{E}  \Big[\sup_{t\le T}|X_t^ \epsilon|^2   \Big]  \le C_T.
$$
\end{proposition}

\proof

We prove tightness in two steps.

 \noindent \textbf{Step 1: Tightness of $\{Z^\epsilon\}_{0< \epsilon   \le1}$.}
  Let $\tau \le T$ be a stopping time and $0 \le h \le 1$.
From \eqref{eq:Z-decomp},
$$
Z_{(\tau+h)  \wedge   T}^\epsilon - Z_\tau^\epsilon
= -\int_\tau^{(\tau+h)  \wedge   T} F_{\rho_s^\epsilon}(X_s^\epsilon) \mathrm{d}  s
+ \int_\tau^{(\tau+h)  \wedge   T} \mathcal{Q}(q_s^\epsilon;X_s^\epsilon,Y_s^\epsilon) \mathrm{d}   s
+ \int_\tau^{(\tau+h)  \wedge   T} \lambda_{\rho_t^\epsilon}^{-1} (X_s^\epsilon) \sigma \mathrm{d}   B_s.
$$
Jensen's inequality, It\^{o}'s isometry, and Lemma \ref{lem:kernel-sum} yield
\begin{align*}
 &\mathbb{E}  \Big[   |Z_{(\tau+h)  \wedge   T}^\epsilon - Z_\tau^\epsilon|^{2}\Big  ] \\
\leq  {} &   3\mathbb{E}  \Big[   \Big|\int_\tau^{(\tau+h)  \wedge   T} F_{\rho_s^\epsilon}(X_s^\epsilon) \mathrm{d}  s\Big|^{2}\Big  ]
+ 3\mathbb{E}  \Big[   \Big|\int_\tau^{(\tau+h)  \wedge   T} \mathcal{Q}(q_s^\epsilon;X_s^\epsilon,Y_s^\epsilon) \mathrm{d}   s\Big|^{2}\Big  ]
+ 3   \mathbb{E}  \Big[   \Big|\sigma \int_\tau^{(\tau+h)  \wedge   T} \lambda_{\rho_t^\epsilon}^{-1} (X_s^\epsilon) \mathrm{d}   B_s   \Big|^{2}\Big  ]\\
\leq  {} &   3  h\int_{0}^{   T}  \mathbb{E}  \Big[   \Big|    F_{\rho_s^\epsilon}(X_s^\epsilon) \Big|^{2}\Big  ] \mathrm{d}  s
+ 3 h   \int_{0}^{   T}   \mathbb{E}  \Big[   \Big| \mathcal{Q}(q_s^\epsilon;X_s^\epsilon,Y_s^\epsilon)\Big|^{2}\Big  ] \mathrm{d}   s     + 3\sigma^{2}  h   \\
\leq  {}  &   C_T h ,
\end{align*}
where  we have used  Proposition  \ref{prop:moment-closure}.
Hence for any $ \xi > 0$,
$$
 \mathbb{P}\big( |Z_{(\tau+h)  \wedge   T}^\epsilon - Z_\tau^\epsilon| \geq \xi \big)
\le \frac{C_T}{\xi^2} h  \rightarrow   0,   \qquad    h  \rightarrow  0  
$$
uniformly in $\epsilon$ and $\tau$. Together with 
$\sup_{   0    <   \epsilon   \leq   1    } \mathbb{E} \Big[  \sup_{t   \le T} |Z_t^\epsilon|^2 \Big  ]     \le C_T$,
the Aldous criterion \cite[Theorem 16.10]{Billingsley1999} yields tightness of $\{Z^\epsilon\}_{0< \epsilon   \le1}$ in 
$\mathcal{C}([0,T];  \mathbb{R}^d)$.

\medskip
 \noindent \textbf{Step 2: Tightness of $\{X^\epsilon\}_{0<\epsilon\le 1}$.}
 Set the grid $t_k := k\epsilon \wedge T$, $k = 0, \ldots, K_\epsilon := \lfloor T/\epsilon \rfloor $.
 For   grid points,   from  \eqref{eq:Xi-pointwise-fourth},  we have 
\begin{align}\label{eq:grid-points}
  \mathbb{E} \Big  [   \max_k   |   \epsilon^{ 1/2}      \Xi   (q_{t_k}^\epsilon;X_{t_k}^\epsilon,Y_{t_k}^\epsilon)      |^2\Big]
 \le  {}  & \left(\sum_k  \mathbb{E}  \Big  [  |  \epsilon^{ 1/2}      \Xi   (q_{t_k}^\epsilon;X_{t_k}^\epsilon,Y_{t_k}^\epsilon) |^4\Big] \right)^{1/2}   \notag  \\
 \le {} &(C_T   \epsilon^{2}     (1+   K_\epsilon ))^{1/2}   \notag \\
 \le  {}&   C_T\epsilon^{1/2}.
\end{align}
 For  between grid points,  
 \begin{align}
 \label{eq:U-block}
\epsilon^{ 1/2}    \Xi   (q_{t }^\epsilon;X_{t }^\epsilon,Y_{t }^\epsilon) -\epsilon^{ 1/2}    \Xi   (q_{t_k}^\epsilon;X_{t_k}^\epsilon,Y_{t_k}^\epsilon)
 =    {}&   -  \epsilon^{-1/2} \int_{t_k}^t   Y_s ^\epsilon  \mathrm{d}  s
 -\int_{t_k}^t  F_{\rho_s^\epsilon}(X_s^\epsilon)\mathrm{d}  s
 +\int_{t_k}^t   \mathcal{Q} (q_{t }^\epsilon;X_{t }^\epsilon, Y_{t }^\epsilon)\mathrm{d}  s    \notag \\
 {} &   +\sigma\int_{t_k}^t\lambda_{\rho_s ^\epsilon }^{-1} (X_s^\epsilon)  \mathrm{d}  B_s.
\end{align}
 By   H\"{o}lder's inequality   and    Proposition  \ref{prop:moment-closure},  we have 
 \begin{align*}
 \mathbb{E}  \Big  [   \sup_{t\in[t_k,t_{k+1}]}
 \Big|  \epsilon^{-1/2}    \int_{t_k}^t   Y_s^\epsilon   \mathrm{d} s   \Big|^4   \Big  ]     
 \le  {}& \epsilon^{-2}  \mathbb{E}  \Big[    
  \Big( \int_{t_k}^{t_{k+1}}  |Y_s^\epsilon |^{4}  \mathrm{d} s   \Big)  \Big(\int_{t_k}^{t_{k+1}}     \mathrm{d} s \Big)^{3}     \Big]   \\
  =   {}&     \epsilon^{-2}\epsilon^3   \int_{t_k}^{t_{k+1}}    \mathbb{E}   \Big[   |Y_s^\epsilon|^4 \Big]     \mathrm{d}  s\\
 \le    {}&   C_T   \epsilon^2.
\end{align*}
Applying the inequality $  \mathbb{E} [ \max_k A_k^2 ]\le \left(\sum_k  \mathbb{E}[ A_k^4]  \right)^{1/2}$, it follows that
 \begin{align*} 
  \mathbb{E} \Big[ \max_k   \sup_{t  \in [t_k,t_{k+1}]}
    \Big|   \epsilon^{-1/2}    \int_{t_k}^t   Y_s^\epsilon   \mathrm{d} s   \Big|^2   \Big]
      \le   {}&      \left(\sum_k         \mathbb{E}  \Big  [   \sup_{t\in[t_k,t_{k+1}]}
 \Big|\epsilon^{-1/2}\int_{t_k}^t   Y_s^\epsilon   \mathrm{d} s   \Big|^4   \Big  ]                       \right)^{1/2}        \\ 
 \le   {}&   \left(\sum_k C_T   \epsilon^2 \right)^{1/2}   \\
 \le  {}&  C_T   \epsilon^{1/2}.
\end{align*}
Using $\max_k A_k^2\le\sum_kA_k^2$, we obtain
 \begin{align*}
 &   \mathbb{E}  \Big[  \max_k   \Big(
 \int_{t_k}^{t_{k+1}}  \bigl(|F_{\rho_s^\epsilon}(X_s^\epsilon)|+| \mathcal{Q}(q_{s }^\epsilon;X_{s }^\epsilon, Y_{s }^\epsilon)|\bigr)\mathrm{d} s
 \Big)^2  \Big  ] \\
 \leq   {}&  \sum_k   \mathbb{E}\left(\int_{t_k}^{t_{k+1}}\bigl(|F_{\rho_s^\epsilon}(X_s^\epsilon)|+| \mathcal{Q}(q_{s }^\epsilon;X_{s }^\epsilon, Y_{s }^\epsilon)|\bigr)    \mathrm{d} s\right)^2\\
 \leq  {}&C\epsilon   \int_0^T      \mathbb{ E}\bigl(|F_{\rho_s^\epsilon}(X_s^\epsilon)| ^2 +| \mathcal{Q}(q_{s }^\epsilon;X_{s }^\epsilon, Y_{s }^\epsilon)|   ^2   \bigr)  \mathrm{d} s  \\
 \leq   {}& C_T\epsilon.
\end{align*}
Similarly, we obtain
 \begin{align*}
&    \mathbb{E}   \Big[ \max_k\sup_{t\in [t_k,t_{k+1}]}
  \left|\int_{t_k}^t\lambda_{\rho_s}^{-1}  (X_s^\epsilon)\sigma   \mathrm{d} B_s\right|^2  \Big  ] \\
  \leq {} & \Big(\sum_k  \mathbb{E}\sup_{t\in[t_k,t_{k+1}]}
 \left|\int_{t_k}^t\lambda_{\rho_s}^{-1}(X_s^\epsilon)\sigma   \mathrm{d} B_s\right|^4  \Big)^{1/2} \\
   \leq {} & C\Big(\sum_k     \mathbb{E} 
\Big( \int_{t_k}^{t_{k+1}}\left|\lambda_{\rho_s}^{-1} (X_s  ^\epsilon)\sigma\right|^2 \mathrm{d} s  \Big)^{2} \Big)^{1/2} \\
 \leq   {}& C_T \epsilon^{1/2}.
\end{align*}
 Then   we  get  
$$  \mathbb{E}   \Big[ \sup_{t\le T}  | \epsilon^{ 1/2}      \Xi   (q_{t }^\epsilon;X_{t }^\epsilon,Y_{t }^\epsilon) -   \epsilon^{ 1/2}     \Xi   (q_{t_k}^\epsilon;X_{t_k}^\epsilon,Y_{t_k}^\epsilon)| ^{2}   \Big  ]\leq    C_T \epsilon^{1/2}.
 $$
 Hence  
\begin{align*}
 \mathbb{E}  \Big[    \sup_{t\le T}| \epsilon^{ 1/2}    \Xi(q_t^\epsilon;X_t^\epsilon,Y_t^\epsilon) |^2   \Big]   \leq  {} &
   \mathbb{E}   \Big[ \sup_{t\le T}  | \epsilon^{ 1/2}    \Xi   (q_{t }^\epsilon;X_{t }^\epsilon,Y_{t }^\epsilon) -   \epsilon^{ 1/2}    \Xi   (q_{0}^\epsilon;X_{0}^\epsilon,Y_{0}^\epsilon)| ^{2}   \Big  ]   +   \mathbb{E}   \Big[   |  \epsilon^{ 1/2}   \Xi   (q_{0}^\epsilon;X_{0}^\epsilon,Y_{0}^\epsilon)| ^{2}   \Big  ]  \\
  \leq {} & C_T \epsilon^{1/2}.   
 \end{align*}
 For any $\xi  >0$,
$$
\mathbb{P}   \left(\sup_{t\le T}\big|  \epsilon^{ 1/2}    \Xi (q_{t }^\epsilon;X_{t }^\epsilon,Y_{t }^\epsilon) \big|>\xi  \right)
\le \frac{\mathbb{E}\big[\sup_{t\le T}|  \epsilon^{ 1/2}     \Xi (q_{t }^\epsilon;X_{t }^\epsilon,Y_{t }^\epsilon)  |^2\big]}{\xi^2}
\le \frac{C_T \epsilon^{1/2}}{\xi^2} \to 0,   \qquad    \epsilon  \rightarrow  0. 
$$
So $\{\epsilon^{ 1/2}   \Xi (q^\epsilon;X^\epsilon,Y^\epsilon)\}_{0< \epsilon   \le1}$ is tight in $\mathcal{C}([0,T];\mathbb{R}^d)$. Since $X^\epsilon=Z^\epsilon-  \epsilon^{ 1/2}   \Xi (q^\epsilon;X^\epsilon,Y^\epsilon)$ and $\{Z^\epsilon\}_{0< \epsilon   \le1}$ is tight, $\{X^\epsilon\}_{0< \epsilon   \le1}$ is also tight.
Moreover,  we  have 
$$
 \sup_{0<\epsilon\le1}  \mathbb{E}  \Big[\sup_{t\le T}|X_t^\epsilon|^2\Big]
 \le 2\sup_{0<\epsilon\le1}    \mathbb{E}\Big[\sup_{t\le T}|Z_t^\epsilon|^2\Big]
 +2\sup_{0<\epsilon\le1}   \mathbb{E}\Big [\sup_{t\le T}|  \epsilon^{ 1/2}    \Xi (q_{t }^\epsilon;X_{t }^\epsilon,Y_{t }^\epsilon)|^2\Big  ]<   C_{T}.
$$
\endproof

Next, we provide an estimate for the second moment of the increments.

\begin{lemma}\label{lem:X-increment}
Under \textup{\textbf{(A1)}-\textbf{(A2)} }
\begin{equation}\label{eq:X-increment}
\mathbb{E}  \Big[     |X_t^ \epsilon-X_s^ \epsilon|^2\Big]    \le C_T  |t-s|   , \qquad 0   \le s\le t\le T .
\end{equation}
\end{lemma}

\proof
Let  $0\le s\le t\le T$  and  set  $h:= t-s$.  
For   $ h  \leq  \epsilon   ,$  by   $
 X_t^\epsilon=  X_0^\epsilon+  \epsilon^{-1/2}  \int_0^t   Y_s^\epsilon\,\mathrm{d}  s
$
and H\"older's inequality,    we have 
\begin{align*}
\mathbb{E}  \Big[     |X_t^ \epsilon-X_s^ \epsilon|^2\Big]    \le  {}&  \epsilon^{-1 }   h   \int_s^t   \mathbb{E}  \Big[     |  Y_r^\epsilon|^2\Big]     \,\mathrm{d}  r   \\
\le  {}&    C_{T}  \epsilon^{- 1  }  h^{2} \\
\le  {}&    C_{T}    h. 
\end{align*}
 For   $ h   >  \epsilon   ,$    analogous to the argument in Proposition   \ref{prop:tightness-complete}, we obtain
$$
 \mathbb{ E}|Z_t^\epsilon-Z_s^\epsilon|^2\le C_T(|t-s|+|t-s|^2)\le C_T|t-s|.
$$
  Using
$$
X_t^\epsilon-X_s^\epsilon
=Z_t^\epsilon-Z_s^\epsilon-   \epsilon^{ 1/2}    \Xi(q_t^\epsilon;X_t^\epsilon,Y_t^\epsilon)     + \epsilon^{ 1/2}    \Xi(q_s^\epsilon;X_s^\epsilon,Y_s^\epsilon) 
$$
and 
$$
 \sup_{r\le T}   \mathbb{E}   \Big[   |\epsilon^{ 1/2}    \Xi(q_r^\epsilon;X_r^\epsilon,Y_r^\epsilon) |^2   \Big] 
 \le \epsilon\sup_{r\le T}(  \mathbb{E} \Big[   |\Xi(q_r^\epsilon;X_r^\epsilon,Y_r^\epsilon) |^4     \Big   ]  )^{1/2}
 \le C_T  \epsilon,
$$
we obtain 
$$
\mathbb{E}  \Big[     |X_t^ \epsilon-X_s^ \epsilon|^2\Big]    \le C_T(|t-s|+ \epsilon) \le C_T   h .
$$
\endproof

\begin{remark}\label{rem:2p-moment-closure}
 Assume that       $ \sup \limits_{0<  \epsilon  \le1}   \mathbb{E}   \bigl[|X_0^\epsilon |^{2p}+|Y_0^\epsilon |^{2p}\bigr]<\infty$   for   $p\geq 2 $.
 Under \textup{\textbf{(A1)}-\textbf{(A2)}}, for each $T>0$ 
\begin{equation*} 
 \sup_{0<  \epsilon  \le1}\sup_{t\le T}
 \mathbb{E}[|Y_t^\epsilon |^{2p}]
 +\sup_{0<\epsilon \le1}  \mathbb{ E}\Big[\sup_{t\le T}|Z_t^\epsilon |^{p}   \Big]
 +\sup_{0<\epsilon \le1}  \mathbb{ E}\Big[\sup_{t\le T}|X_t^\epsilon |^{p}   \Big]
 \le C_T,
\end{equation*}
\begin{equation*} 
 \sup_{\epsilon,t\le T}    \mathbb{E}\Big[|\Xi(q_t^\epsilon;X_t^\epsilon,Y_t^\epsilon)|^{2p} \Big]    \le C_T,       \qquad    \sup_{\epsilon,t\le T}    \mathbb{E}\Big[   \mathcal{Q}(q_t^\epsilon;X_t^\epsilon,Y_t^\epsilon)|^{p} \Big]    \le C_T,  
\end{equation*}
\begin{equation*} 
 \mathbb{E}  \Big[    \sup_{t\le T}|\epsilon^{ 1/2}     \Xi(q_t^\epsilon;X_t^\epsilon,Y_t^\epsilon) |^{p}   \Big] \leq  C_T \epsilon^{(p-1)/2}    
\end{equation*}
and 
\begin{equation*} 
\mathbb{E}  \Big[     |X_t^ \epsilon-X_s^ \epsilon|^{p}\Big]    \le C_T  |t-s|^{p /2}   , \qquad 0   \le s\le t\le T .
\end{equation*}

\end{remark}

\section{Smoluchowski--Kramers approximation}\label{Smoluchowski-Kramers-approximation}

In this section, we prove  Theorems \ref{thm:sequence-identification}  and  \ref{thm:strong-convergence}.
In Section \ref{frozen-C-S}, we mainly establish several results on the fast time scale, including the decomposition of the quadratic functional $\mathcal{Q}$ into the equilibrium average $\overline{\mathcal{Q}}$, local quadratic fluctuations, and exponential transients, and also establish the Poisson--Lyapunov equation. In Section \ref{averaging-lemma}, we   establish the averaging lemma. In Section \ref{proof-main-results}, we complete the proofs of Theorems \ref{thm:sequence-identification}  and  \ref{thm:strong-convergence}.

\subsection{ The  frozen stochastic Cucker--Smale equation}  \label{frozen-C-S}

Take the block length $\Delta_\epsilon$, and set
$$
K_\epsilon=\left\lceil\frac{T}{\Delta_\epsilon}\right\rceil,
\qquad t_k=(k\Delta_\epsilon)\wedge T,
\qquad I_k=[t_k,t_{k+1}], \quad  0\le k  <K_\epsilon .
$$
Define
$$
 \mathscr{A}(q;x,y)
 :=-y+\int_{   \mathbb{R}^{2d}    }   \phi(x-x')(y'-y)q(\mathrm{d}  x',\mathrm{d}  y').
$$
For the fast time 
$0\le r\le (t_{k+1}-t_k)/\epsilon     $, define
\begin{equation*} 
  \check{X}_r^{\epsilon,k}:=   X_{t_k+\epsilon r}^\epsilon,
 \qquad
\check{Y}_r^{\epsilon,k}:=Y_{t_k+\epsilon r}^\epsilon,
 \qquad
 W_r^k:=\frac{B_{t_k+\epsilon r}-B_{t_k}}{\sqrt{\epsilon}}.
\end{equation*}
 Set 
$  \check{q}_r^{\epsilon,k}   := \mathrm{Law}(\check{X}_r^{\epsilon,k},\check{Y}_r^{\epsilon,k})$
and let its spatial marginal be $\check{\rho}_r^{\epsilon,k}$.  Then  we have 
\begin{equation*} 
\left\{
\begin{aligned}
 \mathrm{d}  \check{X}_r^{\epsilon,k}
 &=\sqrt{\epsilon}\,\check{Y}_r^{\epsilon,k}   \mathrm{d}  r,\\
 \mathrm{d}   \check{Y}_r^{\epsilon,k}
 &=\mathscr{A}(\check{q}_r^{\epsilon,k};
       \check{X}_r^{\epsilon,k},\check{Y}_r^{\epsilon,k})\mathrm{d}    r
   -\sqrt{\epsilon}\,
       G_{\check{\rho}_r^{\epsilon,k}}(\check{X}_r^{\epsilon,k})\mathrm{d}    r
   +\sigma  \mathrm{d}  W_r^k,
\end{aligned}
\right.
\end{equation*}
where $W^k$ is a standard Brownian motion.
Note that the damping and noise are both $O(1)$ on the fast time scale; the position moves only $O(\sqrt{\epsilon})$ per fast time unit; and the potential term in the fast equation carries a factor $\sqrt{\epsilon}$. Consequently, the frozen Cucker--Smale  equation freezes the spatial label at $X_{t_k}^{\epsilon}$, freezes the spatial marginal at $\rho_{t_k}^{\epsilon}$, and drops the $O(\sqrt{\epsilon})$ potential term in the fast equation.

Next, we discuss some results for the abstract frozen Cucker--Smale equation on the fast time scale within each block.
Fix $\rho\in   \mathcal{P}_2( \mathbb{R}^d)$ and a probability measure $q_0(\mathrm{d}  x,\mathrm{d}  y)=\rho(\mathrm{d}  x)q_0^x(\mathrm{d}  y)$ whose spatial marginal is $\rho$. On the fast time scale $r$, we keep $X$ fixed with $ \mathrm{Law}(X)=\rho$, and consider the abstract frozen Cucker--Smale equation 
\begin{equation}\label{eq:frozen-fast}
\mathrm{d}  Y_r
=\left[-Y_r+\int_{   \mathbb{R}^{2d}    }   \phi(X-x')\big(y'-Y_r\big)q_r(\mathrm{d}  x',\mathrm{d}  y')\right]\mathrm{d}  r
+\sigma\,\mathrm{d}    W_r,
\qquad q_r= \mathrm{Law}(X,  Y_r).
\end{equation}
Let
$
u_r(x)= \mathbb{E}[Y_r   \mid X=x].
$
Taking conditional expectations with respect to \eqref{eq:frozen-fast} yields
\begin{equation}\label{eq:u-semigroup}
\partial_{r}u_r=-A_\rho u_r,
\qquad u_r=e^{-rA_\rho}u_0,
\end{equation}
where   $A_\rho$ is defined in \eqref{eq:A-operator}.
Define
$
\eta_r :  =Y_r-u_r(X).
$
Subtracting \eqref{eq:frozen-fast} and \eqref{eq:u-semigroup} exactly gives
\begin{equation}\label{eq:centered-OU}
\mathrm{d}   \eta_r=-\lambda_\rho(X)\eta_r\,\mathrm{d}  r+\sigma\,\mathrm{d}  W_r.
\end{equation}

\begin{lemma}\label{lem:mean-decay}
Under \textup{\textbf{(A1)}-\textbf{(A2)}}, for any $p\ge2$, if $u_0\in L^p(\rho)$, then
\begin{equation*}
\|u_r\|_{L^p(\rho)}\le e^{-r}\|u_0\|_{L^p(\rho)}.
\end{equation*}
\end{lemma}

\proof

For $p\ge 2$, write \eqref{eq:u-semigroup} as
$$
\partial_ru_r(x)=-u_r(x)+\int_{   \mathbb{R}^{ d}    }    \phi(x-x')\big(u_r(x')-u_r(x)\big)\rho (\mathrm{d} x').
$$
Multiplying by $|u_r(x)|^{p-2}u_r(x)$ and integrating with respect to $\rho $, by the evenness of $\phi$,
\begin{align*}
&\int_{   \mathbb{R}^{ d}    }  |u_r(x)|^{p-2}u_r(x)\cdot
\int_{   \mathbb{R}^{ d}    }    \phi(x-x')\big(u_r(x')-u_r(x)\big)\rho (\mathrm{d} x')\rho (\mathrm{d} x)\\
= {}&   -\frac{1}{2}  \int _{   \mathbb{R}^{2d}    }   \phi(x-x')
\Big(|u_r(x)|^{p-2}u_r(x)-|u_r(x')|^{p-2}u_r(x')\Big)
\cdot\big(u_r(x)-u_r(x')\big)\rho (\mathrm{d} x)\rho (\mathrm{d} x')\\
\le   {}&0,
\end{align*}
where   we have used that the map $z\mapsto |z|^{p-2}z$ is monotone. Therefore,
$$
\frac{1}{p}\frac{\mathrm{d}}{\mathrm{d} r}\|u_r\|_{L^p(\rho)}^p
\le-\|u_r\|_{L^p(\rho)}^p,
$$
which yields $\|u_r\|_{L^p(\rho)}\le e^{-r}\|u_0\|_{L^p(\rho)}$.

\endproof

Define
\begin{equation}\label{eq:C-rho}
J_\rho(x):=\frac{\sigma ^{2}   \mathbf{I}_{d}}{2\lambda_\rho(x)}.
\end{equation}
It is the unique solution of the Lyapunov equation
\begin{equation}\label{eq:Lyapunov}
\lambda_\rho(x)J_\rho(x)+J_\rho(x)\lambda_\rho(x)=\sigma ^{2}   \mathbf{I}_{d}.
\end{equation}
Let
$$
\pi_\rho^x   =\mathcal{N}(0,J_\rho(x)),
\qquad
\Pi_\rho(\mathrm{d}x,\mathrm{d}y)=\rho(\mathrm{d}x)\pi_\rho^x(\mathrm{d}y).
$$
If
$
S_r(x)= \mathbb{E}[\eta_r\eta_r^{ \top }\mid X=x],
$
then, by \eqref{eq:centered-OU},  we have 
$
\frac{\mathrm{d}}{\mathrm{d}r}S_r(x)=-2\lambda_\rho(x)S_r(x)+  \sigma ^{2}   \mathbf{I}_{d}.
$
Consequently,
\begin{equation}\label{eq:covariance-explicit}
S_r(x)=e^{-2\lambda_\rho(x)r}S_0(x)
+\big(1-e^{-2\lambda_\rho(x)r}\big)J_\rho(x).
\end{equation}
 Define
$$
\mathcal{P}_{2,\rho}(\mathbb{R}^d   \times \mathbb{R}^d  ):=\left\{q \in\mathcal{P}(\mathbb{R}^d \times \mathbb{R}^d):   q(\mathrm{d} x, \mathbb{R}^d  )=  \rho (\mathrm{d} x),  \int_{\mathbb{R}^  {2d}   }|y|^2  \, q(\mathrm{d} x,  \mathrm{d} y)   <\infty\right\}.
$$
For $q,\widetilde q  \in   \mathcal{P}_{2,\rho} $,  from  $q (\mathrm{d}  x,\mathrm{d}  y)=\rho(\mathrm{d}  x)q ^x(\mathrm{d}  y)$  and  $\widetilde{q} (\mathrm{d}  x,\mathrm{d}  y)=\rho(\mathrm{d}  x)\widetilde{q}^x(\mathrm{d}  y)$, we can define the fibered Wasserstein distance  ( see, e.g., \cite[Definition 3.1]{Peszek2023}) 
$$
 \mathbb{W}_{2,\rho}^2(q,\widetilde{q})
:=\int_{   \mathbb{R}^{ d}      }   \mathbb{W}_2^2(q^x,\widetilde{q}^x)\,\rho(\mathrm{d}x).
$$

\begin{proposition} \label{prop:frozen-mixing}
If $\int _{   \mathbb{R}^{2d}    }           |y|^2    q_0(\mathrm{d} x,\mathrm{d} y)<\infty$, then  
\begin{equation}\label{eq:fiber-mixing}
 \mathbb{W}_{2,\rho}^2(q_r,\Pi_\rho)
\le Ce^{-2r}\left(1+\int_{   \mathbb{R}^{2d}    }   |y|^2    q_0(\mathrm{d} x,\mathrm{d} y)\right),
\end{equation}
where $C$ depends   not on $\rho$.
 $\Pi_\rho$ is the unique invariant probability measure of the nonlinear distribution evolution induced by the frozen Cucker--Smale   \eqref{eq:frozen-fast}  with initial distribution $q_0$ in $\mathcal P_{2,\rho}$. 
In addition, if  $\int _{   \mathbb{R}^{2d}    }    |y|^4    q_0(\mathrm{d} x,\mathrm{d} y)  <\infty$, then
\begin{equation}\label{eq:frozen-fourth}
\sup_{r\ge0}  \int_{   \mathbb{R}^{2d}    }    |y|^4   q_r(\mathrm{d}x,\mathrm{d} y)
\le C\left(1+\int_{   \mathbb{R}^{2d}    }    |y|^4    q_0(\mathrm{d} x,\mathrm{d} y)\right).
\end{equation}
\end{proposition}

\proof

Choose $X  \sim\rho$, and conditional on $X=x$, draw $Y_0^x  \sim q_0^x$, set $\eta_0^x=Y_0^x-u_0(x)$. Also let $\xi$ be a standard Gaussian vector conditionally independent, and define $\bar{\eta}_0^x=J_\rho^{1/2}(x)\xi    \sim   \mathcal{N}(0,J_\rho(x))$. 
The solution to the corresponding Ornstein--Uhlenbeck equation is
$$
 \bar{\eta}_r^x=e^{-\lambda_\rho(x)  r}\bar{\eta}_0^x
 +\int_0^r    e^{-\lambda_\rho(x)(r-s)}\sigma\,\mathrm{d} W_s.
$$
It remains centered Gaussian, with covariance
\begin{align*}
 \mathrm{Cov}(\bar{\eta}_r^x)
= {}&e^{-2\lambda_\rho(x) r}J_\rho(x)+
 \int_0^re^{-2\lambda_\rho(x)(r-s)}\sigma ^{2}   \mathbf{I}_{d}\,\mathrm{d} s\\
= {}&e^{-2\lambda_\rho(x) r}J_\rho(x)+
 \frac{1-e^{-2\lambda_\rho(x) r}}{2\lambda_\rho(x)}\sigma ^{2}   \mathbf{I}_{d}\\
 =   {}&   J_\rho(x).
\end{align*}
Thus $\pi_\rho^x    =\mathcal{N}(0,J_\rho(x))$ is invariant for the fiberwise Ornstein--Uhlenbeck process, i.e., $\bar{\eta}_r^x \sim \mathcal{N}(0,J_\rho(x))$.
Driving both processes $\eta_r^x$ and $\bar{\eta}_r^x$ by the same Brownian motion yield
$$
 \eta_r^x-\bar{\eta}_r^x
 =e^{-\lambda_\rho(x)  r}
 \bigl(\eta_0^x-\bar{\eta}_0^x   \bigr).$$
Let
$ Y_r^x   :=u_r(x)+\eta_r^x$    and 
 $\bar{Y}_r^x    :=\bar{\eta}_r^x$.
 Then  
 \begin{align*}
 \mathbb{W}_2  ^2 (q_r^x,\pi_\rho^x)
 &\le   \mathbb{E}\bigl[|Y_r^x-\bar Y_r^x|^2 \bigr]\notag\\
 &\le  2|u_r(x)|^2+
 2e^{-2\lambda_\rho(x)r}
  \mathbb{E}\bigl[|\eta_0^x-\bar\eta_0^x|^2 \bigr].
\end{align*}
Since $\eta_0^x$ and $\bar\eta_0^x$ are conditionally independent and both have conditional mean zero, we have
\begin{align*}
   \mathbb{E}\bigl[|\eta_0^x-\bar\eta_0^x|^2  \bigr] 
 =   {} & \mathbb{E}[|\eta_0^x|^2 ]+ \mathbb{E}[|\bar \eta_0^x|^2 ] \\
  = {}  & \mathbb{E}[|\eta_0^x|^2 ]+\mathrm{tr} (  \mathbb{E}[ \bar{\eta}_{0}^{x}  (\bar{\eta}_{0}^{x})^{\top} ] )\\
 = {}  & \mathbb{E}[|\eta_0^x|^2 ]+\mathrm{tr}  (J_\rho(x)).
\end{align*}
Using Jensen's inequality and $\eta_0^x=Y_0^x-u_0(x)$, we have
\begin{align*}
 \int _{   \mathbb{R}^{ d}    }      |u_0(x)|^2\rho(\mathrm{d} x)
 &\le\int_{   \mathbb{R}^{2d}    }      |y|^2q_0(\mathrm{d} x,\mathrm{d} y),\\
 \int _{   \mathbb{R}^{ d}    }      \mathbb{E}[|\eta_0^x|^2 ]\rho(\mathrm{d} x)
 &=\int_{   \mathbb{R}^{2d}    }      |y|^2q_0(\mathrm{d}x,\mathrm{d}y)
   -\int_{   \mathbb{R}^{ d}    }    |u_0(x)|^2\rho(\mathrm{d}x)\\
 &\le\int_{   \mathbb{R}^{2d}    }      |y|^2q_0(\mathrm{d} x,\mathrm{d}y) 
\end{align*}
 and 
$$
  \mathrm{tr} (J_\rho(x))
 =\frac{ \mathrm{tr}(   \sigma ^{2}    \mathbf{I}_{d}) }{2\lambda_\rho(x)}
 \le  \frac{\mathrm{tr}   (   \sigma ^{2}    \mathbf{I}_{d})  }{2}.
$$
By   Lemma \ref{lem:mean-decay},  we obtain 
\begin{align*}
  \mathbb{W}_{2,\rho}^2 (q_r,\Pi_\rho) =  {} & \int_{   \mathbb{R}^{ d}    }   \mathbb{W}_2^2   (q_r^x,\pi_\rho^x) \rho(\mathrm{d} x)  \\
 \le  {}&  2\int _{   \mathbb{R}^{ d}    }  |u_r(x)|^2   \rho(\mathrm{d} x)+
 2e^{-2\lambda_\rho(x)r}
 \int _{   \mathbb{R}^{ d}    }  ( \mathbb{E}[|\eta_0^x|^2 ]+    \mathrm{tr}(J_\rho(x)) )   \rho(\mathrm{d} x) \\
 \le  {}&    Ce^{-2r}\left(
 1+\int_{   \mathbb{R}^{ 2  d}    }    |y|^2q_0(\mathrm{d} x,\mathrm{d}  y)\right).
\end{align*}

If $q_0=\Pi_\rho$ and since $\pi_\rho^x=\mathcal{N}(0,J_\rho(x))$,  we have $u_0=0$, and hence $u_r=0$. Moreover, $\pi_\rho^x$ is invariant, so $\Pi_\rho$ is unique  invariant measure.

Using Jensen's inequality and  Lemma  \ref{lem:mean-decay},  we  have  
\begin{equation}\label{eq:frozen-u4-detail}
 \sup_{r\ge0}\int_{   \mathbb{R}^{ d}    }    |u_r(x)|^4\rho(\mathrm{d} x)
 \le\int_{   \mathbb{R}^{ d}    }    |u_0(x)|^4\rho(\mathrm{d} x)
 \le\int_{   \mathbb{R}^{2d}    }   |y|^4q_0(\mathrm{d} x,\mathrm{d} y).
\end{equation}
Due to   $|a-b|^4\le8(|a|^4+|b|^4)$,   we  get 
\begin{equation}\label{eq:frozen-eta0-fourth}
  \mathbb{E}\Big[|\eta_0|^4\Big]
 \le   8    \mathbb{E}\Big[|Y_0|^4\Big]   +8    \mathbb{E}\Big[|u_0(X)|^4\Big]
 \le   16\int_{   \mathbb{R}^{2d}    }    |y|^4     q_0(\mathrm{d} x,\mathrm{d} y).
\end{equation}
Define   $\zeta_r^x     :=  \int_0^r    e^{-\lambda_\rho(x)(r-s)}    \sigma   \mathrm{d}  W_s.$
Given $X=x$, the conditional distribution is centered Gaussian with covariance
$$
 \Gamma_r(x)
 =\frac{1-e^{-2\lambda_\rho(x)r}}
 {2\lambda_\rho(x)}  \sigma ^{2}    \mathbf{I}_{d},
\qquad
  \mathrm{tr}(\Gamma_r(x))  \le\frac{1}{2}   \mathrm{tr}(\sigma ^{2}    \mathbf{I}_{d}).
$$
The fourth moment formula for centered Gaussians gives
$$
  \mathbb{E}\Big[|\zeta_r^x|^4 \Big]
 =\bigl( \mathrm{tr}(\Gamma_r(x))\bigr)^2
 +2 \mathrm{tr}\bigl(\Gamma_r^2(x)\bigr)
 \le C .
$$
By  applying  $ \eta _r =e^{-\lambda_\rho(X)  r} \eta _0 
 +\int_0^r    e^{-\lambda_\rho(X)(r-s)}\sigma  \,\mathrm{d} W_s $  and \eqref{eq:frozen-eta0-fourth}, 
we derive
\begin{equation}\label{eq:frozen-fourth-eta}
   \sup_{r\ge0}   \mathbb{E} \Big[|\eta_r|^4  \Big]
 \le C\left(1+\int _{   \mathbb{R}^{2d}    }    |y|^4q_0(\mathrm{d} x,\mathrm{d} y)\right).
\end{equation}
 Finally, using   $Y_r = u_r(X) + \eta_r$   and   \eqref{eq:frozen-u4-detail}, we get 
$$
\sup_{r\ge0} \int_{\mathbb{R}^{2d}} |y|^4 \, q_r(\mathrm{d}x,\mathrm{d}y)
\le C\left(1+\int_{\mathbb{R}^{2d}} |y|^4 \, q_0(\mathrm{d}x,\mathrm{d}y)\right).
$$

\endproof

Set  $\beta_r(x)  :=\int_{   \mathbb{R}  ^{d}}    \nabla\phi(x-x')\cdot u_r(x')\rho(\mathrm{d} x')    $.   
Using the fact that the copies $Y^0,\ldots,Y^n$ are conditionally independent given their spatial labels $x_0,\ldots,x_n$, we have
$$
 \mathbb{E}[Y_r^n(Y_r^i)^{ \top }\mid x_0,\ldots,x_n]
=u_r(x_n)(u_r(x_{i}))^{ \top},\quad i\ne n,
$$
and  
$$
 \mathbb{E}   [Y_r^n(Y_r^n)^{\top  }\mid x_n]=S_r(x_n)+u_r(x_n)    (u_r(x_n))^{\top}.
$$
Hence, \eqref{eq:Q-def} can be written as
\begin{equation}\label{eq:frozen-Q-formula}
\begin{aligned}
 &\mathcal{Q}(q_r;x_0,y_0)\\
={}&y_0  y_0^{\top}g_{0,0, \rho }(x_0)
+y_0  w_{0, \rho}(x_0)\frac{\beta_r(x_0)}{\lambda_\rho(x_0)} +\sum_{n=1}^\infty    \int_{  \mathbb{R}  ^{ nd}}  \Bigg[
u_r(x_{n})    y_0^{\top}    g_{n,0, \rho}
+\sum_{i=1}^{n-1}  u_r(x_{n})    (u_r(x_{i}))^{\top}    g_{n,i, \rho}\\
&+\big(S_r(x_{n})  +u_r(x_{n})      (u_r(x_{n}))^{\top}\big)   g_{n,n, \rho}   +u_r(x_{n})   w_{n, \rho}\sum_{i=0}^n     \frac{\beta_r(x_i)}{\lambda_\rho(x_i)}
\Bigg]\rho^{\otimes n}(\mathrm{d} x_1,  \dots,    \mathrm{d} x_n).
\end{aligned}
\end{equation}
 Under $\Pi_\rho$, we obtain  
\begin{equation}\label{eq:equilibrium-conditional-moments}
  \mathbb{E}^{\Pi _\rho}[Y^i\mid X^i=x_i]=0,
 \qquad
  \mathbb{E}^{\Pi _\rho}[Y^i(Y^i)^{\top}\mid X^i=x_i]=J_\rho(x_i) 
\end{equation}
and 
\begin{equation}\label{eq:equilibrium-two-marks}
  \mathbb{E} ^{\Pi _\rho} \left[Y^n(Y^i)^{ \top }\mid x_0,\ldots,x_n\right]
 =\begin{cases}
 0,&0\le i<n,\\
 J_\rho(x_n),&i=n,
 \end{cases}
\end{equation}
where  $\mathbb{E}^{\Pi_\rho}$ denotes expectation under $\Pi_\rho$. 
With the help of
\eqref{eq:equilibrium-conditional-moments}   and   \eqref{eq:equilibrium-two-marks},
we define 
\begin{align}\label{eq:Qbar}
 \overline{\mathcal{Q}}_\rho(x_0)
:=  {}& \mathbb{E} ^{\Pi _\rho}[Q(  \Pi _\rho  ;  x_{0},    Y^{0} )     \mid   X_{0} =  x_0]  \notag\\
= {}& \int_{\mathbb{R}  ^{n}} \mathcal{ Q}(  \Pi _\rho  ;  x_{0},    y_{0} ) \pi_\rho^{x_{0}} (\mathrm{d}  y_{0})  \notag\\
 = {}&  \sum_{n=0}^\infty\int_{  \mathbb{R}  ^{ nd}}
J_\rho(x_n)g_{n,n,\rho}(x_0,\ldots,x_n)
\rho^{\otimes n}(\mathrm{d} x_1,  \dots,  \mathrm{d} x_n),
\end{align}
where we adopt the convention that $x_0$ is not integrated and the $n=0$ term is interpreted as $J_\rho(x_0)   \nabla R_\rho(x_0)$.   Then   we have 
\begin{equation}\label{eq:Qbar-uniform-bound-early}
 \sup_{\rho\in\mathcal{P}_2( \mathbb{R}^d)}
 \sup_{x\in  \mathbb{R}^d}|\overline{\mathcal{Q}}_\rho(x)|
 \le C\left(1+\sum_{n=1}^\infty r_*^{n-1}\right)
 \le\frac{C}{1-r_*}<\infty.
\end{equation}

Next, we show that for each $1\le \ell\le d$, the $\ell$-th component of $\overline{\mathcal{Q}}_\rho(x)$ can be represented as the unique solution to a contraction-type Fredholm equation and can also be expressed explicitly. 

\begin{proposition} 
\label{prop:Qbar-Lyapunov-resolvent}
For each $1  \le   \ell   \le d$, define \begin{equation*}
 \mathfrak{H}_{\rho,\ell}
 :=\sum_{a=1}^d\left( B_{\rho,a}R_\rho J_{\rho,\ell a}
 +M_{\rho,a}R_\rho J_{\rho,\ell a}\right), 
\end{equation*}  then  
\begin{equation*}
  2\lambda_\rho(x)J_\rho(x)=   \sigma ^{2}    \mathbf{I}_{d}, \qquad
 \overline{\mathcal{Q}}_{\rho,\ell}=-(I-P_\rho)^{-1}\mathfrak{H}_{\rho,\ell}.
\end{equation*}
Equivalently, $ \overline{\mathcal{Q}}$ is the unique   solution to the following equation
\begin{equation*}
 \overline{\mathcal{Q}} _{\rho,\ell}=-\mathfrak{H}_{\rho,\ell}
 +P_\rho   \overline{\mathcal{Q}}  _{\rho,\ell} 
\end{equation*}
where the ordinary position-dependent friction part is exactly \begin{equation}\label{eq:Qbar-local-part}
 \overline{\mathcal{Q}}_{\rho,\ell}^{\mathrm{loc}}(x)
 :=\sum_{a=1}^d   J_{\rho,\ell a}(x)\partial_aR_\rho(x)
 =-\sum_{a=1}^d(M_{\rho,a}R_\rho J_{\rho,\ell a})(x),
\end{equation}
 while the additional nonlocal part generated by the velocity alignment is \begin{equation}\label{eq:Qbar-alignment-part}
 \overline{\mathcal{Q}}_{\rho,\ell}^{  \mathrm{align}}
 := \overline{\mathcal{Q}}_{\rho,\ell}- \overline{\mathcal{Q}}_{\rho,\ell}^{\rm loc}
 =-(I-P_\rho)^{-1}\sum_{a=1}^d
 \left(B_{\rho,a}R_\rho J_{\rho,\ell a}
 +P_\rho M_{\rho,a}R_\rho J_{\rho,\ell a}\right).
\end{equation} 

\end{proposition}

\proof

By   equation \eqref{eq:Lyapunov},   we  have   $2\lambda_\rho(x)J_\rho(x)=   \sigma ^{2}    \mathbf{I}_{d} .$
By \eqref{eq:Qbar}, the \(\ell\)-th component   of  $ \overline{\mathcal{Q}} _{\rho }$ is given by
$$
 \overline{\mathcal{Q}} _{\rho,\ell}
 =\sum_{n=0}^\infty    \sum_{a=1}^d
 \mathscr{D}_{n,n}^a  J_{\rho,\ell a},
$$
where  $ \mathscr{D}_{n,n}^a$  is  defined  in  Lemma  \ref{lem:operator-insertion}.
For $n\ge1$, taking $i=n$ in \eqref{eq:exact-insertion} gives
\begin{equation}\label{eq:Qbar-terminal-insertion}
 \mathscr{D}_{n,n}^a   J_{\rho,\ell a}
 =-P_\rho^{n-1}B_{\rho,a}R_\rho J_{\rho,\ell a}
 -P_\rho^nM_{\rho,a}R_\rho J_{\rho,\ell a}.
\end{equation}
For $n=0$, we only have
$$
 \mathscr{D}_{0,0}^a   J_{\rho,\ell a}
 =-M_{\rho,a}R_\rho J_{\rho,\ell a}
 =J_{\rho,\ell a}\partial_aR_\rho.
$$
As \(\|P_\rho\|_{{L^\infty(\mu_\rho)}   \to {L^\infty(\mu_\rho)} }   \le r_*<1\), the series above is absolutely convergent in \(L^\infty\). Summing \eqref{eq:Qbar-terminal-insertion} over \(n\), we obtain
$$
 \overline{\mathcal{Q}}_{\rho,\ell}
 =-\sum_{m=0}^\infty P_\rho^m
 \sum_{a=1}^d\left( B_{\rho,a}R_\rho J_{\rho,\ell a}
 +M_{\rho,a}R_\rho J_{\rho,\ell a}\right),
$$
that is, $\overline{\mathcal{Q}}_{\rho,\ell}=-(I-P_\rho)^{-1}\mathfrak H_{\rho,\ell}$. 
Multiplying on the left by $I-P_\rho$ yields $\overline{\mathcal{Q}} _{\rho,\ell}=-\mathfrak{H}_{\rho,\ell}
 +P_\rho   \overline{\mathcal{Q}}  _{\rho,\ell}$, while the contraction property simultaneously ensures uniqueness of the bounded solution.
Finally, \eqref{eq:Qbar-local-part} is the \(n=0\) term. Using the identity \((I-P)^{-1}-I=(I-P)^{-1}P\), subtracting the local part from the full expression yields \eqref{eq:Qbar-alignment-part}.

\endproof

 Define
\begin{equation}\label{eq:H-def}
H_\rho(x,\eta)
:=\big(\eta\eta^{ \top}-J_\rho(x)\big)
g_{0,0,\rho}(x).
\end{equation}
Next, we decompose the quadratic functional $\mathcal{Q}$ into the equilibrium average $\overline{\mathcal{Q}}$, local quadratic fluctuations $H_\rho$, and  exponential transient remainder terms $\mathfrak{T}_r$.

\begin{lemma}\label{lem:Q-decomp}
Under \textup{\textbf{(A1)}-\textbf{(A2)}},    there exists some $\mathfrak{T}_r$ such that
\begin{equation}\label{eq:Q-decomp}
 \mathcal{Q}(q_r;X,Y_r)- \overline{\mathcal{Q}}_\rho(X)
=H_\rho(X,\eta_r)+\mathfrak{T}_r 
\end{equation}
and 
\begin{equation}\label{eq:transient-bound}
 \mathbb{E}[|\mathfrak{T}_r|]
\le Ce^{-r}\left(1+\int _{   \mathbb{R} ^{2d}  }|y|^4q_0(\mathrm{d} x,\mathrm{d} y)\right)^{1/2}.
\end{equation}
\end{lemma}

\proof

Set $b_r(x):=\frac{\beta_r(x)}{\lambda_\rho(x)}$   and    $    Y_r  =u_r(X )+\eta_r        $.
Then 
\begin{align*}
Y_r  (Y_r  )^{ \top }g_{0,0 , \rho}
=\eta_r     (\eta_r )^{\top}g_{0,0, \rho}    +\bigl[u_r(X )(\eta_r )^{\top}
 +\eta_r    (u_r(X ))^{\top}+u_r(X )  (u_r(X ))^{\top}\bigr]g_{0,0, \rho} 
\end{align*}
 and    
$$
 S_r(x_n)  +u_r(x_n)   ( u_r(x_n))^{\top}
 =J_\rho(x_n)+(S_r(x_n)-   J_\rho(x_n))+u_r(x_n)  (u_r(x_n))^{\top}  \quad \text{for} \quad    n \ge  1 .
$$
Hence    we obtain  
$$
 \mathcal{Q}(q_r;X ,Y_{r} )- \overline{\mathcal{Q}}_\rho(X )
=H_\rho(X ,\eta_r )+\mathfrak{T}_r, 
$$
where 
$$
H_\rho(X ,\eta_r )
=\big(\eta_r   ( \eta_r )^{ \top}-J_\rho(X )\big)
g_{0,0,\rho}(X )
$$
and  
\begin{equation}\label{eq:transient-exact-decomposition}
\begin{aligned}
\mathfrak{T}_r:={}&
 \bigl[u_r(X )(\eta_r )^{\top}
 +\eta_r    (u_r(X ))^{\top}+u_r(X )    (u_r(X ))^{\top}\bigr]g_{0,0,\rho}
 +Y_r     w_{0,\rho}    b_r(X )\\
&+\sum_{n=1}^\infty  \int _{    \mathbb{R}^{dn}         }              \Bigg[
 u_r(x_n)(Y_r )^{\top}g_{n,0,\rho}
 +\sum_{i=1}^{n-1}  u_r(x_n) ( u_r(x_i))^{\top}g_{n,i,\rho}\\
& 
 +\bigl(S_r(x_n)-J_\rho( x_n )\bigr)g_{n,n,\rho}
 +u_r(x_n)   (u_r(x_n))^{\top}g_{n,n,\rho}
 +u_r(x_n)   w_{n,\rho } \sum_{i=0}^n   b_r(x_i)
 \Bigg]\rho^{\otimes n}(\mathrm{d}x_{1},  \dots,  \mathrm{d}  x_{ n}).
\end{aligned}
\end{equation}

By   Jensen's inequality  and Lemma \ref{lem:mean-decay}, for $p=2,4$, we have
\begin{equation}\label{eq:transient-u-basic}
 \|u_r\|_{L^p(\rho)}
 \le e^{-r}\|u_{0}\|_{L^p(\rho)}
 \le e^{-r}m_4(q_0)^{1/4}  \le e^{-r}(  1  +   m_4(q_0))^{1/4} .
\end{equation}
Hence  
\begin{equation}\label{eq:transient-beta-basic}
 \|b_r\|_\infty\le\|\beta_r\|_\infty
 \le L_\phi\|u_r\|_{L^1(\rho)}
 \le Ce^{-r}(  1  +   m_4(q_0))^{1/4}.
\end{equation}
By  \eqref{eq:frozen-fourth} and \eqref{eq:frozen-fourth-eta},  we have 
\begin{equation}\label{eq:transient-Yeta-basic}
 \|Y_r\|_{L^4(q_r)}+\|\eta_r\|_{L^4(q_r)}
 \le C  (  1  +   m_4(q_0))^{1/4}.
\end{equation}
By   Jensen's inequality,  we get  
$$
 \|S_0\|_{L^2(\rho)}^2
 \le\int_{    \mathbb{R} ^{d }  }    \mathbb{E}[|\eta_0|^4\mid X=x]\rho(\mathrm{d} x)
 =  \mathbb{ E}  \Big[ |\eta_0|^4  \Big]  \le C   (  1  +   m_4(q_0)).
$$
Due  to   $
 S_r(x)-J_\rho(x)
 =e^{-2\lambda_\rho(x)r}(S_0(x)-J_\rho(x)), 
$
  we  get  
  \begin{equation}\label{eq:transient-cov-basic}
 \|S_r-J_\rho\|_{L^2(\rho)}
 \le Ce^{-2r}  (  1  +   m_4(q_0))^{1/2}.
\end{equation} 
Cauchy--Schwarz's inequality and \eqref{eq:transient-u-basic}-\eqref{eq:transient-Yeta-basic}, together with the uniform boundedness of $g_{0,0,  \rho}$ and $w_{0,  \rho}$, yield
\begin{align*}
  \mathbb{E}[|u_r(X )(\eta_r )^{\top}g_{0,0,  \rho}|]
 &\le C\|u_r\|_{L^2(\rho)}  \|\eta_r\|_{L^2(q_r)} 
 \le Ce^{-r}(  1  +   m_4(q_0))^{1/2},\\
  \mathbb{E}  [|u_r(X )    (u_r(X ))^{\top}g_{0,0,  \rho}|]
 &\le C\|u_r\|_{L^2(\rho)}^2
 \le Ce^{-2r}(  1  +   m_4(q_0))^{1/2},\\
  \mathbb{E}[|Y_r   w_{0,  \rho}  b_r(X )|]
 &\le C\|Y_r\|_{L^2(q_r)} \|b_r\|_\infty
 \le Ce^{-r}(  1  +   m_4(q_0))^{1/2}.
\end{align*}
Thus 
$$
 \mathbb{E}[  |\bigl[u_r(X )(\eta_r )^{\top}
 +\eta_r    (u_r(X ))^{\top}+u_r(X )    (u_r(X ))^{\top}\bigr]g_{0,0,\rho}
 +Y_r     w_{0,\rho}    b_r(X )|]\le Ce^{-r}(  1  +   m_4(q_0))^{1/2}.
$$
The $\ell$-th component   of
$\int _{    \mathbb{R}^{dn}         }              
 u_r(x_n)(Y_r )^{\top}g_{n,0,\rho}
  \rho^{\otimes n}(\mathrm{d}x_{1},  \dots,  \mathrm{d}  x_{ n})$
 is given by  $\sum_{a=1}^d  Y_r^{ a}(\mathscr D_{n,0}^au_r^\ell)(X )  .$
From \eqref{eq:D-bound}, Cauchy--Schwarz's inequality, and \eqref{eq:transient-u-basic}-\eqref{eq:transient-Yeta-basic}, we obtain
\begin{equation*}
\begin{aligned}
  &  \bigg\| \int _{    \mathbb{R}^{dn}         }              
 u_r(x_n)(Y_r )^{\top}g_{n,0,\rho}
  \rho^{\otimes n}(\mathrm{d}x_{1},  \dots,  \mathrm{d}  x_{ n})\bigg\|_{L^1(q_r)}\\
 \le   {}  & Cr_*^{n-1}\|Y_r\|_{L^2(q_r)}\|u_r\|_{L^2(\rho)} \\
 \le  {}  &Cr_*^{n-1}e^{-r}(  1  +   m_4(q_0))^{1/2}.
 \end{aligned}
\end{equation*}
Secondly, for each $1\le i<n$, applying the argument in Lemma \ref{lem:kernel-sum} together with \eqref{eq:two-marked-internal},  \eqref{eq:terminal-i-le-n-mark-bound} and  \eqref{eq:transient-u-basic}   yields
\begin{equation*}
\begin{aligned}
 &\left\|\int_{    \mathbb{R}^{dn}         }  u_r(x_n) ( u_r(x_i))^{\top}g_{n,i,\rho}
 \rho^{\otimes n}(\mathrm{d}x_{1},  \dots,  \mathrm{d}  x_{ n}) \right\|_{L^1(\rho)} \\
 \le  {}& Cr_*^{n-1}\|u_r\|_{L^4(\rho)} ^2 \\
 \le {} &  Cr_*^{n-1}e^{-2r}(  1  +   m_4(q_0))^{1/2}.
  \end{aligned}
\end{equation*}
  Applying   \eqref{eq:D-bound}  and  \eqref{eq:transient-u-basic}   yields
\begin{equation*}
\begin{aligned}
&\left\|\int_{    \mathbb{R}^{dn}         }     u_r(x_n)   (u_r(x_n))^{\top}g_{n,n,\rho}
\rho^{\otimes n}(\mathrm{d}x_{1},  \dots,  \mathrm{d}  x_{ n}) \right\|_{L^1(\rho)} \\
 \le  {}&  Cr_*^{n-1}e^{-2r}(  1  +   m_4(q_0))^{1/2}.
  \end{aligned}
\end{equation*}
From   \eqref{eq:D-bound} and  \eqref{eq:transient-cov-basic},  we get  
\begin{equation*}
 \left\|\int_{    \mathbb{R}^{dn}         }   (S_r(x_n)-J_\rho(x_n))g_{n,n, \rho}
 \rho^{\otimes n}(\mathrm{d}x_{1},  \dots,  \mathrm{d}  x_{ n})\right\|_{L^1(\rho)}
 \le Cr_*^{n-1}e^{-2r}(  1  +   m_4(q_0))^{1/2}.
\end{equation*}
Finally, for each $0\le i\le n$, applying the argument in Lemma \ref{lem:kernel-sum} together with \eqref{eq:beta-chain-bound},  \eqref{eq:transient-u-basic} and     \eqref{eq:transient-beta-basic} yields
\begin{equation*}
\begin{aligned}
& \left\|\int _{    \mathbb{R}^{dn}         }   u_r(x_n)   w_{n,\rho }     b_r(x_i) \rho^{\otimes n}(\mathrm{d}x_{1},  \dots,  \mathrm{d}  x_{ n})
 \right\|_{L^1(\rho)}  \\
 \le  {}&Cr_*^n\|b_r\|_\infty\|u_r\|_{L^2(\rho)}\\
 \le  {}& Cr_*^n    e^{-2r}    (  1  +   m_4(q_0))^{1/2}.
   \end{aligned}
\end{equation*}
Hence, summing all the estimates above and using    
$
\sum_{n=1}^\infty(n+1)r_*^{n-1} = \frac{2-r_*}{(1-r_*)^2}<\infty,
$
we obtain
\begin{align*}
  \mathbb{E}[|\mathfrak{T}_r|]  
 \le {} & C\sum_{n=1}^\infty   (n+1)r_*^{n-1}e^{-r} (  1  +   m_4(q_0))^{1/2}\\
 \le  {} &Ce^{-r}\left(1+\int_{\mathbb{R} ^{2}}   |y|^4q_0(\mathrm{d} x,\mathrm{d}  y)\right)^{1/2}.
\end{align*}

\endproof

For fixed $(\rho,x)$, the generator of the centered Ornstein--Uhlenbeck equation \eqref{eq:centered-OU} is
\begin{equation*}
(\mathscr{L}_{\rho,x}f)(\eta)
=-\lambda_\rho(x)\eta\cdot\nabla_\eta f(\eta)
+\frac{\sigma ^{2}}{2}  \Delta_\eta   f(\eta).
\end{equation*}
By a direct computation of the It\^{o} generator,
$$
\mathscr{L}_{\rho,x}(\eta\eta^{ \top })
=-2\lambda_\rho(x)\eta\eta^{ \top }+   \sigma ^{2}   \mathbf{I}_{d},
$$
which   implies
$$
\mathscr{L}_{\rho,x}(\eta\eta^{ \top}-J_\rho   (x))
=-2\lambda_\rho(x)(\eta\eta^{ \top}-J_\rho (x)).
$$
Let  
\begin{equation*}
\Psi_\rho(x,\eta)
:=-\frac1{2\lambda_\rho(x)}
\big(\eta   \eta^{\top}-J_\rho(x)\big)g_{0,0,\rho}(x).
\end{equation*}
Then   
\begin{equation*}
\mathscr{L}_{\rho,x}\Psi_\rho(x,\eta)=H_\rho(x,\eta).
\end{equation*}
There exists a constant $C$, independent of $\rho$ and $x$, such that 
\begin{equation}\label{eq:Psi-bounds}
 |\Psi_\rho(x,\eta) |   \le C(1+ |\eta  |^2),
\qquad
\|\nabla_\eta\Psi_\rho(x,\eta)\|_{F}   \le C|\eta|,
\qquad
\|\nabla_\eta^2\Psi_\rho(x,\eta)\|_{F}   \le C,
\end{equation}
where  $\|\cdot\|_{F}$ denotes the Frobenius norm.
If the left endpoint of the frozen interval in real time is $t_k$, then the centered process satisfies
\[
\mathrm{d}\eta_s=-\epsilon^{-1}\lambda_\rho(X)\eta_s\,\mathrm{d}s
+\epsilon^{-1/2}\sigma\,\mathrm{d}B_s.
\]
Applying It\^{o}'s formula to $\Psi_\rho(X,\eta_s)$, we obtain
\begin{align*}
\mathrm{d}   \Psi_\rho(X,\eta_s)
={}&-\epsilon^{-1}\lambda_\rho(X)\eta_s\cdot
       \nabla_\eta\Psi_\rho(X,\eta_s)\mathrm{d} s\\
 &+\frac{\sigma ^{2}   }{2\epsilon} \Delta_\eta^2\Psi_\rho(X,\eta_s)\mathrm{d} s
 +\epsilon^{-1/2}\nabla_\eta\Psi_\rho(X,\eta_s)
       \sigma\mathrm{d}B_s\\
={}&\epsilon^{-1}
   (\mathscr{L}_{\rho,X}\Psi_\rho)(X,\eta_s)\mathrm{d} s
 +\epsilon^{-1/2}\nabla_\eta\Psi_\rho(X,\eta_s)
       \sigma  \mathrm{d} B_s\\
={}&\epsilon^{-1}H_\rho(X,\eta_s)\mathrm{d} s
 +\epsilon^{-1/2}\nabla_\eta\Psi_\rho(X,\eta_s)
       \sigma\mathrm{d} B_s.
\end{align*}
For any $t_k\le a<b$,   integrating from $a$ to $b$  yields
\begin{equation}\label{eq:block-Poisson}
\int_a^b    H_\rho(X,\eta_s)\,\mathrm{d}  s
=\epsilon\big[\Psi_\rho(X,\eta_b)-\Psi_\rho(X,\eta_a)\big]
-\sqrt{\epsilon}   \int_a^b
\nabla_\eta   \Psi_\rho(X,\eta_s)\sigma\,\mathrm{d}  B_s.
\end{equation}

\subsection{ The averaging lemma}\label{averaging-lemma}
This subsection is mainly devoted to presenting the averaging lemma, which will be used in the proof of the main theorems.

Take the block length $\Delta_\epsilon$, and set
$$
K_\epsilon=\left\lceil\frac{T}{\Delta_\epsilon}\right\rceil,
\qquad t_k=(k\Delta_\epsilon)\wedge T,
\qquad I_k=[t_k,t_{k+1}], \quad t_{K_\epsilon}=T,  \quad     0\le k  <K_\epsilon .
$$
On $I_k$, define the frozen process
$$
\widehat{X}_s^{\epsilon,k}=X_{t_k}^\epsilon,
\qquad \widehat{Y}_{t_k}^{\epsilon,k}=Y_{t_k}^\epsilon,
$$
and $\widehat{Y}_s^{\epsilon,k}$ satisfies the following equation
\begin{equation}\label{eq:frozen-block}
\mathrm{d}    \widehat{Y}_s^{\epsilon,k}
=\epsilon^{-1}\left[-\widehat{Y}_s^{\epsilon,k}
+\int _{   \mathbb{R}   ^{2d}}  \phi(\widehat X_s^{\epsilon,k}-x')
(y'-\widehat Y_s^{\epsilon,k})\widehat q_s^{\epsilon,k}(\mathrm{d}  x',\mathrm{d}  y')\right]\mathrm{d}   s
+\epsilon^{-1/2}\sigma\,\mathrm{d}   B_s,
\end{equation}
where $\widehat{q}_s^{\epsilon,k} =  \mathrm{Law}(\widehat{X}_s^{\epsilon,k},\widehat{Y}_s^{\epsilon,k})$, and its spatial marginal is  $\rho_{t_k}^\epsilon$.
Using a similar argument as in Proposition \ref{prop:moment-closure}, we obtain
\begin{equation}\label{eq:frozen-Y4}
\sup_{0<\epsilon\le1}\sup_{0\le k<K_\epsilon}
\sup_{s\in I_k} \mathbb{E}   \Big[   |\widehat{Y}_s^{\epsilon,k}|^4  \Big]  \le C_T.
\end{equation}

Now we present the averaging lemma.

\begin{lemma}\label{lem:core-averaging}
Under \textup{\textbf{(A1)}-\textbf{(A2)}},  for any block length satisfying
\begin{equation}\label{eq:block-scale}
\Delta_\epsilon  \rightarrow0,
\qquad
\frac{\epsilon}{\Delta_\epsilon}  \rightarrow0,  \quad  \text{as}     \quad     \epsilon  \to  0,
\end{equation} 
\begin{equation}\label{eq:ucp}
 \mathbb{E} \Big[  \sup_{0\le t\le T}|\mathfrak{R}_t^\epsilon|   \Big]   
\rightarrow0,  \quad  \text{as}     \quad     \epsilon  \to  0,
\end{equation}
where    
\begin{equation}\label{eq:R-eps}
\mathfrak{R}_t^\epsilon
:= \int_0^t\left[ \mathcal{Q}(q_s^\epsilon;X_s^\epsilon,Y_s^\epsilon)
-\overline{\mathcal{Q}}_{\rho_s^\epsilon}(X_s^\epsilon)
\right] \mathrm{d} s.
\end{equation}
\end{lemma}

\proof

The solution of \eqref{eq:frozen-block} on $I_k:=[t_k,t_{k+1}]$ is denoted by
$\widehat{Y}^{\epsilon,k}$, and its law is denoted by $\widehat q_s^{\epsilon,k}$. In section \ref{frozen-C-S}, we discuss some results on the fast time scale within each block. The abstract fast-time frozen flow for the $k$-th block is written as
$$
 (X^{(k)},Y_r^{(k)}),\qquad
 q_r^{(k)}= \mathrm{Law}(X^{(k)},Y_r^{(k)}),\qquad
 0\le r    \le\frac{t_{k+1}-t_k}{\epsilon}, 
$$
where
$
 (X^{(k)},Y_0^{(k)})=(X_{t_k}^{\epsilon},Y_{t_k}^{\epsilon}),$
$q_0^{(k)}=q_{t_k}^{\epsilon},$
$ \mathrm{Law}(X^{(k)})=\rho_{t_k}^{\epsilon}.
$
Through the relation
\begin{align}\label{eq:core-two-clock-dictionary}
 r_k(s):=\frac{s-t_k}{\epsilon},\qquad
 \widehat{Y}_s^{\epsilon,k}=Y_{r_k(s)}^{(k)},\qquad
 \widehat{q}_s^{\epsilon,k}=q_{r_k(s)}^{(k)},\quad s\in I_k,
\end{align}
we   define, from
$ u_r^{(k)}(x):= \mathbb{E}[Y_r^{(k)}\mid X^{(k)}=x]$  and 
 $\eta_r^{(k)}:=Y_r^{(k)}-u_r^{(k)}(X^{(k)})$ in section \ref{frozen-C-S},
the quantities in actual time 
\begin{equation}\label{eq:frozen-block-centered-definition}
 \widehat u_s^{\epsilon,k}(x)
 :=u_{r_k(s)}^{(k)}(x)
 = \mathbb{E} \Big[\widehat{Y}_s^{\epsilon,k}    \mid                 X_{t_k}^{\epsilon}=x\Big],
 \qquad
 \widehat{\eta}_s^{\epsilon,k}
 :=\eta_{r_k(s)}^{(k)}
 =\widehat{Y}_s^{\epsilon,k}  -\widehat{u}_s^{\epsilon,k}(X_{t_k}^{\epsilon}).
\end{equation}
Applying Lemma \ref{lem:Q-decomp} to the abstract fast flow, we obtain
\begin{equation}\label{eq:abstract-block-Q-decomposition}
\begin{aligned}
  \mathcal{Q}(q_r^{(k)};X^{(k)},Y_r^{(k)})
 -  \overline{\mathcal{Q}}_{\rho_{t_k}^{\epsilon}}(X^{(k)})
 =H_{\rho_{t_k}^{\epsilon}}(X^{(k)},\eta_r^{(k)})
  +\mathfrak T_r^{(k)}.
\end{aligned}
\end{equation}
Define
$
 \widehat{\mathfrak{T}}_s^{\epsilon,k}
 :=\mathfrak{T}_{r_k(s)}^{(k)}
 =\mathfrak{T}_{(s-t_k)/\epsilon}^{(k)},  s\in I_k.
$
From \eqref{eq:core-two-clock-dictionary} and \eqref{eq:abstract-block-Q-decomposition}, we get
\begin{equation}\label{eq:actual-block-Q-decomposition}
 \mathcal{Q}(\widehat{q}_s^{\epsilon,k};
X_{t_k}^{\epsilon},\widehat{Y}_s^{\epsilon,k})  - \overline{\mathcal{Q}}_{\rho_{t_k}^{\epsilon}}(X_{t_k}^{\epsilon}) 
 =H_{\rho_{t_k}^{\epsilon}}
   (X_{t_k}^{\epsilon},\widehat{\eta}_s^{\epsilon,k})
  +\widehat{\mathfrak{T}}_s^{\epsilon,k}.
\end{equation}
We decompose \eqref{eq:R-eps} as
\begin{equation}\label{eq:R-decompose}
\mathfrak{R}_t^\epsilon   = :  I_t^\epsilon+II_t^\epsilon+III_t^\epsilon,
\end{equation}
where  
\begin{align*}
I_t^\epsilon
&=\sum_{k=0}^{K_\epsilon-1}     \int_{I_k\cap[0,t]}
\big(\mathcal{Q}(q_s^\epsilon;X_s^\epsilon,Y_s^\epsilon)
-\mathcal{Q}(\widehat{q}_s^{\epsilon,k};
X_{t_k}^{\epsilon},\widehat{Y}_s^{\epsilon,k})   \mathrm{d} s,\\
II_t^\epsilon
&=\sum_{k=0}^{K_\epsilon-1}    \int_{I_k\cap[0,t]}
\big(\mathcal{Q}(\widehat{q}_s^{\epsilon,k};
X_{t_k}^{\epsilon},\widehat{Y}_s^{\epsilon,k})   -    \overline{\mathcal{Q}}_{\rho_{t_k}^\epsilon}(X_{t_k}^\epsilon)    \big)   \mathrm{d} s,\\
III_t^\epsilon
&=\sum_{k=0}^{K_\epsilon-1}    \int_{I_k\cap[0,t]}
\big(   \overline{\mathcal{Q}}_{\rho_{t_k}^\epsilon}(X_{t_k}^\epsilon)    -
 \overline{\mathcal{Q}}_{\rho_s^\epsilon}(X_s^\epsilon)\big)   \mathrm{d} s.
\end{align*}

By Lemmas \ref{lem:X-increment} and \ref{lem:freeze-coupling}, we have
$$
\sup_{0\le k<K_\epsilon}\sup_{s\in I_k}
\left(
 \mathbb{E} \Big[|X_s^\epsilon-X_{t_k}^\epsilon|^2  \Big]
+\mathbb{E} \Big[|Y_s^\epsilon-\widehat{Y}_s^{\epsilon, k}|^2\Big]
\right)  \rightarrow0,  \quad  \text{as}     \quad     \epsilon  \to  0.
$$
By Proposition \ref{prop:moment-closure} and  \eqref{eq:frozen-Y4}, the velocities possess uniform fourth moments, while the positions possess uniform second moments.
If the supremum over $0\le k<K_\epsilon$ and $s\in I_k$ of 
$ \mathbb{E}\Big[| \mathcal{Q}(q_s^\epsilon;X_s^\epsilon,Y_s^\epsilon)-\mathcal{Q}(\widehat{q}_s^{\epsilon,k};
X_{t_k}^{\epsilon},\widehat{Y}_s^{\epsilon,k}) |\Big]$ 
does not tend to $0$, then there exist $\xi>0$, a sequence $\epsilon_j\downarrow0$, 
and indices $k_j$ with times $s_j\in I_{k_j}$   such that the corresponding expectation is greater than $\xi$.
However,   this family of couplings satisfies all the assumptions of Lemma \ref{lem:Q-stability}, leading to a contradiction. Therefore,
$$
\sup_{0\le k<K_\epsilon}    \sup_{s\in I_k}
\mathbb{E}\Big[| \mathcal{Q}(q_s^\epsilon;X_s^\epsilon,Y_s^\epsilon)-\mathcal{Q}(\widehat{q}_s^{\epsilon,k};
X_{t_k}^{\epsilon},\widehat{Y}_s^{\epsilon,k}) |\Big] \rightarrow 0,  \quad  \text{as}     \quad     \epsilon  \to  0.
$$
 Thus we have
\begin{align}\label{eq:I-vanish}
\mathbb{E}\Big[   \sup_{t\le T}|I_t^\epsilon|\Big]
\le  {}& \sum_{k=0}^{K_\epsilon-1}  \int_{I_k }
\mathbb{E}\Big[|\mathcal{Q}(q_s^\epsilon;X_s^\epsilon,Y_s^\epsilon)-\mathcal{Q}(\widehat{q}_s^{\epsilon,k};
X_{t_k}^{\epsilon},\widehat{Y}_s^{\epsilon,k}) |   \Big]  \mathrm{ d} s    \notag\\
\le  {}&   T   \sup_{0\le k<K_\epsilon}    \sup_{s\in I_k}
\mathbb{E}\Big[| \mathcal{Q}(q_s^\epsilon;X_s^\epsilon,Y_s^\epsilon)-\mathcal{Q}(\widehat{q}_s^{\epsilon,k};
X_{t_k}^{\epsilon},\widehat{Y}_s^{\epsilon,k}) |\Big] \notag\\
&   \rightarrow0,  \quad  \text{as}     \quad     \epsilon  \to  0.
\end{align}
By Lemma  \ref{lem:X-increment}   
$$
\sup_{0\le k<K_\epsilon}    \sup_{s\in I_k}  \mathbb{E}\Big[|X_s^\epsilon-X_{t_k}^\epsilon|^2\Big]
\le C_T(\Delta_\epsilon ) 
 \rightarrow0,  \quad  \text{as}     \quad     \epsilon  \to  0.
$$
Similarly, by Lemma \ref{lem:Q-stability} 
$$
\sup_{0\le k<K_\epsilon}    \sup_{s\in I_k}
\mathbb{E}\Big[  | \overline{\mathcal{Q}}_{\rho_{t_k}^\epsilon}(X_{t_k}^\epsilon)    -
 \overline{\mathcal{Q}}_{\rho_s^\epsilon}(X_s^\epsilon)|\Big]   \rightarrow0,  \quad  \text{as}     \quad     \epsilon  \to  0,
$$
which  implies   
\begin{equation}\label{eq:III-vanish}
\mathbb{E}\Big[   \sup_{t\le T}|III_t^\epsilon| \Big]    \rightarrow0,  \quad  \text{as}     \quad     \epsilon  \to  0.
\end{equation}
From  \eqref{eq:actual-block-Q-decomposition}, we have 
$$
 \mathcal{Q}(\widehat{q}_s^{\epsilon,k};
X_{t_k}^{\epsilon},\widehat{Y}_s^{\epsilon,k})  - \overline{\mathcal{Q}}_{\rho_{t_k}^{\epsilon}}(X_{t_k}^{\epsilon}) 
 =H_{\rho_{t_k}^{\epsilon}}
   (X_{t_k}^{\epsilon},\widehat{\eta}_s^{\epsilon,k})
  +\widehat{\mathfrak{T}}_s^{\epsilon,k}.
$$
Lemma \ref{lem:Q-decomp} gives, for the $k$-th abstract fast flow,
$$
 \mathbb{E} \Big[     |\mathfrak{T}_r^{(k)}|\Big] 
 \le Ce^{-r}\bigl(1+m_4(q_0^{(k)})\bigr)^{1/2},
 \qquad r\ge0.
$$
Since $q_0^{(k)}=q_{t_k}^\epsilon$, for $s\in I_k$ we have
\begin{align}
 \mathbb{E} \Big[   |\widehat{\mathfrak{T}}_s^{\epsilon,k}|\Big] 
=  {}& \mathbb{E} \Big[  |
 \mathfrak T_{(s-t_k)/\epsilon}^{(k)} |\Big]   \notag\\
\le  {}&
 C\exp \left(-\frac{s-t_k}{\epsilon}\right)
 \bigl(1+m_4(q_{t_k}^\epsilon)\bigr)^{1/2}   \notag\\
 \le   {}&  C_Te^{-(s-t_k)/\epsilon},
\label{eq:block-transient-time-change-exact}
\end{align}
where the last line follows from Proposition \ref{prop:moment-closure}.
Then,  by \eqref{eq:block-transient-time-change-exact},  we  have  
\begin{equation}\label{eq:transient-sum}
\begin{aligned}
& \mathbb{E} \bigg[\sup_{t\le T}
 \bigg|\sum_{k=0}^{K_\epsilon-1}\int_{I_k\cap[0,t]}
 \widehat{\mathfrak{T}}_s^{\epsilon,k} \mathrm{d} s \bigg|\bigg]  \notag \\
\le  {}& C_T\sum_{k=0}^{K_\epsilon-1}
 \int_{t_k}^{t_{k+1}}e^{-(s-t_k)/\epsilon}   \mathrm{d} s \notag\\
\le   {} & C_T K_\epsilon  \epsilon\\
 \le   {}& C_T\left(1+\frac T{\Delta_\epsilon}\right)\epsilon  \notag\\
&  \rightarrow0,  \quad  \text{as}     \quad     \epsilon  \to  0.
\end{aligned}
\end{equation}
For $0\le j\le K_\epsilon$, it follows from \eqref{eq:block-Poisson} that
\begin{equation}\label{eq:grid-H-Poisson}
 \sum_{k=0}^{j-1}\int_{t_k}^{t_{k+1}}
 H_{\rho_{t_k}^ \epsilon}
 (X_{t_k}^ \epsilon,\widehat{\eta}_s^{ \epsilon,k}) \mathrm{d} s
 =\epsilon  \sum_{k=0}^{j-1}\left[
 \Psi_{\rho_{t_k}^\epsilon}
 (X_{t_k}^\epsilon,\widehat{\eta}_{t_{k+1}}^{\epsilon,k})
 -\Psi_{\rho_{t_k}^\epsilon}
 (X_{t_k}^\epsilon,\widehat{\eta}_{t_k}^{\epsilon,k})
 \right]-\sqrt{ \epsilon} \mathcal{M}_{t_j}^ \epsilon,
\end{equation}
 where
 $$
\mathcal{M}_t^\epsilon
:=\sum_{k=0}^{K_\epsilon-1}\int_{I_k\cap[0,t]}
\nabla_{ \eta }\Psi_{\rho_{t_k}^\epsilon}
(X_{t_k}^\epsilon,\widehat\eta_s^{\epsilon,k})
\sigma  \mathrm{d} B_s.
$$
 By Jensen's inequality, we obtain
\begin{align*}
 \mathbb{E} \Big[|\widehat{u}_s^{\epsilon,k}(X_{t_k}^\epsilon)|^2\Big] 
\le  \mathbb{E} \Big[|\widehat{Y}_s^{\epsilon,k}|^2\Big]  
\end{align*}
and 
\begin{align}\label{eq:frozen-eta-second-uniform}
 \mathbb{E} \Big[|\widehat{\eta}_s^{\epsilon,k}|^2\Big] 
= \mathbb{E} \Big[|\widehat{Y}_s^{\epsilon,k}|^2\Big] 
  -\mathbb{E} \Big[|\widehat{u}_s^{\epsilon,k}(X_{t_k}^\epsilon)|^2\Big] 
 \le\mathbb{E} \Big[|\widehat{Y}_s^{\epsilon,k}|^2\Big] \le C_T.
\end{align}
Further combining with \eqref{eq:Psi-bounds}, we have
$$
 \sup_{\epsilon,k}\sup_{s\in I_k}
  \mathbb{E} \Big[  |\Psi_{\rho_{t_k}^\epsilon}
 (X_{t_k}^\epsilon,\widehat{\eta}_s^{\epsilon,k})|\Big]    \le C_T.
$$
Hence  we have 
\begin{align}\label{eq:boundary-sum}
 \mathbb{E} \bigg[  \max_{0\le j\le K_\epsilon}   \bigg|\epsilon  \sum_{k=0}^{j-1}\left[
 \Psi_{\rho_{t_k}^\epsilon}
 (X_{t_k}^\epsilon,\widehat{\eta}_{t_{k+1}}^{\epsilon,k})
 -\Psi_{\rho_{t_k}^\epsilon}
 (X_{t_k}^\epsilon,\widehat{\eta}_{t_k}^{\epsilon,k})
 \right]\bigg| \bigg]  
 \le C_TK_\epsilon\epsilon   \rightarrow0,  \quad  \text{as}     \quad     \epsilon  \to  0.
\end{align}
By \eqref{eq:Psi-bounds}, \eqref{eq:frozen-eta-second-uniform}, and the Burkholder-Davis-Gundy inequality, we obtain
\begin{align}
 \mathbb{E} \Big[  \sup_{t\le T}|\sqrt{\epsilon } \mathcal{M}_t^\epsilon |^2\Big]
\le   {}& C \epsilon  \mathbb{E} \Big[    \operatorname{tr}\langle \mathcal{M}^\epsilon  \rangle_T\Big]\notag\\
\le {}& C\epsilon  \sum_{k=0}^{K_\epsilon -1}\int_{I_k}
 \mathbb{E} \Big[  \|
 \nabla_\eta\Psi_{\rho_{t_k}^\epsilon }
 (X_{t_k}^\epsilon ,\widehat{\eta}_s^{\epsilon ,k})\sigma
 \|_{F}^2\Big]  \mathrm{d} s\notag\\
\le {}& C\epsilon \sum_{k=0}^{K_\epsilon -1}\int_{I_k}  \mathbb{E} \Big[ |\widehat{\eta}_s^{\epsilon, k  }|^2 \Big] \mathrm{d} s\notag\\
\le {}&  C_T\epsilon \notag\\
&  \rightarrow0,  \quad  \text{as}     \quad     \epsilon  \to  0.
\label{eq:martingale-sum}
\end{align}
Then  
\begin{align*}
 \mathbb{E} \Big[   \max_{0\le j\le K_\epsilon}|II_{t_j}^\epsilon|\Big]  \leq C_T   \Big(\frac{\epsilon}{\Delta_\epsilon}+\epsilon+\sqrt{\epsilon}  \Big)
 \rightarrow  0,  \quad  \text{as}     \quad     \epsilon  \to  0.
\end{align*}
 Thus  
\begin{equation}\label{eq:R-grid-vanish}
 \mathbb{E} \Big[    \max_{0\le j\le K_\epsilon}   |\mathfrak{R}_{t_j}^\epsilon| \Big]
 \rightarrow  0,  \quad  \text{as}     \quad     \epsilon  \to  0.
\end{equation}

Define 
$$
 \kappa_\epsilon(t):=\max\{t_k:t_k\le t\}.
$$
From  \eqref{eq:R-eps}, \eqref{eq:Qbar-uniform-bound-early}, Cauchy--Schwarz's inequality  and  Proposition \ref{prop:moment-closure},  we get 
\begin{align*}
 \mathbb{E} \Big[    \sup_{t\le T}|\mathfrak{R}_t^\epsilon-\mathfrak{R}_{\kappa_\epsilon(t)}^\epsilon|    \Big]  
\le {}&   \mathbb{E} \Bigg[    \max_k\int_{I_k}
 (|\mathcal{Q}(q_s^\epsilon;X_s^\epsilon,Y_s^\epsilon)|+|\overline{\mathcal{Q}}_{\rho_s^\epsilon}(X_s^\epsilon)| ) \mathrm{d} s  \Bigg]  \\
\le {}&  \mathbb{E}\Bigg[    \max_k  |I_k|^{1/2}
 \left(\int_{I_k}   |\mathcal{Q}(q_s^\epsilon;X_s^\epsilon,Y_s^\epsilon)|^2     \mathrm{d} s\right)^{1/2}\Bigg]   +  C\max_k |I_k|  \\
\le {}& (\Delta_\epsilon )^{1/2}
 \left(\int_{0}^{T}  \mathbb{E}\Big[ |\mathcal{Q}(q_s^\epsilon;X_s^\epsilon,Y_s^\epsilon)|^2\Big]     \mathrm{d} s\right)^{1/2}  +  C  \Delta_\epsilon    \\
\le {}&   C_{T}(\Delta_\epsilon )^{1/2}
   +  C  \Delta_\epsilon    \\
  {}&  \rightarrow  0,  \quad  \text{as}     \quad     \epsilon  \to  0.
\end{align*}
Further combining with \eqref{eq:R-grid-vanish}, we have
$$
  \mathbb{E} \Big[  \sup_{t\le T}|\mathfrak{R}_t^ \epsilon |\Big]
 \le    \mathbb{E} \Big[ \max_{0\le j\le K_ \epsilon }|\mathfrak{R}_{t_j}^ \epsilon |\Big]
 +     \mathbb{E} \Big[ \sup_{t\le T}|\mathfrak{R}_t^ \epsilon    -\mathfrak{R}_{\kappa_ \epsilon (t)}^ \epsilon |\Big] \rightarrow  0,  \quad  \text{as}     \quad     \epsilon  \to  0.
$$

\endproof

\begin{remark}\label{rem:core-L1-averaging}
 Under \textbf{(A1\textquotesingle)}-\textbf{(A2\textquotesingle)},  for any block length satisfying
\begin{equation*} 
\Delta_\epsilon  \rightarrow0,
\qquad
\frac{\epsilon}{\Delta_\epsilon}  \rightarrow0,  \quad  \text{as}     \quad     \epsilon  \to  0,
\end{equation*}
\begin{equation*} 
 \mathbb{E} \Big[  \sup_{0\le t\le T}|\mathfrak{R}_t^\epsilon|   \Big]    \leq C_{T}
   \Big(\frac{\epsilon}{\Delta_\epsilon}    +  (\Delta_\epsilon )^{1/4}  +  \epsilon  ^{1/2}\Big).
\end{equation*}
 Moreover, taking  $\Delta_\epsilon= \epsilon ^{4/5}$ 
\begin{equation*} 
 \mathbb{E} \Big[  \sup_{0\le t\le T}|\mathfrak{R}_t^\epsilon|   \Big]    \leq C_{T}
    \epsilon  ^{1/5} .
\end{equation*} 
\end{remark}

\subsection{ Proof of the main results}\label{proof-main-results}

In this subsection, we prove    our main results.

  \noindent {\emph{\textbf{Proof   of   Theorem \ref{thm:sequence-identification}:} }}
Tightness of $\{\mathrm{ Law}(X^\epsilon)\}$ follows from Proposition~\ref{prop:tightness-complete}. Moreover, from  the definition of $ Z ^\epsilon$,   Lemma
\ref{lem:core-averaging}  and  
\eqref{eq:Xi-ucp-L2},  we  have  
$$
\mathbb{E} \Big[  \sup\limits _{0 \leq t  \leq T} |Z_{t}^\epsilon-X_{t}^\epsilon   |^{2}   \Big] \to 0,\qquad 
\mathbb{E} \Big[ \sup\limits _ {0 \leq t  \leq T} |\mathfrak{R}_{t}^\epsilon| \Big] \to 0,   \quad  \text{as}     \quad     \epsilon  \to  0.
$$
Let  $\tilde{Z}_t^\epsilon:=Z_t^\epsilon-\mathfrak{R}_t^\epsilon$.
Then we deduce
$$
\mathbb{E} \Big[  \sup\limits _{0 \leq t  \leq T}  |\widetilde{Z}_{t}^\epsilon-   X_{t}^\epsilon | \Big] \to 0,   \quad  \text{as}     \quad     \epsilon  \to  0.
$$
Hence the joint family $\{\mathrm{ Law}(X^\epsilon, B^\epsilon, \tilde{Z}^\epsilon)\}$ is tight in $\mathcal{P}(\mathcal{C}([0,T];\mathbb{R}^d)\times \mathcal{C}([0,T];\mathbb{R}^d)\times \mathcal{C}([0,T];\mathbb{R}^d))$. By Prokhorov's theorem, there exists a subsequence (still denoted by $\epsilon$) such that
$
\mathrm{ Law}(X^\epsilon, B^\epsilon, \tilde{Z}^\epsilon)  \Longrightarrow  \Lambda,
$
for some probability measure $\Lambda$ on the product space. Let $(\tilde{X},\tilde{B},\tilde{Z})$ denote the canonical process under $\Lambda$. Since $\tilde{Z}^\epsilon$ and $\tilde{X}^\epsilon$ are close in $L^1$, the limit measure $\Lambda$ is supported on the set
$
\{(x,w,z): z=x\}.
$
Thus $\tilde{Z}=\tilde{X}$ a.s. under $\Lambda$. In particular, the marginal law of $\tilde{X}$ is the limit $\mathfrak{P}$ mentioned earlier.

By the Skorokhod representation theorem, there exists a new probability space $(\tilde{\Omega},\tilde{\mathcal{ F}},\tilde{\mathbb{P}})$ and, on it, $(\mathcal{C}([0,T];\mathbb{R}^d)\times \mathcal{C}([0,T];\mathbb{R}^d)\times \mathcal{C}([0,T];\mathbb{R}^d))$-valued random variables $(\tilde{X}^\epsilon,\tilde{B}^\epsilon,\tilde{Z}^\epsilon)$ with the same joint law as the original, such that
\begin{equation}\label{eq:joint-Skorokhod-convergence}
(\tilde{X}^\epsilon,\tilde{B}^\epsilon,\tilde{Z}^\epsilon)  \rightarrow   (\tilde{X},\tilde{B},\tilde{Z}) \quad \tilde{\mathbb{P}}\text{-a.s. in the product topology},
\end{equation}
where $(\tilde{X},\tilde{B},\tilde{Z})$ has law $\Lambda$. Since \(\Lambda\) is supported on $z=x$, we have $\tilde{Z}=\tilde{X}$ a.s. Moreover, $\tilde{B}$ is a Brownian motion because it is the almost sure limit of standard Brownian motions. 
From the estimates in Proposition \ref{prop:tightness-complete}, there exists a constant $C_{T}$ such that
$$
\sup_\epsilon \mathbb{\tilde{E}} \Big[  \sup\limits _{0 \leq t  \leq T}  |\tilde{X}_{t}^\epsilon  |^2    \Big] \le C_{T}.
$$
Fatou's lemma yields $\mathbb{\tilde{E}} \Big[  \sup\limits _{0 \leq t  \leq T}   |\tilde X_{t} | ^2 \Big] \le C_{T}$. Consequently, the family $\Big\{  \sup\limits _{0 \leq t  \leq T}   |\tilde{X}_{t}^\epsilon  -\tilde{X}_{t}  |  \Big\}$ is uniformly integrable. Together with the a.s. convergence in \eqref{eq:joint-Skorokhod-convergence}, Vitali's convergence theorem gives
\begin{equation}\label{eq:path-L1}
\mathbb{\tilde{E}} \Big[   \sup\limits _{0 \leq t  \leq T}   |\tilde{X}_{t}^\epsilon  -\tilde{X}_{t}  |\Big]    \rightarrow 0,   \quad  \text{as}     \quad     \epsilon  \to  0.
\end{equation}
As a direct consequence, for every $t\in[0,T]$, the evaluation map $e_t:\omega\mapsto\omega(t)$ is continuous, hence
$$
\tilde{\rho}_t^\epsilon := \mathrm{ Law}(\tilde{X}_t^\epsilon)  \Longrightarrow   \tilde{\rho}_t := \mathrm{ Law}(\tilde{X}_t).
$$
Moreover, \eqref{eq:path-L1} implies  
\begin{equation}\label{eq:path-W1}
\sup_{0\le t\le T} \mathbb{W}_1(\tilde{\rho}_t^\epsilon,\tilde{\rho}_t)
\le  \mathbb{\tilde{E}}  \Big[   \sup\limits _{0 \leq t  \leq T}  |\tilde{X}_t^\epsilon-\tilde{X}_t |\Big]  \rightarrow 0,   \quad  \text{as}     \quad     \epsilon  \to  0. 
\end{equation}
In particular, the initial law is preserved: $\tilde{\rho}_0^\epsilon \Rightarrow \tilde{\rho}_0 = \rho_0$ and $\mathfrak{P}\circ e_0^{-1}=\rho_0$.

From \eqref{eq:Z-decomp} and 
\eqref{eq:R-eps}
we have 
\begin{equation}\label{eq:Ztilde}
\mathrm{ d}  \tilde{Z}_t^\epsilon
 =b_{\tilde{\rho}_t^\epsilon}(\tilde{X}_t^\epsilon)\mathrm{ d}  t
 + \sigma  \lambda_{\tilde{\rho}_t^\epsilon}^{-1}   (\tilde{X}_t^\epsilon)  \mathrm{ d} \tilde{B}_t^\epsilon.
\end{equation}
For any $h\in \mathcal{C}_c^\infty$, It\^o's formula gives a square-integrable   martingale
\begin{align}
 \mathcal{\tilde{M}} _t^{\epsilon,h}
 :={}&h(\tilde{Z}_t^\epsilon)-h(\tilde{Z}_0^\epsilon)\notag\\
 &-\int_0^t\left[
 b_{\tilde{\rho}_s^\epsilon}(\tilde{X}_s^\epsilon)\cdot\nabla h(\tilde{Z}_s^\epsilon)
 +\frac{\sigma^{2}   \lambda_{\tilde{\rho}_s^\epsilon}^{-2}   (\tilde{X}_s^\epsilon)}{2}   \Delta h(\tilde{Z}_s^\epsilon)
 \right]  \mathrm{d} s\notag\\
={}&\int_0^t  \sigma   \lambda_{\tilde{\rho}_s^\epsilon}^{-1}   (\tilde{X}_s^\epsilon)   \nabla h(\tilde{Z}_s^\epsilon)^{ \top}
 \mathrm{d} \tilde{B}_s^\epsilon
\label{eq:prelimit-test-martingale}
\end{align}
 with quadratic variation
$$
\langle \mathcal{\tilde{M}}^{\epsilon,h}\rangle_t
=
\int_0^t \sigma^{2}   \lambda_{\tilde{\rho}_s^\epsilon}^{-2}   (\tilde{X}_s^\epsilon)   |\nabla h(\tilde{Z}_s^\epsilon)|^{2} \mathrm{d}s.
$$
Since the stochastic integrand is uniformly bounded, the Burkholder-Davis-Gundy inequality also gives, for each $p\ge2$,
\begin{equation}\label{eq:prelimit-martingale-pmom}
 \sup_\epsilon  \mathbb{\tilde{E}} \Big[ \sup_{t\le T}| \mathcal{M}_t^{\epsilon,h}|^p\Big]  <\infty.
\end{equation}

Next,  we   shall show that $\mathcal{\tilde{M}}^{\epsilon,h}$ converges in $L^1$ to
 \begin{align}
 \mathcal{\tilde{M}}_t^{ h}
 :={}&h(\tilde{X}_t )-h(\tilde{X}_0 )\notag\\
 &-\int_0^t\left[ b_{\tilde{\rho}_s }(\tilde{X}_s )\cdot\nabla h(\tilde{X}_s )
 +\frac{\sigma^{2}   \lambda_{\tilde{\rho}_s }^{-2}   (\tilde{X}_s )}{2}   \Delta h(\tilde{X}_s )
 \right]  \mathrm{d} s.
\label{eq:limit-test-martingale}
\end{align}
Using \eqref{eq:Qbar-stability},  \eqref{eq:G-stability}  and  \eqref{eq:inverse-2-lambda-stability},  we have 
\begin{align*}
  \mathbb{\tilde{E}}  [|b_{\tilde{\rho}_s^\epsilon}(\tilde{X}_s^\epsilon)-b_{\tilde{\rho}_s}(\tilde{X}_s)|]   &\to0,   \quad  \text{as}     \quad     \epsilon  \to  0, \\
 \mathbb{\tilde{E}}  [   |\sigma^{2}   \lambda_{\tilde{\rho}_s^\epsilon}^{-2}(\tilde{X}_s^\epsilon)  -\sigma^{2}   \lambda_{\tilde{\rho}_s }^{-2}(\tilde{X}_s ) |]   &\to0,   \quad  \text{as}     \quad     \epsilon  \to  0. 
\end{align*}
From Proposition \ref{prop:moment-closure} 
\begin{align*}
 \mathbb{\tilde{E}}  [|  b_{\tilde{\rho}_s^\epsilon}(\tilde{X}_s^\epsilon)-b_{\tilde{\rho}_s}(\tilde{X}_s)| ]\le C_T.
\end{align*}
Therefore, applying the dominated convergence theorem, we obtain
\begin{equation}\label{eq:time-integrated-coeff-convergence}
\begin{aligned}
 &\int_0^T \mathbb{\tilde{E}}   [|b_{\tilde{\rho}_s^\epsilon}(\tilde{X}_s^\epsilon)-b_{\tilde{\rho}_s}(\tilde{X}_s)|]    \mathrm{d}s\to0,   \quad  \text{as}     \quad     \epsilon  \to  0, \\
 &\int_0^T \mathbb{\tilde{E}}  [  |\sigma^{2}   \lambda_{\tilde{\rho}_s^\epsilon}^{-2}(\tilde{X}_s^\epsilon)  -\sigma^{2}   \lambda_{\tilde{\rho}_s }^{-2}(\tilde{X}_s ) | ] \mathrm{d}s\to0,   \quad  \text{as}     \quad     \epsilon  \to  0. 
\end{aligned}
\end{equation}
From \eqref{eq:joint-Skorokhod-convergence},  we  have 
$$
 \sup_{s\le T}|\tilde{Z}_s^\epsilon -\tilde{X}_s|   \rightarrow0
 \quad\text{almost surely, as } \epsilon \to 0.
$$
Since $\nabla h$ is bounded and continuous,  
$$
 \Big|\nabla h(\tilde{Z}_s^\epsilon )-\nabla h(\tilde{X}_s)\Big|^2
 \rightarrow0
 \quad\text{$\mathbb{P}  \otimes \mathrm{d}s$-almost everywhere, as } \epsilon \to 0,
$$
using the dominated convergence theorem, we obtain
\begin{equation}\label{eq:test-gradient-L2-convergence}
 \int_0^T    \mathbb{\tilde{E}} \Big[  |
 \nabla h(\tilde{Z}_s^\epsilon )-\nabla h(\tilde{X}_s)
   |^2  \Big]  \mathrm{ d} s   \rightarrow0
 \quad\text{as } \epsilon \to 0.
\end{equation}
By   \eqref{eq:time-integrated-coeff-convergence},  \eqref{eq:test-gradient-L2-convergence} and  Cauchy--Schwarz's inequality, we  get 
\begin{align*}
&\int_0^T \mathbb{\tilde{E}}   [|   b_{\tilde{\rho}_s^\epsilon}(\tilde{X}_s^\epsilon)\cdot\nabla h(\tilde{Z}_s^\epsilon)
-b_{\tilde{\rho}_s}(\tilde{X}_s)\cdot\nabla h(\tilde{X}_s)|]    \mathrm{d}s\\
\leq  {} &\int_0^T \mathbb{\tilde{E}}   [|(b_{\tilde{\rho}_s^\epsilon}(\tilde{X}_s^\epsilon)-b_{\tilde{\rho}_s}(\tilde{X}_s))
 \cdot\nabla h(\tilde{X}_s)|]    \mathrm{d}s\\
 {}& +\int_0^T \mathbb{\tilde{E}}   [|b_{\tilde{\rho}_s^\epsilon}(\tilde{X}_s^\epsilon)\cdot
 (\nabla h(\tilde{Z}_s^\epsilon)-\nabla h(\tilde{X}_s))|]    \mathrm{d}s\\
 \leq   {}&   \sup_{s\le T} |\nabla h(\tilde{X}_s) | 
 \int_0^T   \mathbb{\tilde{E}}   [|(b_{\tilde{\rho}_s^\epsilon}(\tilde{X}_s^\epsilon)-b_{\tilde{\rho}_s}(\tilde{X}_s))| ]\mathrm{d} s \\
  {}&  +   \Big(\int_0^T   \mathbb{\tilde{E}}\Big[|b_{\tilde{\rho}_s^\epsilon}(\tilde{X}_s^\epsilon)|^2 \Big]\mathrm{d} s\Big)^{1/2}
 \Big(\int_0^T \mathbb{\tilde{E}}   \Big[   |
 \nabla h(\tilde{Z}_s^\epsilon)-\nabla h(\tilde{X}_s)
  |^2  \Big] \mathrm{d} s\Big)^{1/2}\\
 & \rightarrow0,
 \quad\text{as } \epsilon \to 0.
\end{align*}
Similarly, using \eqref{eq:time-integrated-coeff-convergence}  and the dominated convergence theorem, we obtain
\begin{align*}
 \int_0^T \mathbb{\tilde{E}}   \bigg[\bigg|  \frac{\sigma^{2}   \lambda_{\tilde{\rho}_s^\epsilon}^{-2}   (\tilde{X}_s^\epsilon)}{2}   \Delta h(\tilde{Z}_s^\epsilon)  
-\frac{\sigma^{2}   \lambda_{\tilde{\rho}_s}^{-2}   (\tilde{X}_s)}{2}   \Delta h(\tilde{X}_s)\bigg|\bigg]    \mathrm{d}s 
 \rightarrow0,
 \quad\text{as } \epsilon \to 0.
\end{align*}
Since $h$ is bounded and continuous, and the paths converge uniformly almost everywhere,
\begin{equation}\label{eq:test-endpoint-uniform-L1}
  \mathbb{\tilde{E}}   \Big[\sup_{t\le T}|h(\tilde{Z}_t^\epsilon)-h(\tilde{X}_t)| \Big]
 \rightarrow0,
 \qquad
  \mathbb{\tilde{E}} [|h(\tilde{Z}_0^\epsilon)-h(\tilde{X}_0)|  ]  \rightarrow0,
 \quad\text{as } \epsilon \to 0.
\end{equation}
Hence 
\begin{equation}\label{eq:test-martingale-L1-convergence}
  \mathbb{\tilde{E}}   \Big[\sup_{t\le T}|\mathcal{\tilde{M}}_t^{\epsilon,h}-\mathcal{\tilde{M}}_t^h|  \Big]
 \rightarrow 0,
 \quad\text{as } \epsilon \to 0.
\end{equation}

Let $0\le s<t\le T$ and let $ \Upsilon$ be any bounded continuous function on $\mathcal{C}([0,T];\mathbb{R}^d)$ that is measurable with respect to the past up to time $s$. Since $\mathcal{\tilde{M}}^{\epsilon,h}$ is a martingale with respect to the filtration generated by $(\tilde{X}^\epsilon,\tilde{B}^\epsilon)$, we have
\begin{equation}\label{eq:prelimit-martingale-law-identity}
\tilde{\mathbb{E}}\left[(\mathcal{\tilde{M}}_t^{\epsilon,h}-\mathcal{\tilde{M}}_s^{\epsilon,h})\Upsilon(\tilde{X}^\epsilon|_{[0,s]})\right]=0. 
\end{equation}
Using the a.s. convergence $\tilde{X}^\epsilon  \to  \tilde{X}$ in the uniform topology, we have $\Upsilon(\tilde{X}^\epsilon|_{[0,s]})\to  \Upsilon(\tilde{X}|_{[0,s]})$ a.s. Also, by \eqref{eq:test-martingale-L1-convergence}, $\mathbb{\tilde{E}}   \Big[ \mathcal{\tilde{M}}_t^{\epsilon,h} -\mathcal{\tilde{M}}_s^{\epsilon,h}   - (   \mathcal{\tilde{M}}_t^h- \mathcal{\tilde{M}}_s^h)  \Big]  \rightarrow 0$,  as  $ \epsilon \to 0$. Since both sequences are uniformly integrable, we may pass to the limit in \eqref{eq:prelimit-martingale-law-identity}  and obtain
$$
\tilde{\mathbb{E}}\left[(\mathcal{\tilde{M}}_t^h-\mathcal{\tilde{M}}_s^h)\Upsilon(\tilde X|_{[0,s]})\right]=0.
$$
By the monotone class theorem, this holds for all bounded functions measurable with respect to the natural filtration of $\tilde{X}$ up to time $s$; hence $\mathcal{\tilde{M}}^h$ is a martingale with respect to that filtration.

For $h,  \hbar\in \mathcal{C}_c^\infty$, define
$$
\mathbf{\tilde{A}}_t^{\epsilon,h,  \hbar}:=
\int_0^t  \sigma^{2}   \lambda_{\tilde{\rho}_s^\epsilon}^{-2} (\tilde{X}_s^\epsilon)   \nabla h(\tilde{Z}_s^\epsilon)^\top \nabla \hbar(\tilde{Z}_s^\epsilon) \mathrm{d}s.
$$
By the same argument as in the proof of the $L^1$-convergence of test-function martingales, we have
\begin{equation}\label{eq:cross-bracket-L1-convergence}
  \mathbb{\tilde{E}}\Big[  \sup_{t\le T}|\mathbf{\tilde{A}}_t^{\epsilon,h,  \hbar}-    \mathbf{\tilde{A}}_t^{ h,  \hbar}|\Big]
 \rightarrow0,   \quad\text{as } \epsilon \to 0, 
\end{equation}
where
$$  \mathbf{\tilde{A}}_t^{ h,  \hbar}:=
\int_0^t  \sigma^{2}   \lambda_{\tilde{\rho}_s }^{-2} (\tilde{X}_s )   \nabla h(\tilde{X}_s )^\top \nabla \hbar(\tilde{X}_s ) \mathrm{d}s.$$
Ito's product formula gives that
$
\mathcal{\tilde{M}}_t^{\epsilon,h}   \mathcal{\tilde{M}}_t^{\epsilon,\hbar}-   \mathbf{\tilde{A}}_t^{\epsilon,h,  \hbar}
$
is a martingale. 
First, \eqref{eq:test-martingale-L1-convergence} gives that
$(\mathcal{\tilde{M}}_t^{\epsilon,h},\mathcal{\tilde{M}}_t^{\epsilon, \hbar})\to(\mathcal{\tilde{M}}_t^h,\mathcal{\tilde{M}}_t^ \hbar)$ in probability,
 as $ \epsilon \to 0$;
hence, by the continuous mapping theorem,  we  have 
\begin{equation}\label{eq:martingale-product-probability}
\mathcal{\tilde{M}}_t^{\epsilon,h} \mathcal{\tilde{M}}_t^{\epsilon, \hbar}  \rightarrow \mathcal{\tilde{M}}_t^h \mathcal{\tilde{M}}_t^ \hbar
 \quad\text{in probability},
 \quad\text{as } \epsilon \to 0.
\end{equation}
Using \eqref{eq:prelimit-martingale-pmom}, the family $\{\mathcal{\tilde{M}}_t^{\epsilon,h}   \mathcal{\tilde{M}}_t^{\epsilon,\hbar}\}$ is uniformly integrable.
 Consequently, \eqref{eq:martingale-product-probability}
and Vitali's theorem yield
\begin{equation}\label{eq:martingale-product-L1}
  \mathbb{\tilde{E}}  [|\mathcal{\tilde{M}}_t^{\epsilon,h} \mathcal{\tilde{M}}_t^{\epsilon, \hbar}  \rightarrow \mathcal{\tilde{M}}_t^h \mathcal{\tilde{M}}_t^ \hbar |  ]
 \rightarrow 0,
 \quad\text{as } \epsilon \to 0.
\end{equation}
Combining this with \eqref{eq:cross-bracket-L1-convergence}, for each fixed $t$,  we obtain
$$
  \mathbb{\tilde{E}}[ |\mathcal{\tilde{M}}_t^{\epsilon,h}   \mathcal{\tilde{M}}_t^{\epsilon,\hbar}-   \mathbf{\tilde{A}}_t^{\epsilon,h,  \hbar} - (\mathcal{\tilde{M}}_t^{ h}   \mathcal{\tilde{M}}_t^{ \hbar}-   \mathbf{\tilde{A}}_t^{  h,  \hbar} )| ]\rightarrow 0,
 \quad\text{as } \epsilon \to 0.
$$
Applying the same martingale transfer argument  to this product martingale yields that
$
\mathcal{\tilde{M}}_t^{ h}   \mathcal{\tilde{M}}_t^{ \hbar}-   \mathbf{\tilde{A}}_t^{  h,  \hbar} 
$
is a martingale. Therefore, by the uniqueness of the predictable quadratic covariation for continuous martingales, we have
\begin{equation}\label{eq:limit-bracket}
   \langle \mathcal{\tilde{M}}_t^{ h} , \mathcal{\tilde{M}}_t^{ \hbar} \rangle_t =  \mathbf{\tilde{A}}_t^{  h,  \hbar}
.
\end{equation}

Note that  $\tilde{X}$ satisfies the martingale problem with generator
$$
(\mathscr{L} f)(x)=b_{\tilde{\rho} }(x)\cdot \nabla f(x)  + \frac{\sigma^{2}   \lambda_{\tilde{\rho} }^{-2}   (x)}{2}   \Delta f(x).
$$
 Consequently, there exists a Brownian motion $\tilde{B}$   such that
$$
\tilde{X}_t = \tilde{X}_0 + \int_0^t b_{\tilde{\rho}_s}(\tilde{X}_s)\mathrm{ d}s + \int_0^t \sigma_{\tilde{\rho}_s}(\tilde{X}_s) \mathrm{d}\tilde{B}_s,
$$
with $\tilde{\rho}_t=\mathrm{ Law}(\tilde X_t)$. This is precisely the weak   form of the McKean--Vlasov SDE \eqref{eq:limit-SDE}.
Any convergent subsequence of $\{\mathrm{ Law}(X^\epsilon)\}$ has a limit that satisfies the McKean--Vlasov SDE \eqref{eq:limit-SDE}. By the assumed uniqueness of the weak solution, all subsequential limits coincide with the unique distribution $\mathfrak{P}$. Since the family is tight, the full sequence converges weakly:
$
\mathrm{ Law}(X^\epsilon)\Longrightarrow \mathfrak{P}.
$
Let $\rho_t := \mathfrak{P}\circ X_t^{-1}$. Since $\tilde X$ has distribution $\mathfrak{P}$, we have $\tilde{\rho}_t = \rho_t$ for every $t$. Combining this with \eqref{eq:path-W1}, we obtain  
$
\sup_{t\le T} \mathbb{W}_1(\rho_t^\epsilon,\rho_t) \longrightarrow 0,
$
where $\rho_t^\epsilon = \mathrm{ Law}(X_t^\epsilon)$. 
 The proof is finished.

\vspace{0.5cm}

By Lemma \ref{lem:limit-Lipschitz}, the coefficients of the McKean--Vlasov SDE \eqref{eq:limit-SDE} are globally Lipschitz. Consequently, the contraction mapping theorem yields the existence and uniqueness of solutions to the McKean--Vlasov SDE \eqref{eq:limit-SDE}. We now proceed to give the proof of strong convergence.

 \noindent{  \emph{\textbf{Proof   of   Theorem \ref{thm:strong-convergence}:} }}
 Set $\Delta_\epsilon = \epsilon ^{4/5}$.
 Let
$$
 Z_t^\epsilon-X_t^\epsilon=\epsilon ^{1/2}\Xi_t^\epsilon,
 \qquad
 \widetilde{Z}_t^\epsilon=Z_t^\epsilon-\mathfrak{R}_t^\epsilon.
$$
By Proposition  \ref{prop:tightness-complete} and Remark \ref{rem:core-L1-averaging},  we  have 
\begin{equation}\label{eq:UR-L1}
  \mathbb{E}\Big[\sup_{0  \le   t\le T} |\epsilon ^{1/2}\Xi_t^\epsilon|^2 \Big]\leq C_{T} \epsilon^{1/2}, \quad \mathbb{E}\Big[\sup_{0  \le   t\le T} |\mathfrak{R}_t^\epsilon|  \Big]\leq C_{T} \epsilon^{1/5}.
\end{equation}
By   Cauchy--Schwarz's inequality,  \eqref{eq:UR-L1} and    $\mathfrak{R}_0^\epsilon  =0$,  we get
\begin{align}\label{eq:original-Ztilde-X-L1}
           \mathbb{E} \Big[ \sup_{0  \le r\le  T}   | \widetilde{Z}_r^\epsilon   - X_r^\epsilon      | \Big]    =  {}&    \mathbb{E} \Big[ \sup_{0  \le r\le  T}   | \epsilon ^{1/2}\Xi_r^\epsilon   - \mathfrak{R}_r^\epsilon    | \Big] \notag\\ 
           \leq  {}&\Big( \mathbb{E}\Big[\sup_{0  \le r\le  T}  |\epsilon ^{1/2}\Xi_r^\epsilon|^2 \Big]\Big)^{1/2} + \mathbb{E}\Big[\sup_{0  \le r\le  T} |\mathfrak{R}_r^\epsilon|  \Big]   \notag \\  
           \leq  {}&C_{T} \epsilon^{1/5}      
\end{align}
and 
\begin{align}\label{eq:original-Ztilde-X-0-L1}
       \mathbb{E} \Big[  |\widetilde{Z}_0^\epsilon-X_0  |  \Big]    \leq  {}&     \mathbb{E} \Big[  | X _0^\epsilon -X_0  |  \Big]     +    \mathbb{E} \Big[ |\epsilon ^{1/2}\Xi_0^\epsilon|  \Big]    \notag \\ 
     \leq   {}&  \Big(  \mathbb{E} \Big[  | X _0^\epsilon -X_0  |^{2}  \Big]   \Big)^{1/2}  +   \Big( \mathbb{E} \Big[ |\epsilon ^{1/2}\Xi_0^\epsilon|^{2}   \Big] \Big)^{1/2}   \notag \\ 
     \leq   {} &  \Big(  \mathbb{E} \Big[  | X _0^\epsilon -X_0  |^{2}  \Big]   \Big)^{1/2}  + C_{T} \epsilon^{1/4}.  
\end{align}
Let
$
  \mathcal{U} _t^\epsilon:=\widetilde{Z}_t^\epsilon-X_t $  and 
$ \rho_t= \mathrm{Law}(X_t).
$
Take $h\in(0,1]$ to be determined later, and set
$
 s_k=kh \wedge T,$    $0\le k\le K_h:=\lceil T/h\rceil.
$
Note  that 
$$
  \mathrm{d}  \widetilde{Z}_t^\epsilon
 =b_{\rho_t^\epsilon}(X_t^\epsilon) \mathrm{d} t
 + \lambda_{\rho_t  ^\epsilon }^{-1}  (X_t^\epsilon)  \sigma  \mathrm{d} B_t 
$$
and 
$$
  \mathrm{ d} X_t=b_{\rho_t}(X_t) \mathrm{d} t+\lambda_{\rho_{t}}^{-1}  (X_t)  \sigma     \mathrm{d} B_t.
$$
For $s\in[s_k,s_{k+1}]$, pathwise it holds that
$$
 |\mathcal{U}_s^\epsilon|
 \le |\mathcal{U}_{s_k}^\epsilon|
 +\sup_{s_k\le r\le s_{k+1}}
  |\mathcal{U}_r^\epsilon    -   \mathcal{U}_{s_k}^\epsilon|.
$$
Using Lemma  \ref{lem:limit-Lipschitz},  we  have 
\begin{align}\label{eq:direct-L1-coeff-block-sup}
 & \mathbb{E} \Big[ \sup_{s_k\le r\le s_{k+1}}
 \bigl(| b_{\rho_r^\epsilon}(X_r^\epsilon)  - b_{\rho_r}(X_r)|+ |\lambda_{\rho_r^\epsilon}^{-1}(X_r^\epsilon) \sigma  -  \lambda_{\rho_r}^{-1}  (X_r)  \sigma |\bigr)  \Big]  \notag  \\
 \le   {}& C   \mathbb{E} \Big[ \sup_{s_k\le r\le s_{k+1}}
 \bigl(| X_{r}^\epsilon     - X _{r} |+\mathbb{W}_1(\rho_{r}^\epsilon, \rho _{r})\bigr)  \Big]  \notag  \\
  \le  {} & C   \mathbb{E} \Big[ \sup_{s_k\le r\le s_{k+1}}
 \bigl(| X_{r}^\epsilon     - X _{r} | \bigr)  \Big]  \notag  \\
   \le   {}& C   \Big(     \mathbb{E} \Big[ \sup_{s_k\le r\le s_{k+1}} |\mathcal{U}_r^\epsilon  |  \Big]                 
      +   \mathbb{E} \Big[ \sup_{s_k\le r\le s_{k+1}}   | \widetilde{Z}_r^\epsilon   - X_r^\epsilon      | \Big]\Big)       \notag  \\
  \le   {}&  C   \Big(     \mathbb{E} \Big[ \sup_{0  \le r\le s_{k }} |\mathcal{U}_r^\epsilon  |  \Big]   +              \mathbb{E} \Big[ \sup_{s_k\le r\le s_{k+1}}
  |\mathcal{U}_r^\epsilon    -   \mathcal{U}_{s_k}^\epsilon|  \Big]    +   \mathbb{E} \Big[ \sup_{0  \le r\le  T}   | \widetilde{Z}_r^\epsilon   - X_r^\epsilon      | \Big]\Big).  
\end{align}
 From   Burkholder-Davis-Gundy's inequality and  \eqref{eq:direct-L1-coeff-block-sup},  we obtain
\begin{align}
& \mathbb{E} \Big[ \sup_{s_k\le r\le s_{k+1}}
  |\mathcal{U}_r^\epsilon    -   \mathcal{U}_{s_k}^\epsilon|  \Big] \notag \\
 =  {}&  \mathbb{E} \Big[  \sup_  {    s_k    \le   r\le s_{k  +1}  }|\widetilde{Z}_r^\epsilon-X_r   -(\widetilde{Z}_{s_k}^\epsilon-X_{s_k})| \Big]  \notag\\
 \leq  {}&    \mathbb{E}  \bigg[ \sup_{ s_k    \le   r\le s_{k  +1}} \bigg|  \int_{{s_k}}^{r }   b_{\rho_s^\epsilon}(X_s^\epsilon)  - b_{\rho_s}(X_s)    \mathrm{d} s   \bigg|   \bigg]  +   \mathbb{E}  \bigg[ \sup_{ s_k    \le   r\le s_{k  +1}} \bigg|\int_{ s_k }^{r} \lambda_{\rho_s^\epsilon}^{-1}  (X_s^\epsilon) \sigma  -  \lambda_{\rho_s}^{-1}  (X_s)  \sigma\mathrm{d} B_s \bigg|   \bigg]\notag\\
\leq  {}&  h \mathbb{E} \Big[ \sup_{s_k\le r\le s_{k+1}}
  | b_{\rho_r^\epsilon}(X_r^\epsilon)  - b_{\rho_r}(X_r)|    \Big]    +  C \mathbb{E}   \bigg(\int_{s_{k  }}^{s_{k  +1}}  \Big(\sup_{    s_k    \le   r\le s_{k  +1}  } |\lambda_{\rho_r^\epsilon}^{-1}(X_r^\epsilon) \sigma  -  \lambda_{\rho_r}^{-1}  (X_r)  \sigma |\Big)^{2}\mathrm{d} s\bigg)^{1/2}\notag\\
\leq  {}&    h \mathbb{E} \Big[ \sup_{s_k\le r\le s_{k+1}}
  | b_{\rho_r^\epsilon}(X_r^\epsilon)  - b_{\rho_r}(X_r)|    \Big] +  C\sqrt{h} \mathbb{E} \Big[ \sup_{s_k\le r\le s_{k+1}}
  |\lambda_{\rho_r^\epsilon}^{-1}(X_r^\epsilon) \sigma  -  \lambda_{\rho_r}^{-1}  (X_r)  \sigma |  \Big]   \notag\\
\leq  {}&      \theta_h   \Big(     \mathbb{E} \Big[ \sup_{0  \le r\le s_{k }} |\mathcal{U}_r^\epsilon  |  \Big]   +              \mathbb{E} \Big[ \sup_{s_k\le r\le s_{k+1}}
  |\mathcal{U}_r^\epsilon    -   \mathcal{U}_{s_k}^\epsilon|  \Big]    +   \mathbb{E} \Big[ \sup_{0  \le r\le  T}   | \widetilde{Z}_r^\epsilon   - X_r^\epsilon      | \Big]\Big),\notag
\end{align}
 where  $\theta_h:=C(h+\sqrt{h}).
 $
 Choose a fixed $h_0=h_0(C)>0$ such that $\theta_{h_0}\le 1/2$, and use this same $h=h_0$ for all blocks.
 Then   
 \begin{align*}
\mathbb{E} \Big[ \sup_{s_k\le r\le s_{k+1}}
  |\mathcal{U}_r^\epsilon    -   \mathcal{U}_{s_k}^\epsilon|  \Big]  
\leq      2  \theta_{h_0}   \Big(     \mathbb{E} \Big[ \sup_{0  \le r\le s_{k }} |\mathcal{U}_r^\epsilon  |  \Big]    +   \mathbb{E} \Big[ \sup_{0  \le r\le  T}   | \widetilde{Z}_r^\epsilon   - X_r^\epsilon      | \Big]\Big).
\end{align*}
 Due  to
$
 \sup  \limits_{0\le r\le s_{k+1}}|\mathcal{U}_r^\epsilon|
 \le\sup   \limits_{0\le r\le s_k}|\mathcal{U}_r^\epsilon|
 +\sup   \limits_{s_k\le r\le s_{k+1}}
   |\mathcal{U}_r^\epsilon-\mathcal{U}_{s_k}^\epsilon|,
$
 we obtain
  \begin{align*}
\mathbb{E} \Big[ \sup_{0 \le r\le s_{k+1}}
 |\mathcal{U}_r^\epsilon  |   \Big]  
\leq    (1+  2  \theta_{h_0}  )        \mathbb{E} \Big[ \sup_{0  \le r\le s_{k }} |\mathcal{U}_r^\epsilon  |  \Big]    +  
 2  \theta_{h_0}    \mathbb{E} \Big[ \sup_{0  \le r\le  T}   | \widetilde{Z}_r^\epsilon   - X_r^\epsilon      | \Big] .
\end{align*}
 Upon iteration up to $K_{h_0}\le    1+ \frac{T}{ h_0}$, we get
   \begin{align*}
\mathbb{E} \Big[ \sup_{0 \le r\le T}
 |\widetilde{Z}_r^\epsilon-X_r  |   \Big]    =  &   \mathbb{E} \Big[ \sup_{0 \le r\le K_{h_0}}
 |\mathcal{U}_r^\epsilon  |   \Big]   \\
\leq   {}&  (1+  2  \theta_{h_0}  )        \mathbb{E} \Big[ \sup_{0  \le r\le {K_{h_0}  -1}} |\mathcal{U}_r^\epsilon  |  \Big]    +  
 2  \theta_{h_0}    \mathbb{E} \Big[ \sup_{0  \le r\le  T}   | \widetilde{Z}_r^\epsilon   - X_r^\epsilon      | \Big] \\
 \leq   {}&  (1+  2  \theta_{h_0}  )^{    K_{h_0}   }        \mathbb{E} \Big[  |\mathcal{U}_0^\epsilon  |  \Big]    +  
 2  \theta_{h_0}    \mathbb{E} \Big[ \sup_{0  \le r\le  T}   | \widetilde{Z}_r^\epsilon   - X_r^\epsilon      | \Big]
 \sum _{j =  0}^{  K_{h_0}  -1     }   (1+  2  \theta_{h_0}  )^{j }     \\
   \leq  {} &  C_{T}   \Big(          \mathbb{E} \Big[  |\widetilde{Z}_0^\epsilon-X_0  |  \Big]    +      \mathbb{E} \Big[ \sup_{0  \le r\le  T}   | \widetilde{Z}_r^\epsilon   - X_r^\epsilon      | \Big]                \Big)   .
\end{align*}
 By virtue of   \eqref{eq:original-Ztilde-X-L1}  and   \eqref{eq:original-Ztilde-X-0-L1},   we derive 
    \begin{align*}
\mathbb{E} \Big[ \sup_{0 \le r\le T}
 |X_r^\epsilon-X_r  |   \Big]     \leq {} &  \mathbb{E} \Big[ \sup_{0 \le r\le T}
 |X_r^\epsilon-\widetilde{Z}_r^\epsilon |   \Big]   +   \mathbb{E} \Big[ \sup_{0 \le r\le T}
 |\widetilde{Z}_r^\epsilon-X_r  |   \Big]   \\
   \leq  {} &  C_{T}   \Big(          \mathbb{E} \Big[  |\widetilde{Z}_0^\epsilon-X_0  |  \Big]    +      \mathbb{E} \Big[ \sup_{0  \le r\le  T}   | \widetilde{Z}_r^\epsilon   - X_r^\epsilon      | \Big]                \Big)  \\
 \leq    {} &  C_{T} \Big( \Big(  \mathbb{E} \Big[  | X _0^\epsilon -X_0  |^{2}  \Big]   \Big)^{1/2}  + \epsilon^{1/5}\Big).
\end{align*}
 The proof is finished.

 \begin{remark}
Under  $ \sup  \limits_{0<\epsilon  \le  1}  \mathbb{E}\bigl[|X_0^\epsilon|^4+|Y_0^\epsilon |^4\bigr]<\infty$  and  Assumption \textup{\textbf{(A1\textquotesingle)}-\textbf{(A2\textquotesingle)}}, the existing Khasminskii blocking proof yields an upper bound of order $\epsilon^{1/5}$; however, the following steps lead to a loss of order.
\begin{enumerate}
\item If one can prove the uniform fourth-order increment estimate
$
\mathbb{E}|X_t^\epsilon - X_s^\epsilon|^4 \leq C(|t-s| + \epsilon)^2,
$
then it is promising to improve the mean-square error of the frozen velocity   in   Lemma \ref{lem:freeze-coupling}   to
$
\sup  \limits_{0\le k<K_\epsilon}\sup  \limits_{s\in I_k}
\mathbb{E}\left[|Y_s^\epsilon-\widehat{Y}_s^{\epsilon,k}|^2\right]
\le C(\Delta_\epsilon + \epsilon),
$
thereby improving the freezing error of $\mathcal{Q}$ from $\Delta_\epsilon^{1/4}$ to $\Delta_\epsilon^{1/2}$.
\item To avoid direct differentiation with respect to the measure parameter, one combines Khasminskii blocking with the Poisson equation. However, in the proof of the averaging lemma, a term $\epsilon/\Delta_\epsilon$ arises. 
    If one can avoid block freezing and instead use the global Poisson equation directly, then $\epsilon/\Delta_\epsilon$ can be eliminated.
\item We have only obtained
$
\mathbb{E}\Big[\sup\limits  _{t\le T}|\epsilon^{1/2}\Xi(q_t^\epsilon;X_t^\epsilon,Y_t^\epsilon)|^2\Big]
\le C_T \epsilon^{1/2}.
$
\end{enumerate}
\end{remark}


\section*{Acknowledgements}\label{Acknowledgements}
The authors are grateful to the anonymous referees for their fruitful comments to improve this paper.


\section*{Declaration of competing interest}

The authors declare no conflict of interest.

%


\appendix

\section{ Some  estimates of  several operators}

In section,  we  provide estimates for several operators.   
 First, recall the definitions of $B_\rho$ and $M_\rho$. For each coordinate
direction $a=1,\dots,d$, define the corresponding components by
\begin{align*}
 (B_{\rho,a}h)(x)
 &:=\frac{1}{\lambda_\rho(x)}
   \int_{\mathbb{R}^d}\partial_a\phi(x-x')\,h(x')\,\rho(\mathrm{d}x'),\\
 (M_{\rho,a}h)(x)
 &:=\frac{\partial_a\lambda_\rho(x)}{\lambda_\rho(x)}\,h(x).
\end{align*}

\begin{lemma}\label{lem:operator-insertion}
Set $\mu_\rho(\mathrm{d}  x) := \lambda_\rho(x) \rho(\mathrm{d}  x)$. Under \textup{\textbf{(A1)}-\textbf{(A2)}},   for any $ 1\leq   p \leq   \infty$  and each coordinate $a=1,\dots,d$ 
\begin{equation}\label{eq:PBM-bounds}
 \|B_{\rho,a}\|_{L^p(\mu_\rho) \to L^p(\mu_\rho)}    \le L_\phi,
 \qquad \|M_{\rho,a}\|_{L^p(\mu_\rho) \to L^p(\mu_\rho)}       \le L_\phi,
 \qquad \|R_\rho\|_{L^p(\mu_\rho) \to L^p(\mu_\rho)}   \le1,
\end{equation}
where $p<\infty$ refers to the $L^p(\mu_\rho)$ norm, and $p=\infty$ refers to the $L^\infty(\mu_\rho)$ norm.
Moreover, for any suitable function $h$, setting 
\[
 (\mathscr{D}_{n,i}^{a}h)(x_0)
 :=\int _{   \mathbb{R} ^{nd}}  \partial_{x_i^a}w_{n,\rho}(x_0,\ldots,x_n)
 h(x_n)\rho^{\otimes n}(\mathrm{d}  x_1,  \dots , \mathrm{d}  x_n),
\] 
 we have
\begin{equation}\label{eq:exact-insertion}
\begin{aligned}
\mathscr{D}_{n,i}^{a}
={}&\boldsymbol1_{\{i<n\}}P_\rho^iB_{\rho,a}
       P_\rho^{n-i-1}R_\rho\\
&-\boldsymbol1_{\{i>0\}}P_\rho^{i-1}B_{\rho,a}
       P_\rho^{n-i}R_\rho
-P_\rho^iM_{\rho,a}P_\rho^{n-i}R_\rho .
\end{aligned}
\end{equation}
In particular, for \(n\ge 1\),
\begin{equation}\label{eq:D-bound}
\|\mathscr{D}_{n,i}^{a}h\|_{L^p(\mu_\rho)}
\le C r_*^{n-1}\|h\|_{L^p(\mu_\rho)},
\end{equation}
while for $n=0$, the right-hand side is bounded by $C\|h\|_{L^p(\mu_\rho)}$.
\end{lemma}
\proof 

The kernel of $B_{\rho,a}$ with respect to $\mu_\rho$ is
$$
 b_{\rho,a}(x,x')
 =\frac{\partial_a\phi(x-x')}{\lambda_\rho(x)\lambda_\rho(x')}.
$$
Under Assumption \textbf{(A1)}, it holds that
$$
 \int _{ \mathbb{R}^d}  |b_{\rho,a}(x,x')| \mu_\rho(\mathrm{d} x')
 =\frac{1}{\lambda_\rho(x)}     \int_{ \mathbb{R}^d} |\partial_a\phi(x-x') |  \rho(\mathrm{d} x') \le \frac{L_\phi}{\lambda_\rho(x)} \le L_\phi .
$$
Then the absolute row sums and column sums are bounded by $L_\phi$. Similarly, for $b_{\rho,a}$, by the same $L^1\)-\(L^\infty$ argument as in Lemma \ref{lem:contraction}, we obtain $\|B_{\rho,a}\|_{L^p(\mu_\rho) \to L^p(\mu_\rho)} \le L_\phi$.
Furthermore, since
$$
|\partial_a\lambda_\rho|  \le \int_{ \mathbb{R}^d} |\partial_a\phi(x-x') |\,\rho(\mathrm{d}x') \le L_\phi
$$
and $\lambda_\rho\ge1$, we have
$$
\|M_{\rho,a}\|_{L^p(\mu_\rho) \to L^p(\mu_\rho)} \le \sup _ { x \in   \mathbb{R}^d}  \left|\frac{\partial_a\lambda_\rho(x)}{\lambda_\rho(x)} \right| \le L_\phi;
$$
the bound for $R_\rho$ is obvious.

Applying the iteration of $P_\rho$ gives
\begin{equation}\label{eq:chain-operator}
 (P_\rho^n   R_\rho h)(x_0)
 =\int_{ \mathbb{R}^{n d}}     w_{n,\rho}(x_0,\ldots,x_n)h(x_n)
 \rho^{\otimes n}(\mathrm{d} x_1, \dots,  \mathrm{d} x_n).
\end{equation}
Then  we have   for $0\le i  \le n$ 
\begin{equation}\label{eq:w-derivative-exact}
\begin{aligned}
 \partial_{x_i^a}w_{n,\rho}
={}&\boldsymbol1_{\{i<n\}}
 \frac{\partial_a\phi(x_i-x_{i+1})
       \prod  \limits_{\substack{0\le j\le n-1\\j\ne i}}\phi(x_j-x_{j+1})}
      {\prod_{j=0}^{n}\lambda_\rho(x_j)}\\
&-\boldsymbol1_{\{i>0\}}
 \frac{\partial_a\phi(x_{i-1}-x_i)
       \prod  \limits  _{\substack{0\le j\le n-1\\j\ne i-1}}\phi(x_j-x_{j+1})}
      {\prod_{j=0}^{n}\lambda_\rho(x_j)}
 -w_{n,\rho}\frac{\partial_a\lambda_\rho(x_i)}
                         {\lambda_\rho(x_i)},
\end{aligned}
\end{equation}
where the negative sign of the second term comes from  
$
\partial_{x_i^a}\phi(x_{i-1}-x_i)=-\partial_a  \phi(x_{i-1}-x_i), 
$
 and    $\partial_a $ denotes differentiation with respect to the 
$a$-th coordinate of the argument of  $\phi.$ 
Multiplying \eqref{eq:w-derivative-exact} by $h(x_n)$ and integrating with respect to $x_1,\ldots,x_n$, we obtain
\begin{align*} 
\mathscr{D}_{n,i}^{a}
={}&\boldsymbol1_{\{i<n\}}P_\rho^iB_{\rho,a}
       P_\rho^{n-i-1}R_\rho\\
&-\boldsymbol1_{\{i>0\}}P_\rho^{i-1}B_{\rho,a}
       P_\rho^{n-i}R_\rho
-P_\rho^iM_{\rho,a}P_\rho^{n-i}R_\rho ,
\end{align*}
where we used the fact that
\begin{align*}
&\int_{     \mathbb{R}^{nd}} \boldsymbol1_{\{i<n\}}
 \frac{\partial_a\phi(x_i-x_{i+1})
       \prod_{j\ne i}\phi(x_j-x_{j+1})}{\prod_{j=0}^{n}\lambda_\rho(x_j)}
 h(x_n)\rho^{\otimes n}(\mathrm{d}   x_1,   \dots,   \mathrm{d} x_n)\\
=&\boldsymbol1_{\{i<n\}}
  (P_\rho^iB_{\rho,a}P_\rho^{n-i-1}R_\rho h)(x_0),
\end{align*}
\begin{align*}
&-\int_{     \mathbb{R}^{nd}} \boldsymbol1_{\{i>0\}}
 \frac{\partial_a\phi(x_{i-1}-x_i)
       \prod_{j\ne i-1}\phi(x_j-x_{j+1})}{\prod_{j=0}^{n}\lambda_\rho(x_j)}
 h(x_n)\rho^{\otimes n}(\mathrm{d}   x_1,   \dots,   \mathrm{d} x_n)\\
=&-\boldsymbol1_{\{i>0\}}
  (P_\rho^{i-1}B_{\rho,a}P_\rho^{n-i}R_\rho h)(x_0) 
\end{align*}
and
$$
 -\int_{     \mathbb{R}^{nd}} w_{n,\rho}
 \frac{\partial_a\lambda_\rho(x_i)}{\lambda_\rho(x_i)}
 h(x_n)\rho^{\otimes n}(\mathrm{d}   x_1,   \cdots,   \mathrm{d} x_n)
 =-(P_\rho^iM_{\rho,a}P_\rho^{n-i}R_\rho h)(x_0).
$$
For   $n\ge1$,   by virtue of Lemma \ref{lem:contraction} and the preceding estimate, we derive   
\begin{align*}
 \|\mathscr{D}_{n,i}^{a}h\|_{L^p(\mu_\rho)}
 &\le L_\phi\bigl(
 \boldsymbol1_{\{i<n\}}r_*^{n-1}
 +\boldsymbol1_{\{i>0\}}r_*^{n-1}+r_*^n\bigr)\|h\|_{L^p(\mu_\rho)}\\
 &\le C r_*^{n-1}\|h\|_{L^p(\mu_\rho)}.
\end{align*}
For $n=0$, we have $\mathscr{D}_{0,0}^a h = -M_{\rho,a} R_\rho h$, which is therefore bounded by $C\|h\|_{L^p(\mu_\rho)}$.

\endproof

Set
$$
 S(x):=\int_{\mathbb{R}^{d}} yy^{ \top }q^x(\mathrm{d} y),\qquad
 m_4(q):=\int_{\mathbb{R}^{2d}}    |y|^4q(\mathrm{d}  x,\mathrm{d} y).
$$
Applying the   Jensen inequality together with $1\le\lambda_\rho\le1+\phi_M$ yields
\begin{equation}\label{eq:uS-norms}
 \|u\|_{L^4(\mu_\rho)}\le C m_4(q)^{1/4},
 \qquad
 \|S\|_{L^2(\mu_\rho)}\le C m_4(q)^{1/2}.
\end{equation}

\begin{lemma} \label{lem:kernel-sum}
Assume that   \textup{\textbf{(A1)}-\textbf{(A2)}} hold and that  $m_4(q) < \infty$.
Then 
\begin{equation}\label{eq:Q-tail-L2}
 \| \mathcal{Q} (q;\cdot)   -\mathcal{Q}^{(N)}(q;\cdot)  \|_{L^2(q)}
 \le C\bigl(1+m_4(q)^{1/2}\bigr)
 \sum_{n>N}(n+1)^2r_{*}^{(n-2)_{+}} 
\end{equation}
and 
\begin{equation}\label{eq:Q-L2-final}
 \|\mathcal{ Q}(q;\cdot)\|_{L^2(q)}
 \le C\bigl(1+m_4(q)^{1/2}\bigr).
\end{equation}
Moreover 
\begin{equation}\label{eq:Xi-tail-L4}
 \|\Xi(q;\cdot)-\Xi^{(N)}(q;\cdot)\|_{L^4(q)}
 \le C r_*^{N+1}m_4(q)^{1/4},\qquad N\ge0 
\end{equation}
 and 
\begin{equation}\label{eq:Xi-L4-final}
 \|\Xi(q;\cdot)\|_{L^4(q)}
 \le C\bigl(1+m_4(q)^{1/4}\bigr).
\end{equation}
\end{lemma}

\proof

In the following, all estimates are performed coordinatewise, and summing over the finite number of coordinates only changes the constant.   Let $\mathsf{M}_f h := f h$, and set
$
\mathsf{m}_2 (x):=\int_{   \mathbb{R}  ^{d}} |y|^2 q^x(\mathrm{d}y).
$
Since $1\le\lambda_\rho\le1+\phi_M$, the $L^p(\rho)$ and $L^p(\mu_\rho)$ norms are equivalent with a uniform constant.

For fixed output coordinate $\ell=1,\dots,d$ and inner product coordinate $a=1,\dots,d$, we define
\begin{equation}\label{eq:Tni-definition}
 T_{n,i}^{\ell a}(x_0,y_0)
 := \int_{   \mathbb{R}^{2nd}  } y_n^\ell y_i^a\,
 \partial_{x_i^a}w_{n,\rho}(x_0,\ldots,x_n)
 q^{\otimes n}(\mathrm{d} z_1,\dots  ,\mathrm{d} z_n).
\end{equation}
In fact, the $\ell$-th component   of  $\int  _{  \mathbb{R}  ^{2nd}}    y_n  y_i \cdot g_{n,i,\rho}(x_0,\ldots,x_n)      q^{\otimes n}( \mathrm{d} z_1,  \dots ,  \mathrm{d} z_n)$ is given by  $\sum_{a=1}^d  T_{n,i}^{\ell a}(x_0,y_0)  .$
Since, given $x_0,\ldots,x_n$, the velocities of the independent copies are conditionally independent, we have
$$
  \mathbb{E}_{1:n}[y_n^\ell y_i^a\mid x_0,\ldots,x_n]
 =u^\ell(x_n)u^a(x_i),\quad 1\le i<n, 
  $$ 
 $$   
  \mathbb{E}_{1:n} [y_n^\ell y_0^a\mid x_0,\ldots,x_n]
 =u^\ell(x_n)  y_0^a
$$
and
$$
  \mathbb{E}_{1:n}[y_n^\ell y_n^a\mid x_0,\ldots,x_n]
 =S_{\ell a}(x_n).
$$
\noindent  {For $1\le i<n$.}   From equation  \eqref{eq:exact-insertion}  
\begin{equation}\label{eq:two-marked-internal}
\begin{aligned}
T_{n,i}^{\ell a}(x_0)
={}&\bigl(P_\rho^i\mathsf{M}_{u^a}B_{\rho,a}
P_\rho^{n-i-1}R_\rho u^\ell\bigr)(x_0)\\
&-\bigl(P_\rho^{i-1}B_{\rho,a}\mathsf M_{u^a}
P_\rho^{n-i}R_\rho u^\ell\bigr)(x_0)\\
&-\bigl(P_\rho^i\mathsf{M}_{u^a}  M_{\rho,a}
P_\rho^{n-i}R_\rho u^\ell\bigr)(x_0).
\end{aligned}
\end{equation}
By Lemmas \ref{lem:contraction} and \ref{lem:operator-insertion},  we have 
\begin{align*}
&\|P_\rho^i\mathsf{M}_{u^a}B_{\rho,a}
P_\rho^{n-i-1}R_\rho u^\ell\|_{L^2(\rho)}  \notag\\
\le  {}& r_*^i\|u^a\|_{L^4(\mu_\rho)}
L_\phi r_*^{n-i-1}\|u^\ell\|_{L^4(\mu_\rho)}   \notag\\
\le  {}& C r_*^{n-1}m_4(q)^{1/2}. 
\end{align*}
Similarly, we   get
\begin{align*}
\|P_\rho^{i-1}B_{\rho,a}\mathsf{M}_{u^a}
P_\rho^{n-i}R_\rho u^\ell\|_{L^2(\rho)}
\le C r_*^{n-1}m_4(q)^{1/2} 
\end{align*}
and 
\begin{align*}
\|P_\rho^i\mathsf M_{u^a}M_{\rho,a}
P_\rho^{n-i}R_\rho u^\ell\|_{L^2(\rho)}
\le C r_*^{n}m_4(q)^{1/2}.
\end{align*}
Hence    
\begin{equation} \label{eq:terminal-i-le-n-mark-bound}
\|   T_{n,i}^{\ell a} \|_{L^2(\rho)}
  \le C r_*^{n-1}m_4(q)^{1/2} .
\end{equation}
\noindent  {For $i=n\ge1  .$}  Note  that 
\begin{equation*}
 T_{n,n}^{\ell a}
 =\mathscr{D}_{n,n}^aS_{\ell a}.
\end{equation*}
By virtue of  \eqref{eq:D-bound} and  \eqref{eq:uS-norms},   we get 
\begin{equation}\label{eq:terminal-mark-bound}
 \|T_{n,n}^{\ell a}\|_{L^2(\rho)}
 \le Cr_*^{n-1}\|S_{\ell a}\|_{L^2(\mu_\rho)}
 \le Cr_*^{n-1}m_4(q)^{1/2}.
\end{equation}
\noindent  {For $  i=0<n$.}   By the definitions  of   $\mathscr{D}_{n,i}^{a} $  and $ T_{n,i}^{\ell a}$, we obtain 
\begin{equation*}
 T_{n,0}^{\ell a}(x_0,y_0)
 =y_0^a   (\mathscr{D}_{n,0}^au^\ell)(x_0).
\end{equation*}
Using H\"{o}lder's inequality, we have
\begin{align*}
 \|y_0^a     (\mathscr{D}_{n,0}^au^\ell)(x_0)  \|_{L^2(q)}^2
 &\le\int _{  \mathbb{R}^{d}}\mathsf{m}_2(x)| (\mathscr{D}_{n,0}^a    u^\ell)(x)|^2\rho(\mathrm{d} x)\notag\\
  &\le\|\mathsf{m}_2\|_{L^2(\rho)}   \|(\mathscr{D}_{n,0}^a    u^\ell)^2 \|_{L^2(\rho)} \notag\\
 &=  \|\mathsf{m}_2\|_{L^2(\rho)}
       \|  \mathscr{D}_{n,0}^a    u^\ell   \|_{L^4(\rho)}^2.
\end{align*}
Applying   Jensen's inequality   gives
$$
 \|\mathsf{m}_2\|_{L^2(\rho)}^2
 =\int_{  \mathbb{R}^{d}}   \left (\int_{   \mathbb{R}  ^{d}} |y|^2 q^x(\mathrm{d}y)  \right)^2\rho(\mathrm{d} x)
 \le m_4(q).
$$
 Equation \eqref{eq:D-bound}    yields
$$
 \|\mathscr{D}_{n,0}^a   u^\ell\|_{L^4(\rho)}
 \le Cr_*^{n-1}\|u^\ell\|_{L^4(\mu_\rho)}
 \le Cr_*^{n-1}m_4(q)^{1/4}.
$$
Then  \begin{equation}\label{eq:i-le-n-derivative-bound}
 \|T_{n,0}^{\ell a}\|_{L^2(q)}
 \le Cr_*^{n-1}m_4(q)^{1/2}.
\end{equation}
\noindent  {For $n=i=0$.}  By the definition   of   $ T_{n,i}^{\ell a}$,    \eqref{eq:weights}  and \textbf{(A1)},   we  obtain 
$$
 T_{0,0}^{\ell a}
 =y_0^\ell y_0^a\,\partial_aR_\rho(x_0),\qquad
 |\partial_aR_\rho(x_0)|\le L_\phi,
$$
and subsequently, this yields
\begin{equation}\label{eq:nzero-derivative-bound}
 \|T_{0,0}^{\ell a}\|_{L^2(q)}
 \le C m_4(q)^{1/2}.
\end{equation}

  Let   $b_q:=\beta_q/\lambda_\rho$  and  by   Cauchy--Schwarz's inequality,   we  have
\begin{equation}\label{eq:beta-bound-detailed}
 \|\beta_q\|_\infty
 \le L_\phi\int_{  \mathbb{R}^{2d}}     |y|   q(\mathrm{d} x,\mathrm{d} y)
 \le C m_4(q)^{1/4},\qquad
 \|u\|_{L^2(\mu_\rho)}\le C m_4(q)^{1/4}.
\end{equation}
Define  
$$
 U_{n,i}^{\ell}(x_0,y_0)
 :=\int_{    \mathbb{R}^{2nd}     } y_n^\ell w_{n,\rho}(x_0,\ldots,x_n)  b_q(x_i)     
 q^{\otimes n}(\mathrm{d} z_1,\dots  ,\mathrm{d} z_n).
$$
For $n\ge1$,   we get      
\begin{equation*}
 U_{n,0}^{\ell}
 =\mathsf{M}_{b_q}   P_\rho^nR_\rho u^\ell,\qquad
 U_{n,i}^{\ell}
 =P_\rho^i\mathsf M_{b_q}P_\rho^{n-i}R_\rho u^\ell,
 \quad1\le i\le n.
\end{equation*}
By Lemmas \ref{lem:contraction} and \ref{lem:operator-insertion},  we have
\begin{equation}\label{eq:beta-chain-bound}
 \|U_{n,i}^{\ell}\|_{L^2(\rho)}
 \le r_*^n\|b_q\|_\infty\|u^\ell\|_{L^2(\mu_\rho)}
 \le Cr_*^n m_4(q)^{1/2}.
\end{equation}
For $n=0$,
$$
 U_{0,0}^{\ell}(x_0,y_0)
 =y_0^\ell R_\rho(x_0)b_q(x_0),
$$
hence, using \eqref{eq:beta-bound-detailed} and $\|Y\|_{L^2(q)}\le m_4(q)^{1/4}$, we obtain
\begin{equation}\label{eq:beta-nzero-bound}
 \|U_{0,0}^{\ell}\|_{L^2(q)}
 \le C m_4(q)^{1/2}.
\end{equation}

From 
\eqref{eq:terminal-i-le-n-mark-bound}-\eqref{eq:nzero-derivative-bound}
and  \eqref{eq:beta-chain-bound}-\eqref{eq:beta-nzero-bound},  we have 
\begin{equation}\label{eq:Q-level-strong}
 \| \mathcal{Q}_n(q;\cdot)\|_{L^2(q)}
 \le C(n+1)r_*^{(n-1)_+}m_4(q)^{1/2}  \le C(n+1)^2r_*^{(n-2)_+}
 \bigl(1+m_4(q)^{1/2}\bigr)  ,
\end{equation}
where   $\mathcal{Q}_{n}$ is the $n$-th layer of $\mathcal{Q}$  and 
the sums over the finite-dimensional coordinates $a,\ell$ have been absorbed into the constant $C$; and in order to ensure that it holds for all $n\ge 0$, we enlarge it to the larger upper bound
$
C(n+1)^2 r_*^{(n-2)_+}
 \bigl(1+m_4(q)^{1/2}\bigr).
$
We apply the  Minkowski's inequality   to    obtain
$$
 \left\|\sum_{n>N}   \mathcal{Q}_n\right\|_{L^2(q)}
 \le\sum_{n>N}\|  \mathcal{Q}_n\|_{L^2(q)},
$$
  and combining with \eqref{eq:Q-level-strong} 
\begin{equation*} 
 \| \mathcal{Q} (q;\cdot)   -\mathcal{Q}^{(N)}(q;\cdot)  \|_{L^2(q)}
 \le C\bigl(1+m_4(q)^{1/2}\bigr)
 \sum_{n>N}(n+1)^2r_{*}^{(n-2)_{+}} .
\end{equation*}
Due  to  $r_*<1$,  we get 
\[
 \sum_{n=0}^\infty(n+1)^2r_*^{(n-2)_+}<\infty,
\]
which  implies  
\begin{equation*} 
 \|\mathcal{Q}(q;\cdot)\|_{L^2(q)}
 \le C\bigl(1+m_4(q)^{1/2}\bigr).
\end{equation*}
On the other hand,  by Lemmas \ref{lem:contraction} and \ref{lem:operator-insertion}, we find 
\begin{align*}
 \|\Xi-\Xi^{(N)}\|_{L^4(q)}
=&  \left\|\sum_{n=N+1}^\infty
 P_\rho^nR_\rho u\right\|_{L^4(\rho)}\\
\le& \sum_{n=N+1}^\infty r_*^n
 \|u\|_{L^4(\mu_\rho)}   \\
 \le &\frac{C r_*^{N+1}}{1-r_*}m_4(q)^{1/4}\\
  \le &  C r_*^{N+1}m_4(q)^{1/4}.
\end{align*}
Together with    the  definition   of  $\Xi  $   and     $\|R_\rho Y\|_{L^4(q)}\le m_4(q)^{1/4}$, we obtain
\begin{equation*} 
 \|\Xi(q;\cdot)\|_{L^4(q)}
 \le C\bigl(1+m_4(q)^{1/4}\bigr).
\end{equation*}

\endproof

\begin{corollary} 
Assume that   \textup{\textbf{(A1)}-\textbf{(A2)}} hold and that $m_{2p}(q):=\int_{\mathbb{R}^{2d}}    |y|^{2p}  q(\mathrm{d}  x,\mathrm{d} y) < \infty$  for  $p\geq 2 $.
Then 
\begin{equation*} 
 \| \mathcal{Q} (q;\cdot)   -\mathcal{Q}^{(N)}(q;\cdot)  \|_{L^{p}(q)}
 \le C_{p}\bigl(1+m_{2p}(q)^{1/p}\bigr)
 \sum_{n>N}(n+1)^2r_{*}^{(n-2)_{+}} 
\end{equation*}
and 
\begin{equation*}  
 \|\mathcal{ Q}(q;\cdot)\|_{L^{p}(q)}
 \le C_{p}\bigl(1+m_{2p}(q)^{1/p}\bigr).
\end{equation*}
Moreover 
\begin{equation*} 
 \|\Xi(q;\cdot)-\Xi^{(N)}(q;\cdot)\|_{L^{2p}(q)}
 \le C_{p} r_*^{N+1}m_{2p}(q)^{1/2p},\qquad N\ge0   
\end{equation*}
 and 
\begin{equation*} 
 \|\Xi(q;\cdot)\|_{L^{2p}(q)}
 \le C_{p}  \bigl(1+m_{2p}(q)^{1/2p}\bigr).
\end{equation*}
\end{corollary}

\proof

It suffices to replace every $L^4$ space by $L^{2p}$ and every $L^2$ space by $L^p$ in the proof of Lemma \ref{lem:kernel-sum}; the conclusion then follows in the same manner.
\endproof

\medskip

\section{Stability of $\mathcal Q$ and $\overline{\mathcal{Q}}$}

\begin{lemma}\label{lem:freeze-coupling}
Under \textup{\textbf{(A1)}-\textbf{(A2)}}, if $\Delta_\epsilon \to 0$ as $\epsilon \to 0$, then
\begin{equation}\label{eq:velocity-freeze-error}
\sup_{0\le k<K_\epsilon}\sup_{s\in I_k}
 \mathbb{E}\Big[|Y_s^\epsilon-\widehat{Y}_s^{\epsilon,k}|^2\Big]
\le C_T\big((\Delta_\epsilon )^{1/2}+\epsilon\big)  \rightarrow  0   \quad   \text{as} \quad   \epsilon \rightarrow 0 .
\end{equation}
\end{lemma}
\proof

 Using the  It\^{o}'s formula  and symmetry    we have,    
\begin{align*}
&\mathrm{d}   |Y_s^\epsilon-\widehat{Y}_s^{\epsilon,k}|^2 \\
 =  {}& 2\epsilon^{-1}  (Y_s^\epsilon-\widehat{Y}_s^{\epsilon,k} )   \cdot  \bigg[-(Y_s^\epsilon   -\widehat{Y}_s^{\epsilon,k} )
+    \bigg(\int _{   \mathbb{R}   ^{2d}} \phi(X_s^\epsilon   -x')  \big(y'-Y_s^\epsilon   \big)
q_s^\epsilon    (\mathrm{d}x',\mathrm{d}y') \\
 {}& -\int _{   \mathbb{R}   ^{2d}}  \phi(\widehat X_s^{\epsilon,k}-x')
(y'-\widehat Y_s^{\epsilon,k})\widehat q_s^{\epsilon,k}(\mathrm{d}  x',\mathrm{d}  y')  \bigg)\bigg]  \mathrm{d}t  -2 \epsilon^{-1/2} (Y_s^\epsilon-\widehat{Y}_s^{\epsilon,k} )   \cdot   G_{\rho_s^\epsilon}(X_s^\epsilon  ) \mathrm{d}t 
\end{align*}
and 
\begin{align*}
& \mathbb{E}  \bigg[(Y_s^\epsilon-\widehat{Y}_s^{\epsilon,k} )   \cdot  \bigg(   \int _{   \mathbb{R}   ^{2d}} \phi(X_s^\epsilon   -x')  \big(y'-Y_s^\epsilon   \big)
q_s^\epsilon    (\mathrm{d}x',\mathrm{d}y')   -\int _{   \mathbb{R}   ^{2d}}  \phi(\widehat X_s^{\epsilon,k}-x')
(y'-\widehat Y_s^{\epsilon,k})\widehat q_s^{\epsilon,k}(\mathrm{d}  x',\mathrm{d}  y') \bigg) \bigg]\\
= {}& \mathbb{E}  \bigg[  \mathbb{\tilde{E}}  \bigg[\phi(X_s^{\epsilon }    -  \tilde{X}_s^{\epsilon } )\Big(Y_s^\epsilon-\widehat{Y}_s^{\epsilon,k} \Big)   \cdot \Big(\Big(\tilde{Y}_s^{\epsilon } -Y_s^{\epsilon } \Big)-\Big(\tilde{\widehat{Y}}_s^{\epsilon,k} -   \widehat{Y} _s^{\epsilon,k}\Big)\Big)  \\
 {}&+\Big(\phi(X_s^{\epsilon }-\tilde{X}_s^{\epsilon })-
\phi(X_{t_k}^{\epsilon }-\tilde{X}_{t_k}^{\epsilon } )\Big)\Big(Y_s^\epsilon-\widehat{Y}_s^{\epsilon,k} \Big)   \cdot \Big(\tilde{\widehat{Y}}_s ^{\epsilon,k}-\widehat{Y}_s^{\epsilon,k}\Big)\bigg]\bigg]\\
= {}&  -   \frac{1}{2}\mathbb{E}  \bigg[  \mathbb{\tilde{E}}  \bigg[\phi(X_s^{\epsilon }    -  \tilde{X}_s^{\epsilon } )    \Big|  \Big(Y_s^{\epsilon }-   \widehat{Y} _s^{\epsilon,k} \Big)-\Big(\tilde{Y}_s^{\epsilon } -\tilde{\widehat{Y}}_s^{\epsilon,k}\Big)\Big|^{2}   \bigg]    \bigg] \\
 {}&+\mathbb{E}  \bigg[  \mathbb{\tilde{E}}  \bigg[ \Big(\phi(X_s^{\epsilon }-\tilde{X}_s^{\epsilon })-
\phi(X_{t_k}^{\epsilon }-\tilde{X}_{t_k}^{\epsilon } )\Big)\Big(Y_s^\epsilon-\widehat{Y}_s^{\epsilon,k} \Big)   \cdot \Big(\tilde{\widehat{Y}}_s ^{\epsilon,k}-\widehat{Y}_s^{\epsilon,k}\Big)\bigg]\bigg]\\
\leq   {} &   \mathbb{E}  \bigg[    \Big(\phi(X_s^{\epsilon }-\tilde{X}_s^{\epsilon })-
\phi(X_{t_k}^{\epsilon }-\tilde{X}_{t_k}^{\epsilon } )\Big)\Big(Y_s^\epsilon-\widehat{Y}_s^{\epsilon,k} \Big)   \cdot \Big(\tilde{\widehat{Y}}_s ^{\epsilon,k}-\widehat{Y}_s^{\epsilon,k}\Big) \bigg], 
\end{align*}
respectively.
By Young's inequality, we have
\begin{align*}
&\mathbb{E}  \bigg[    \Big(\phi(X_s^{\epsilon }-\tilde{X}_s^{\epsilon })-
\phi(X_{t_k}^{\epsilon }-\tilde{X}_{t_k}^{\epsilon } )\Big)\Big(Y_s^\epsilon-\widehat{Y}_s^{\epsilon,k} \Big)   \cdot \Big(\tilde{\widehat{Y}}_s ^{\epsilon,k}-\widehat{Y}_s^{\epsilon,k}\Big) \bigg]\\
\leq  {}&    \xi  \mathbb{E}  \Big[|Y_s^\epsilon-\widehat{Y}_s^{\epsilon,k}|^2\Big]   +  C_{\xi}   \mathbb{E}  \bigg[  \Big|\phi(X_s^{\epsilon }-\tilde{X}_s^{\epsilon })-
\phi(X_{t_k}^{\epsilon }-\tilde{X}_{t_k}^{\epsilon } )\Big|^{2} \Big|\tilde{\widehat{Y}}_s ^{\epsilon,k}-\widehat{Y}_s^{\epsilon,k}\Big|^{2} \bigg]
\end{align*}
and  
\begin{align*}
&\epsilon^{-1/2} \mathbb{E}  \bigg[ (Y_s^\epsilon-\widehat{Y}_s^{\epsilon,k} )   \cdot   G_{\rho_s^\epsilon}(X_s^\epsilon  )\bigg]\\
\leq  {}&   \frac{ \xi  }{\epsilon}\mathbb{E}  \Big[|Y_s^\epsilon-\widehat{Y}_s^{\epsilon,k}|^2\Big]   +  C_{\xi}  \mathbb{E}  \Big[|G_{\rho_s^\epsilon}(X_s^\epsilon  )|^2\Big]
\end{align*}
where $\xi  > 0$ is chosen later.
From Cauchy--Schwarz's inequality, Assumption \textbf{(A1)},   \eqref{eq:frozen-Y4}  and Lemma \ref{lem:X-increment}, it follows that
\begin{align*}
&\mathbb{E}  \bigg[    \Big|\phi(X_s^{\epsilon }-\tilde{X}_s^{\epsilon })-
\phi(X_{t_k}^{\epsilon }-\tilde{X}_{t_k}^{\epsilon } )\Big|^{2} \Big|\tilde{\widehat{Y}}_s ^{\epsilon,k}-\widehat{Y}_s^{\epsilon,k}\Big|^{2}\bigg]   \\
\leq   {} &    \bigg(\mathbb{E}  \bigg[    \Big|\phi(X_s^{\epsilon }-\tilde{X}_s^{\epsilon })-
\phi(X_{t_k}^{\epsilon }-\tilde{X}_{t_k}^{\epsilon } )\Big|^{4} \bigg] \bigg)^{1/2}
\bigg( \mathbb{E}  \bigg[    \Big|\tilde{\widehat{Y}}_s ^{\epsilon,k}-\widehat{Y}_s^{\epsilon,k}\Big|^{4} \bigg]\bigg)^{1/2}  \\
\leq   {} &   C_{T}\bigg(\mathbb{E}  \Big[  |X_s^{\epsilon }  - X_{t_k}^{\epsilon }|^{2}  \Big]+ \mathbb{ E}   \Big[ | \tilde{X}_s^{\epsilon }  - \tilde{X}_{t_k}^{\epsilon }   |^{2}\Big]
  \bigg)^{1/2}\\
 \leq  {}&  C_T(\Delta_\epsilon )^{1/2}.
\end{align*}
Taking $\xi$ sufficiently small  and   due to   $\sup_s    \mathbb{E} \big[|G_{\rho_s^\epsilon}(X_s^\epsilon)|^2 \big]  \le C_T$,  we obtain
\begin{align*}
&\mathrm{d}   \mathbb{E} \big[  |Y_s^\epsilon-\widehat{Y}_s^{\epsilon,k}|^2\big]  \\
\leq {}&  -\frac{c }{ \epsilon}    \mathbb{E} \big[  |Y_s^\epsilon-\widehat{Y}_s^{\epsilon,k}|^2\big]\mathrm{d}t
+\frac{  C_T}{\epsilon}(\Delta_\epsilon )^{1/2}\mathrm{d}t
+C_{T}\mathrm{d}t.
\end{align*}
 Then   differential inequality  yields  
\begin{align*}
& \mathbb{E} \big[  |Y_s^\epsilon-\widehat{Y}_s^{\epsilon,k}|^2\big]\\
\le   {} &    \Big(\frac{  C_T}{\epsilon}(\Delta_\epsilon)^{1/2}
+C_{T}  \Big)  \int_{s}^{t_{k}}   e^{   -\frac{c}{\epsilon } (s  -r )     }\mathrm{d}  r  \\
\le  {}&    \frac{\epsilon}{ c}  \Big(\frac{  C_T}{\epsilon}(\Delta_\epsilon)^{1/2}
+C_{T}  \Big)  \\
\le  {}&C_T((\Delta_\epsilon)^{1/2}+ \epsilon).
\end{align*}

\endproof

Next, we give the sequential stability of $\mathcal Q$ and $\overline{\mathcal{Q}}$. Note that the index $j$ in the following lemma is merely a sequential index for the objects to be compared; it is neither a particle label nor time.

\begin{lemma} \label{lem:Q-stability}
For each $j  \in  \mathbb{N}$, let
$
 q_j = \mathrm{Law}(X_j,Y_j),$
 $\widehat{q}_j= \mathrm{Law}(\widehat{X}_j,\widehat{Y}_j),
$
whose spatial marginals are denoted by $\rho_j,\widehat{\rho}_j$, respectively. Assume that   \textup{\textbf{(A1)}-\textbf{(A2)}}    hold, and  the four random variables are placed on the same coupling space, and as $j\to\infty$,
$$
\mathbb{E} \big[ |X_j-\widehat X_j|^2\big] +\mathbb{E} \big[|Y_j-\widehat Y_j|^2 \big]   \rightarrow      0,
$$
and 
$$
\sup_j\mathbb{E} \big[  |X_j|^2+|\widehat X_j|^2+|Y_j|^4+|\widehat Y_j|^4\big]<\infty.
$$
Then
\begin{equation}\label{eq:Q-stability}
\mathbb{E} \big[ | \mathcal{Q}(q_j;X_j,Y_j)-
 \mathcal{Q}(\widehat{q}_j;\widehat{X}_j,\widehat{Y}_j)|  \big] \rightarrow     0, \quad  \text{as}     \quad j  \to  \infty.
\end{equation}
Likewise, if $(X_j,\widehat X_j)$ are coupled with $\rho_j,  \widehat{\rho}_j$ such that
$$
\mathbb{E} \big[  |X_j-\widehat X_j|^2\big] \rightarrow    0,
 \qquad
 \sup_j\mathbb{E} \big[ |X_j|^2+|\widehat X_j|^2\big]<\infty,
$$
then
\begin{equation}\label{eq:Qbar-stability}
\mathbb{E} \big[ |\overline{\mathcal{Q}}_{\rho_j}(X_j)- \overline{\mathcal{Q}}_{\widehat{\rho}_j}(\widehat{X}_j)| \big] \rightarrow   0, \quad  \text{as}     \quad j  \to  \infty.
\end{equation}
\end{lemma}

\proof

 Let
$
 \Gamma_j= \mathrm{Law}(X_j,Y_j,\widehat{X}_j,\widehat{Y}_j) 
$
and $\pi_j := \mathrm{Law}(X_j, \widehat{X}_j)$ be the marginal distribution of $\Gamma_j$ on the coordinates $(X_j,\widehat{X}_j)$.
On a countable product probability space, take independent and identically distributed coupled copies
$
 (X_j^\alpha,Y_j^\alpha,
   \widehat{X}_j^\alpha,\widehat{Y}_j^\alpha),$
$ \alpha=0,1,2,\ldots,$
where each quadruple has law $\Gamma_j$, and the copies for different $\alpha$ are mutually independent. Thus, for each fixed $\alpha$,
\begin{equation}\label{eq:copy-coupling-convergence}
  \mathbb{E} \Big[ |X_j^\alpha-\widehat X_j^\alpha|^2   \Big]  +
  \mathbb{E} \Big[ |Y_j^\alpha-\widehat Y_j^\alpha|^2\Big]  \rightarrow  0, \quad  \text{as}     \quad j  \to  \infty,
\end{equation}
and each fixed finite set of copies inherits 
the uniform second moments of the positions and the fourth moments of the velocities from the original assumptions.
The following fact will be used repeatedly below:  If $A_j-\widehat{A}_j\to0$ in probability as $j  \to  \infty$, and $\{ \mathrm{Law}(A_j), \mathrm{Law}(\widehat{A}_j):   j\ge1\}$ is uniformly tight, then for each $f\in \mathcal{C}_b(\mathbb{R}^m)$ and $1\le p<\infty$,
\begin{equation}\label{eq:Cb-transfer}
 \mathbb{E}  \Big[|f(A_j)-f(\widehat{A}_j)|^p\Big]  \rightarrow  0, \quad  \text{as}     \quad j  \to  \infty.
 \end{equation}

 For any fixed $\alpha,\beta$, set
$
 A_j^{\alpha\beta}:=X_j^\alpha-X_j^\beta,$
 $\widehat{A}_j^{\alpha\beta}
 :=\widehat{X}_j^\alpha-\widehat{X}_j^\beta.
$
By $|a+b|^2\le2|a|^2+2|b|^2$ and \eqref{eq:copy-coupling-convergence},  we have 
\begin{align}
  \mathbb{E}  \Big[  |A_j^{\alpha\beta}-\widehat{A}_j^{\alpha\beta}|^2\Big]
 &\le
 2  \mathbb{E}\Big[ |X_j^\alpha-\widehat{X}_j^\alpha|^2\Big]
 +2  \mathbb{E}\Big[ |X_j^\beta-\widehat{X}_j^\beta|^2\Big]
 \rightarrow   0,  \quad  \text{as}     \quad j  \to  \infty. \label{eq:pair-difference-L2}
\end{align}
Under the assumptions of the lemma, there exists a constant $C >  0$ such that
$
 \sup_j  \mathbb{E}  \Big[ |A_j^{\alpha\beta}|^2\Big]  \le   C$  and 
 $\sup_j  \mathbb{E}  \Big[ |\widehat A_j^{\alpha\beta}|^2\Big]   \le C.
$
Consequently, Markov's inequality gives
$$
 \sup_j\mathbb{P}(|A_j^{\alpha\beta}|>   \xi )
 +\sup_j\mathbb{P}(|\widehat{A}_j^{\alpha\beta}|>\xi)
 \le 2  C   \xi^{-2}  \rightarrow 0,  \quad  \text{as}     \quad \xi  \to  \infty.   
$$
Thus $\{ \mathrm{Law}(A_j^{\alpha\beta}), \mathrm{Law}(\widehat{A}_j^{\alpha\beta}):   j\ge1\}$  is uniformly tight, and \eqref{eq:pair-difference-L2} further implies that $A_j^{\alpha\beta}-\widehat{A}_j^{\alpha\beta}  \to  0$  in probability. Under Assumption \textbf{(A1)}, we have $\phi\in \mathcal{C}_b$ and each $\partial_a \phi  \in \mathcal{C}_b$. Therefore, \eqref{eq:Cb-transfer} yields, for every finite $p$,
\begin{align}
 \mathbb{E}  \Big[|\phi(A_j^{\alpha\beta})
       -\phi(\widehat A_j^{\alpha\beta})|^p\Big]\to 0   \quad  \text{and }     \quad  \mathbb{E}  \Big[|\partial_a\phi(A_j^{\alpha\beta})
       -\partial_a\phi(\widehat A_j^{\alpha\beta})|^p\Big]  \to 0, \quad  \text{as}     \quad j  \to  \infty,
 \label{eq:phi-transfer-detailed}
\end{align}
which  implies  
\begin{align}
 \mathbb{E} \Big[| \nabla \phi(A_j^{\alpha\beta})
       - \nabla \phi(\widehat A_j^{\alpha\beta})|^p\Big]  \to 0, \quad  \text{as}     \quad j  \to  \infty.
 \label{eq:phi-grad-transfer-detailed}
\end{align}

 Fix a  copy \(\alpha\), and take an auxiliary copy \(\gamma\) independent of it, then we have
 \begin{align*}
 \lambda_{\rho_j}(X_j^\alpha)
 &=1+ \mathbb{E}_{\gamma}  [\phi(X_j^\alpha-X_j^\gamma)],
 \\
 \nabla\lambda_{\rho_j}(X_j^\alpha)
 &= \mathbb{E}_{\gamma}[\nabla\phi(X_j^\alpha-X_j^\gamma)] ,
\\
 \beta_{q_j}(X_j^\alpha)
 &= \mathbb{E}_{\gamma} \left[
 \nabla\phi(X_j^\alpha-X_j^\gamma)\cdot Y_j^\gamma\right],
\end{align*}
 where $ \mathbb{E}_\gamma$ denotes expectation with respect to the $\gamma$-th full coupled copy only.
 By   Jensen's inequality  and   \eqref{eq:phi-transfer-detailed}, for any $p \geq 1$, we have
 \begin{align}\label{eq:lambda-stability}
 & \mathbb{E} \Big[|\lambda_{\rho_j}(X_j^\alpha)
 -\lambda_{\widehat\rho_j}(\widehat X_j^\alpha)|^p\Big]   \notag\\
=  {}& \mathbb{E} \Big[\Big|\mathbb{E}_{\gamma}  [\phi(X_j^\alpha-X_j^\gamma)-   \phi(\widehat{X}_j^\alpha-\widehat{X}_j^\gamma)  ]
\Big|^p  \Big]  \notag\\
\le  {}&
  \mathbb{E} \Big [\mathbb{E}_{\gamma}   \Big [ \Big|\phi(X_j^\alpha-X_j^\gamma)
 -\phi(\widehat X_j^\alpha-\widehat X_j^\gamma)\Big|^p\Big]\Big]    \notag \\
 \le  {}&
  \mathbb{E} \Big [  \Big|\phi(X_j^\alpha-X_j^\gamma)
 -\phi(\widehat X_j^\alpha-\widehat X_j^\gamma)\Big|^p \Big] \notag\\
 &  \rightarrow0, \quad  \text{as}     \quad j  \to  \infty.
\end{align}
 Similarly, using \eqref{eq:phi-grad-transfer-detailed},  we obtain
 \begin{align}
  \mathbb{E} \Big [  |\nabla\lambda_{\rho_j}(X_j^\alpha)
 -\nabla\lambda_{\widehat\rho_j}(\widehat X_j^\alpha)|^p   \Big]  \rightarrow0, \quad  \text{as}     \quad j  \to  \infty.
 \label{eq:lambda-grad-stability}
\end{align}
 Since all $\lambda_{\rho}(x)  \ge 1$, the map $a\mapsto a^{-1}$ is 1-Lipschitz on $[1,\infty)$, it follows that
 \begin{equation}\label{eq:inverse-lambda-stability}
\mathbb{E} \Big [    |\lambda_{\rho_j}^{-1}(X_j^\alpha)
 -\lambda_{\widehat{\rho}_j}^{-1}(\widehat{X}_j^\alpha)|^p  \Big]  \rightarrow0, \quad  \text{as}     \quad j  \to  \infty.
\end{equation}
From \eqref{eq:phi-grad-transfer-detailed},   Assumption \textbf{(A1)}  and  Cauchy--Schwarz's inequality,  we have
\begin{align}\label{eq:beta-L2-stability}
&  \mathbb{E} \Big [     | \beta_{q_j}(X_j^\alpha)
 -\beta_{\widehat q_j}(\widehat X_j^\alpha) |^{2}    \Big]   \notag\\
 \leq  {}& 2\mathbb{E} \Big [  \mathbb{E}_\gamma   \Big[
| \bigl(\nabla\phi(X_j^\alpha-X_j^\gamma)
 -\nabla\phi(\widehat X_j^\alpha-\widehat X_j^\gamma)\bigr)
 \cdot Y_j^\gamma   |^{2}   \Big]   \Big] \notag\\
    {}&+ 2  \mathbb{E} \Big [  \mathbb{E}_\gamma    \Big[
| \nabla\phi(\widehat X_j^\alpha-\widehat X_j^\gamma)
 \cdot(Y_j^\gamma-\widehat Y_j^\gamma)  |^{2} \Big]  \Big]   \notag\\
 \leq  {}& 2   \Big( \mathbb{E} \Big [|\nabla\phi(X_j^\alpha-X_j^\gamma)
 -\nabla\phi(\widehat X_j^\alpha-\widehat X_j^\gamma)|^4  \Big]
 \Big)^{1/2}
\Big(            \mathbb{E} \Big [ |Y_j^\gamma|^4         \Big]        \Big)^{1/2}  + 2 L_\phi^2   \mathbb{E} \Big [ |Y_j^\gamma-\widehat Y_j^\gamma|^2  \Big]   \notag\\
  &  \rightarrow0, \quad  \text{as}     \quad j  \to  \infty.
\end{align}
 From    Assumption \textbf{(A1)}, it follows that
\begin{equation}\label{eq:beta-uniform-bound-stability}
 \sup_j   \bigl(\|\beta_{q_j}\|_\infty+
 \|\beta_{\widehat{q}_j}\|_\infty\bigr)
 \le L_\phi\sup_j\Big(
  \mathbb{E} \Big [|Y_j|\Big]+ \mathbb{E} \Big [   |\widehat Y_j|\Big]  \Big)<\infty.
\end{equation}
Set
$
 b(X_j^\alpha) :=\frac{\beta_{q_j}(X_j^\alpha)}
 {\lambda_{\rho_j}(X_j^\alpha)}, $
 $ b (\widehat{X}_j^\alpha):=
 \frac{\beta_{\widehat{q}_j}(\widehat{X}_j^\alpha)}
 {\lambda_{\widehat{\rho}_j}(\widehat{X}_j^\alpha)}.
$
By  using  \eqref{eq:inverse-lambda-stability}, 
\eqref{eq:beta-L2-stability}   and 
\eqref{eq:beta-uniform-bound-stability},   we have  
\begin{align}\label{eq:b-L2-stability-detailed}
& \mathbb{E} \Big [  |b(X_j^\alpha)  -b (\widehat{X}_j^\alpha)|^2\Big]   \notag  \\
 =  {}& \mathbb{E} \bigg [  \bigg|\frac{\beta_{q_j}(X_j^\alpha)}
 {\lambda_{\rho_j}(X_j^\alpha)}   -\frac{\beta_{\widehat q_j}(\widehat X_j^\alpha)}
 {\lambda_{\widehat\rho_j}(\widehat X_j^\alpha)}  \bigg|^2\bigg]  \notag \\ 
  =  {}& \mathbb{E} \bigg [ \Big|\lambda_{\rho_j}^{-1}(X_j^\alpha)( \beta_{q_j}(X_j^\alpha)   -  \beta_{\widehat q_j}(\widehat X_j^\alpha)    )
 + (     \lambda_{\rho_j}^{-1}(X_j^\alpha)  -   \lambda_{\widehat\rho_j}^{-1}(\widehat X_j^\alpha)    )  \beta_{\widehat q_j}(\widehat X_j^\alpha) 
  \Big|^2\bigg]  \notag\\
    \leq    {}& 2 \mathbb{E} \bigg [ \Big|  \beta_{q_j}(X_j^\alpha)   -  \beta_{\widehat q_j}(\widehat X_j^\alpha)     \Big|^2\bigg]
 + 2   \|\beta_{\widehat{q}_j}\|_\infty^{2} \mathbb{E} \bigg [ \Big|      \lambda_{\rho_j}^{-1}(X_j^\alpha)  -   \lambda_{\widehat\rho_j}^{-1}(\widehat X_j^\alpha)      
  \Big|^2\bigg]  \notag\\
 &  \rightarrow0, \quad  \text{as}     \quad j  \to  \infty 
 \end{align}
 and 
 \begin{align}\label{eq:b-uniform-bound-detailed}
 \sup_j\bigl(\| b(X_j^\alpha)   \|_{L^\infty(\Omega)}
 +\|b (\widehat{X}_j^\alpha)  \|_{L^\infty(\Omega)}\bigr)<\infty.
\end{align}

Note that,  for the fixed $m$ product factors,
\begin{equation}\label{eq:finite-product-telescope}
 \prod_{\ell=1}^m \Theta_{ j, \ell} 
 -\prod_{\ell=1}^m  \widehat{\Theta}_{j,\ell}
 =
 \sum_{\ell=1}^m
 \left(\prod_{h<\ell}\Theta_{j,h}\right)
 (\Theta_{j,\ell}-\widehat{\Theta}_{j,\ell})
 \left(\prod_{h>\ell}\widehat{\Theta}_{j,h}\right).
\end{equation}
If all factors are bounded in absolute value by $K$, then Minkowski's inequality yields
\begin{equation}\label{eq:finite-product-Lp}
 \left\|\prod_{\ell=1}^m \Theta_{j,\ell}
 -\prod_{\ell=1}^m\widehat{\Theta}_{j,\ell}\right\|_{L^p}
 \le K^{m-1}\sum_{\ell=1}^m
 \|\Theta_{j,\ell}-\widehat{\Theta}_{j,\ell}\|_{L^p}.
\end{equation}
For $0\le k<n$, $0\le \ell\le n$, and coordinate $1\le a\le d$, denote
\begin{align*}
 \Phi_{j,k}&  :=\phi(X_j^k-X_j^{k+1}),&
 \widehat{\Phi}_{j,k}
 &:=\phi(\widehat{X}_j^k-\widehat{X}_j^{k+1}),\\
 \dot{\Phi}_{j,k}^a
 &:=\partial_a\phi(X_j^k-X_j^{k+1}),&
 \widehat{\dot{\Phi}}_{j,k}^a
 &:=\partial_a\phi(\widehat{X}_j^k-\widehat{X}_j^{k+1}),\\
 R_{j,\ell}&:=\lambda_{\rho_j}^{-1}(X_j^\ell),&
 \widehat{R}_{j,\ell}
 &:=\lambda_{\widehat{\rho}_j}^{-1}(\widehat{X}_j^\ell),\\
 H_{j,\ell}^a&:=  \partial_a\lambda_{\rho_j}(X_j^\ell)    \lambda_{\rho_j}^{-1}(X_j^\ell),&
 \widehat{H}_{j,\ell}^a
 &:=\partial_a\lambda_{\widehat{\rho}_j}(\widehat{X}_j^\ell)  \lambda_{\widehat{\rho}_j}^{-1}(\widehat{X}_j^\ell).
\end{align*}
By   \eqref{eq:lambda-grad-stability},  \eqref{eq:inverse-lambda-stability}  and   Assumption \textbf{(A1)},  we have 
 \begin{align*}
 & \mathbb{E} \Big [  | \partial_a\lambda_{\rho_j}(X_j^\ell)    \lambda_{\rho_j}^{-1}(X_j^\ell)  - 
  \partial_a\lambda_{\widehat{\rho}_j}(\widehat{X}_j^\ell)  \lambda_{\widehat{\rho}_j}^{-1}(\widehat{X}_j^\ell)|^p   \Big] \\
\leq  {}&  p \mathbb{E} \Big [  | \lambda_{\rho_j}^{-1}(X_j^\ell)  (\partial_a\lambda_{\rho_j}(X_j^\ell)    
 -\partial_a\lambda_{\widehat{\rho}_j}(\widehat{X}_j^\ell) ) |^p    \Big]   + p \mathbb{E} \Big [  | \partial_a\lambda_{\widehat{\rho}_j}(\widehat{X}_j^\ell) (  \lambda_{\rho_j}^{-1}(X_j^\ell)   - \lambda_{\widehat{\rho}_j}^{-1}(\widehat{X}_j^\ell))|^p   \Big]\\  
 \leq {} & C \mathbb{E} \Big [  |  \partial_a\lambda_{\rho_j}(X_j^\ell)    
 -\partial_a\lambda_{\widehat{\rho}_j}(\widehat{X}_j^\ell)   |^p    \Big]   + C \mathbb{E} \Big [  | \lambda_{\rho_j}^{-1}(X_j^\ell)   - \lambda_{\widehat{\rho}_j}^{-1}(\widehat{X}_j^\ell) |^p   \Big]\\
  & \rightarrow0, \quad  \text{as}     \quad j  \to  \infty,
\end{align*}
which  implies    
 \begin{align} \label{eq:lambda-grad-stability}
 \|H_{j,\ell}^a  -    \widehat{H}_{j,\ell}^a \|_{L^{p}(\pi_j)}   \rightarrow0, \quad  \text{as}     \quad j  \to  \infty.
 \end{align}
From \eqref{eq:weights},    it follows that
$$
 w_{n,\rho_j}(X_j^{0}, \dots,   X_j^{n})
 =\Bigg(\prod_{k=0}^{n-1}\Phi_{j,k}\Bigg)
  \Bigg(\prod_{\ell=0}^{n}R_{j,\ell}\Bigg),
$$
By 
\eqref{eq:finite-product-telescope},  \eqref{eq:finite-product-Lp} and  Assumption \textbf{(A1)},  we get
\begin{align}\label{eq:w-finite-stability}
 & \Big(\mathbb{E}\Big[|w_{n,\rho_j}(X_j^{0}, \dots,   X_j^{n})
 -w_{n,\widehat\rho_j}(\widehat{X}_j^{0}, \dots,   \widehat{X}_j^{n}) |^p \Big]\Big)^{1/p}   \notag\\
 =  {}&  \Bigg\|\Bigg(\prod_{k=0}^{n-1}\Phi_{j,k}\Bigg)
  \Bigg(\prod_{\ell=0}^{n}R_{j,\ell}\Bigg)      -  \Bigg(\prod_{k=0}^{n-1}\widehat{\Phi}_{j,k}\Bigg)
  \Bigg(\prod_{\ell=0}^{n}\widehat{R}_{j,\ell}\Bigg)   \Bigg\|_{L^{p}(\pi_j^{\otimes (n+1)})}   \notag\\ 
  \leq    {}&  \sum_{i=0}^{n-1}
\left\|  \left(\prod_{k=0}^{i-1}\Phi_{j,k}\right)
  (\Phi_{j,i}-\widehat{\Phi}_{j,i})
  \left(\prod_{k=i+1}^{n-1}\widehat{\Phi}_{j,k}\right)
  \left(\prod_{\ell=0}^{n}R_{j,\ell}\right)
\right\|_{L^{p}\left(\pi_j^{\otimes (n+1)}\right)}  \notag\\
 {}&     +
\sum_{\ell=0}^{n}
\left\|   \left(\prod_{k=0}^{n-1}\widehat{\Phi}_{j,k}\right)
  \left(\prod_{h=0}^{\ell-1}R_{j,h}\right)
  (R_{j,\ell}-\widehat{R}_{j,\ell})
  \left(\prod_{h=\ell+1}^{n}\widehat{R}_{j,h}\right)
\right\|_{L^{p}\left(\pi_j^{\otimes (n+1)}\right)}   \notag\\
   \leq    {}&    C     \sum_{i=0}^{n-1}
\left\|  \Phi_{j,i}-\widehat{\Phi}_{j,i} \right\|_{L^{p} (\pi_j )}      +
   C   \sum_{\ell=0}^{n}
\left\|   R_{j,\ell}-\widehat{R}_{j,\ell} \right\|_{L^{p}(\pi_j )}     \notag\\
 & \rightarrow0, \quad  \text{as}     \quad j  \to  \infty,
\end{align}
where the last line follows from \eqref{eq:phi-transfer-detailed} and \eqref{eq:inverse-lambda-stability}.
Furthermore, from \eqref{eq:w-derivative-exact}, the $a$-th component of $g_{n,i,\rho_j}$ can be written as
\begin{align*}
 g_{n,i,\rho_j}^a(X_j^{0}, \dots,   X_j^{n})
 ={}&
 \boldsymbol1_{\{i<n\}}\dot\Phi_{j,i}^a
 \Bigg(\prod_{\substack{0\le k<n\\k\ne i}}\Phi_{j,k}\Bigg)
 \left(\prod_{\ell=0}^{n}R_{j,\ell}\right)\\
 &-\boldsymbol1_{\{i>0\}}\dot\Phi_{j,i-1}^a
 \Bigg(\prod_{\substack{0\le k<n\\k\ne i-1}}\Phi_{j,k}\Bigg)
 \left(\prod_{\ell=0}^{n}R_{j,\ell}\right)\\
 &-\left(\prod_{k=0}^{n-1}\Phi_{j,k}\right)
 \left(\prod_{\ell=0}^{n}R_{j,\ell}\right)H_{j,i}^a.
\end{align*}
Analogously, by  \eqref{eq:phi-transfer-detailed}, \eqref{eq:phi-grad-transfer-detailed}, \eqref{eq:inverse-lambda-stability} and \eqref{eq:finite-product-telescope}-\eqref{eq:lambda-grad-stability}, we have
\begin{equation}\label{eq:g-finite-stability}
\begin{aligned}
 &\mathbb{ E}\Big[|g_{n,i,\rho_j}(X_j^{0}, \dots,   X_j^{n})
 -g_{n,i,\widehat{\rho}_j}(\widehat{X}_j^{0}, \dots,  \widehat{ X}_j^{n})|^p \Big] \rightarrow0, \quad  \text{as}     \quad j  \to  \infty.
\end{aligned}
\end{equation}

For any fixed copy indices $\alpha$ and $\beta$, under the assumptions of the lemma  and Cauchy--Schwarz's inequality,  we have
\begin{align}
& \mathbb{ E}\Big[  \|Y_j^\alpha\otimes Y_j^\beta-
 \widehat Y_j^\alpha\otimes\widehat Y_j^\beta\|_{F}\Big]  \notag \\
 =   {}& \mathbb{ E}\Big[  \| (Y_j^\alpha-\widehat Y_j^\alpha)\otimes Y_j^\beta
 +\widehat Y_j^\alpha\otimes
 (Y_j^\beta-\widehat Y_j^\beta)\|_{F}  \Big]   \notag \\
\le   {} & \Big( \mathbb{E}\Big[|Y_j^\alpha-\widehat Y_j^\alpha|^2\Big]\Big)^{1/2}
      \Big( \mathbb{E}\Big[|Y_j^\beta|^2\Big]\Big)^{1/2}   +\Big( \mathbb{E}\Big[|\widehat Y_j^\alpha|^2\Big]\Big)^{1/2}
      \Big( \mathbb{E} \Big[  |Y_j^\beta-\widehat Y_j^\beta|^2\Big]\Big)^{1/2}   \notag \\
&   \rightarrow0, \quad  \text{as}     \quad j  \to  \infty
 \label{eq:velocity-product-L1}
\end{align}
 and
 \begin{align}\label{eq:velocity-product-uniform-bound}
\mathbb{ E}\Big[  \|   Y_j^\alpha\otimes Y_j^\beta\|_{F}^{2}   \Big]  
= {} &\mathbb{E}\Big[ |Y_j^\alpha|^2|Y_j^\beta|^2\Big]   \notag  \\
\le   {} &
 \Big(\mathbb{E}\Big[|Y_j^\alpha|^4\Big]\Big)^{1/2}
 \Big(\mathbb{E}\Big[|Y_j^\beta|^4\Big]\Big)^{1/2}\notag   \\
  \le   {}&C,
\end{align}
where  $\|\cdot\|_{F}$ denotes the Frobenius norm.
From \eqref{eq:Q-def}, we have
\begin{equation}\label{eq:QN-copy-formula}
\begin{aligned}
& \mathcal{Q}^{(N)}(q_j;X_j^0,Y_j^0)\\
= {}&\sum_{n=0}^N \mathbb{E}_{1:n} \left[
 Y_j^n\left(
 \sum_{i=0}^nY_j^i\cdot
 g_{n,i,\rho_j}(X_j^{0}, \dots,   X_j^{n})
 +w_{n,\rho_j}(X_j^{0}, \dots,   X_j^{n})\sum_{i=0}^n   b(X_j^i)
 \right)\right],
\end{aligned}
\end{equation}
where   we adopt the convention that   the $n=0$ term is interpreted as $Y_j^0\bigl[Y_j^0\cdot \nabla R_{\rho_j}(X_j^0)+  R_{\rho_j}(X_j^0) b(X_j^0)  \bigr]$. 
In the following steps, we always first fix $0\le n\le N   \in   \mathbb{N}$ and $0\le i\le n$, and then let $j\to\infty$.
By  using  
\eqref{eq:g-finite-stability}-\eqref{eq:velocity-product-uniform-bound},   Assumption \textbf{(A1)} and Cauchy--Schwarz's inequality,   we  obtain 
\begin{align}
& \mathbb{E}\Big[\Big| \mathbb{E}_{1:n}\Big[
  Y_j^n    (Y_j^i  \cdot g_{n,i,\rho_j}(X_j^{0}, \dots,   X_j^{n}))   -\widehat{Y}_j^n ( \widehat{Y}_j^i  \cdot g_{n,i,\widehat{\rho}_j}(\widehat{X}_j^{0}, \dots,   \widehat{X}_j^{n}))\Big]\Big|\Big]\notag\\
  =  {} & \mathbb{E}\Big[\Big| \mathbb{E}_{1:n}\Big[
  Y_j^n\otimes Y_j^i  g_{n,i,\rho_j}(X_j^{0}, \dots,   X_j^{n})   -\widehat{Y}_j^n\otimes\widehat{Y}_j^i  g_{n,i,\widehat{\rho}_j}(\widehat{X}_j^{0}, \dots,   \widehat{X}_j^{n})\Big]\Big|\Big]\notag\\
\le   {} &   \mathbb{E}\Big[|(Y_j^n\otimes Y_j^i) (g_{n,i,\rho_j}(X_j^{0}, \dots,   X_j^{n}) - g_{n,i,\widehat{\rho}_j}(\widehat{X}_j^{0}, \dots,   \widehat{X}_j^{n}))| \Big] \notag\\
     {}&  + \mathbb{ E}\Big[|(Y_j^n\otimes Y_j^i    -\widehat{Y}_j^n\otimes\widehat{Y}_j^i ) g_{n,i,\widehat{\rho}_j}(\widehat{X}_j^{0}, \dots,   \widehat{X}_j^{n})|\Big]  \notag\\
\le   {}&  \Big( \mathbb{E}\Big[ |g_{n,i,\rho_j}(X_j^{0}, \dots,   X_j^{n}) - g_{n,i,\widehat{\rho}_j}(\widehat{X}_j^{0}, \dots,   \widehat{X}_j^{n}) |^{2}   \Big]\Big)^{1/2}    \Big( \mathbb{E}\Big[  \|Y_j^n\otimes Y_j^i \|_{F}^2  \Big]\Big)^{1/2} \notag \\
     {} & +  C   \mathbb{E}\Big[  \|Y_j^n\otimes Y_j^i    -\widehat{Y}_j^n\otimes\widehat{Y}_j^i\|_{F}\Big] \notag \\
& \rightarrow0, \quad  \text{as}     \quad j  \to  \infty,
\label{eq:derivative-chain-fixed-stability}
\end{align}
where   we have  used $  (Y_j^n\otimes Y_j^i)g_{n,i,\rho_j}(X_j^{0}, \dots,   X_j^{n})=Y_j^n(Y_j^i\cdot g_{n,i,\rho_j}(X_j^{0}, \dots,   X_j^{n})).$
By virtue of \eqref{eq:b-L2-stability-detailed}, \eqref{eq:b-uniform-bound-detailed}, \eqref{eq:w-finite-stability},  Assumption \textbf{(A1)} and Cauchy--Schwarz's inequality,  we derive 
  \begin{align}
  & \mathbb{E}\Big[\Big| \mathbb{E}_{1:n}\Big[w_{n,\rho_j}(X_j^{0}, \dots,   X_j^{n})    Y_j^n    b(X_j^i)  -   w_{n,\widehat\rho_j}(\widehat{X}_j^{0}, \dots,   \widehat{X}_j^{n})   \widehat{Y}_j^n  b(\widehat{X}_j^i)\Big]\Big|\Big]\notag\\
\leq  {} & \mathbb{E}\Big[|w_{n,\rho_j}(X_j^{0}, \dots,   X_j^{n})    Y_j^n    b(X_j^i)  -   w_{n,\widehat\rho_j}(\widehat{X}_j^{0}, \dots,   \widehat{X}_j^{n})   \widehat{Y}_j^n  b(\widehat{X}_j^i)|   \Big]
\notag\\
\leq   {}&  \mathbb{E}\Big[|(w_{n,\rho_j}(X_j^{0}, \dots,   X_j^{n})   - w_{n,\widehat\rho_j}(\widehat{X}_j^{0}, \dots,   \widehat{X}_j^{n})  )Y_j^n   b(X_j^i) |   \Big]
 +\mathbb{E}\Big[|w_{n,\widehat\rho_j}(\widehat{X}_j^{0}, \dots,   \widehat{X}_j^{n})(Y_j^n-\widehat{Y}_j^n) b(X_j^i)|   \Big]\notag\\
 {}&+\mathbb{E}\Big[|w_{n,\widehat\rho_j}(\widehat{X}_j^{0}, \dots,   \widehat{X}_j^{n})    \widehat{Y}_j^n   (b(X_j^i)  -  b(\widehat{X}_j^i))|   \Big]\notag  \\
\le {} & C   \Big( \mathbb{E}\Big[ |w_{n,\rho_j}(X_j^{0}, \dots,   X_j^{n})   - w_{n,\widehat\rho_j}(\widehat{X}_j^{0}, \dots,   \widehat{X}_j^{n})|^{2}   \Big]\Big)^{1/2}      \Big( \mathbb{E}\Big[ |Y_j^n|^{2}   \Big]\Big)^{1/2}
 +C\Big( \mathbb{E}\Big[ |Y_j^n-  \widehat{Y}_j^n|^{2}   \Big]\Big)^{1/2}  \notag   \\
 {} &   +   C  \Big( \mathbb{E}\Big[ |\widehat{Y}_j^n|^{2}   \Big]\Big)^{1/2}      \Big( \mathbb{E}\Big[ |b(X_j^i)  -  b(\widehat{X}_j^i)|^{2}   \Big]\Big)^{1/2}    \notag \\
& \rightarrow0, \quad  \text{as}     \quad j  \to  \infty.
\label{eq:beta-chain-fixed-stability}
\end{align}
  From \eqref{eq:QN-copy-formula},  \eqref{eq:derivative-chain-fixed-stability} and \eqref{eq:beta-chain-fixed-stability}, we obtain
\begin{equation}\label{eq:QN-finite-stability}
 \mathbb{E}\Big[ | \mathcal{Q}^{(N)}(q_j;X_j,Y_j)
 - \mathcal{Q}^{(N)}(\widehat q_j;\widehat X_j,\widehat Y_j)| \Big]
\rightarrow 0, \quad  \text{as}     \quad j  \to  \infty,
\end{equation}
which holds for each fixed $N   \in   \mathbb{N}$.
 From
  Lemma   \ref{lem:kernel-sum}, we have 
  \begin{align*}
 \sup_j  \mathbb{E}\Big[| \mathcal{Q}(q_j;X_j,Y_j)- \mathcal{Q}^{(N)}(q_j;X_j,Y_j)|\Big]
 &\le C\Big(1+\Big(   \mathbb{E}\Big[| Y_j   |^{4}\Big]    \Big)^{1/2}\Big)\sum_{n>N}(n+1)^2r_*^{(n-2)_+},\\
 \sup_j    \mathbb{E}\Big[   | \mathcal{Q}(\widehat{q}_j;\widehat{X}_j,\widehat{Y}_j)
 - \mathcal{Q}^{(N)}(\widehat{q}_j;\widehat{X}_j,\widehat{Y}_j)|\Big]
 &\le C\Big(1+\Big(   \mathbb{E}\Big[| \widehat{Y}_j   |^{4}\Big]    \Big)^{1/2}\Big)\sum_{n>N}(n+1)^2r_*^{(n-2)_+}.
\end{align*}
  Together with \eqref{eq:QN-finite-stability}, we get
  \begin{align*}
 &\limsup_{j\to\infty}
 \mathbb{E}\Big[| \mathcal{Q}(q_j;X_j,Y_j)
 - \mathcal{Q}(\widehat{q}_j;\widehat{X}_j,\widehat{Y}_j)|\Big]\\
 \le &  C\Big(1+      \Big(   \mathbb{E}\Big[| Y_j   |^{4}\Big]    \Big)^{1/2} +   \Big(   \mathbb{E}\Big[| \widehat{Y}_j   |^{4}\Big]    \Big)^{1/2}  \Big)\sum_{n>N}(n+1)^2r_*^{(n-2)_+},
\end{align*}
  which holds for each fixed $N   \in   \mathbb{N}$.
  Finally, taking the limit $N\to\infty$, we get
$$
\mathbb{E} \big[ | \mathcal{Q}(q_j;X_j,Y_j)-
 \mathcal{Q}(\widehat{q}_j;\widehat{X}_j,\widehat{Y}_j)|  \big] \rightarrow     0, \quad  \text{as}     \quad j  \to  \infty.
$$

  From \eqref{eq:Qbar}, we have
 $$
  \overline{\mathcal{Q}}_{\rho_j}^{(N)}(X_j^0)
 =\sum_{n=0}^N     \mathbb{E}_{1:n} \big[
J_{\rho_j}(X_j^n)
 g_{n,n,\rho_j}(X_j^0,\ldots,X_j^n)\big].
$$
where we adopt the convention that   the $n=0$ term is interpreted as $J_{\rho_j}(X_j^0)   \nabla R_{\rho_j}(X_j^0)$.  
Due   to \eqref{eq:inverse-lambda-stability},  we  get 
  \begin{align}\label{eq:J-stability-detailed}
  \mathbb{ E}\Big[  \| J_{\rho_j}(X_j^n)-
 J_{\widehat{\rho}_j}(\widehat{X}_j^n)\|_{F}^{p}\Big]  \rightarrow     0, \quad  \text{as}     \quad j  \to  \infty.
\end{align}
Using 
 \eqref{eq:J-stability-detailed},     \eqref{eq:g-finite-stability}  and Cauchy--Schwarz's inequality,  we  have 
\begin{align}
&\mathbb{E}\Big[\Big| \mathbb{E}_{1:n}\Big[    J_{\rho_j}(X_j^n)   g_{n,n,\rho_j}(X_j^0,\ldots,X_j^n)-   J_{\widehat{\rho}_j}(\widehat{X}_j^n)  g_{n,n,\widehat{\rho}_j}(\widehat{X}_j^0,\ldots,\widehat{X}_j^n)\Big]\Big|\Big]\notag\\
\leq  {}&  \mathbb{E}\Big[\Big|    ( J_{\rho_j}(X_j^n)   -    J_{\widehat{\rho}_j}(\widehat{X}_j^n) ) g_{n,n,\rho_j}(X_j^0,\ldots,X_j^n)  \Big|\Big]    \notag    \\
  {}&   +    \mathbb{E}\Big[\Big|  J_{\widehat{\rho}_j}(\widehat{X}_j^n)    (g_{n,n,\rho_j}(X_j^0,\ldots,X_j^n)  -g_{n,n,\widehat{\rho}_j}(\widehat{X}_j^0,\ldots,\widehat{X}_j^n)) \Big|\Big]\notag\\
 \leq  {}&  C   \Big( \mathbb{E}\Big[\Big\|     J_{\rho_j}(X_j^n)   -    J_{\widehat{\rho}_j}(\widehat{X}_j^n)   \Big\|_{F}^{2}\Big] \Big)^{1/2}   \notag    \\
  {}&   +   C   \Big( \mathbb{E}\Big[\Big|   g_{n,n,\rho_j}(X_j^0,\ldots,X_j^n)  -g_{n,n,\widehat{\rho}_j}(\widehat{X}_j^0,\ldots,\widehat{X}_j^n)  \Big|^{2}\Big]  \Big)^{1/2}  \notag\\
 &  \rightarrow     0, \quad  \text{as}     \quad j  \to  \infty.
\end{align}
  Then   we  get  
$$
\mathbb{E}\Big[  |\overline{\mathcal{Q}}_{\rho_j}^{(N)}(X_j)
 -\overline{\mathcal{Q}}_{\widehat\rho_j}^{(N)}(\widehat X_j)|\Big]  \rightarrow     0, \quad  \text{as}     \quad j  \to  \infty,
$$ 
 which holds for each fixed $N   \in   \mathbb{N}$. 
 By  \eqref{eq:Qbar-uniform-bound-early},  we have 
$$
 \sup_\rho\|\overline{\mathcal{Q}}_\rho-\overline{\mathcal{Q}}_\rho^{(N)}\|_\infty
 \le C\sum_{n>N}  r_*^{n-1},
$$
 which holds for each fixed $N   \in   \mathbb{N}$.     Further, we deduce
$$
 \limsup_{j\to\infty}
 \mathbb{E}\Big[  |  \overline{\mathcal{Q}}_{\rho_j}(X_j)
 -\overline{\mathcal{Q}}_{\widehat{\rho}_j}(\widehat{X}_j)|\Big]
 \le    C\sum_{n>N}r_*^{n-1}.
$$
Finally, letting $N\to\infty$, we obtain
$$
\mathbb{E} \big[ |\overline{\mathcal{Q}}_{\rho_j}(X_j)- \overline{\mathcal{Q}}_{\widehat{\rho}_j}(\widehat{X}_j)| \big] \rightarrow   0, \quad  \text{as}     \quad j  \to  \infty.
$$

\endproof

\begin{remark} \label{rem:G-stability}
  From Assumption \textup{\textbf{(A2)}}, we obtain
$$
 |G_{\rho_j}(X_j)-G_{\widehat{\rho}_j}(\widehat X_j)|
 \le \operatorname{Lip}(\nabla K)
      \bigl(|X_j   -  \widehat X_j|+\mathbb{W}_1(\rho_j,\widehat{\rho}_j)\bigr).
$$
If $(X_j,\widehat{X}_j)$ is a coupling of $(\rho_j,\widehat{\rho}_j)$ such that
$X_j-\widehat{X}_j\to0$ in probability, and
$$
 \sup_j  \mathbb{E} \big[ |X_j|^2+|\widehat{X}_j|^2\big]<\infty,
$$
then
\begin{align}\label{eq:G-L1-stability-detailed}
\mathbb{E} \big[|G_{\rho_j}(X_j)-G_{\widehat{\rho}_j}(\widehat X_j)|\big]
 &  \rightarrow   0, \quad  \text{as}     \quad j  \to  \infty.
\end{align}
Note  that 
\begin{align*}
 F_{\rho_{j} }^{(N)} (X_{j} ^0)    =      \sum_{n=0}^N    \mathbb{E}_{1:n}   \left[
 w_{n,\rho_ j  }(X_{j} ^0,\ldots,X_{j} ^n)G_{\rho_{j} }(X_{j} ^n)\right],
 \end{align*}
 where   we adopt the convention that   the $n=0$ term is interpreted as  $R_{\rho_{j}} G_{\rho_{j}} (X_{j} ^0)$. 
 By virtue of    Assumption \textup{\textbf{(A1)}}, \eqref{eq:G-L1-stability-detailed}  and Cauchy--Schwarz's inequality     
  \begin{align}
  & \mathbb{E}\Big[\Big| \mathbb{E}_{1:n}\Big[w_{n,\rho_j}(X_j^{0}, \dots,   X_j^{n})  G_{\rho_{j} }(X_{j} ^n)  -   w_{n,\widehat\rho_j}(\widehat{X}_j^{0}, \dots,   \widehat{X}_j^{n})   G_{\widehat{\rho}_{j} }(\widehat{X}_{j} ^n)\Big]\Big|\Big]\notag\\
\leq   {}& \mathbb{E}\Big[|w_{n,\rho_j}(X_j^{0}, \dots,   X_j^{n})   G_{\rho_{j} }(X_{j} ^n) -   w_{n,\widehat\rho_j}(\widehat{X}_j^{0}, \dots,   \widehat{X}_j^{n})   G_{\widehat{\rho}_{j} }(\widehat{X}_{j} ^n)|   \Big]
\notag\\
\leq  {} &  \mathbb{E}\Big[|(w_{n,\rho_j}(X_j^{0}, \dots,   X_j^{n})   - w_{n,\widehat\rho_j}(\widehat{X}_j^{0}, \dots,   \widehat{X}_j^{n})  ) G_{\rho_{j} }(X_{j} ^n)  |   \Big]\notag \\
  {}&+\mathbb{E}\Big[|w_{n,\widehat\rho_j}(\widehat{X}_j^{0}, \dots,   \widehat{X}_j^{n})(G_{\rho_{j} }(X_{j} ^n) - G_{\widehat{\rho}_{j} }(\widehat{X}_{j} ^n)) |   \Big]\notag\\
\le {} &    C \Big( \mathbb{E}\Big[ |w_{n,\rho_j}(X_j^{0}, \dots,   X_j^{n})   - w_{n,\widehat\rho_j}(\widehat{X}_j^{0}, \dots,   \widehat{X}_j^{n})|^{2}   \Big]\Big)^{1/2}      \Big( \mathbb{E}\Big[ |\widehat{X}_j^n|^{2}   \Big]\Big)^{1/2}
    \notag   \\
  {}&   +   C      \mathbb{E}\Big[ |G_{\rho_{j} }(X_{j} ^n) - G_{\widehat{\rho}_{j} }(\widehat{X}_{j} ^n)|    \Big]     \notag \\
& \rightarrow0, \quad  \text{as}     \quad j  \to  \infty.
\label{eq:G-chain-fixed-stability}
\end{align}
 Similar to the discussion in Lemma \ref{lem:Q-stability} 
  \begin{align}\label{eq:G-stability}
\mathbb{E} \big[ | F_{\rho_{j} } (X_{j} ^0)-  F_{\widehat{\rho}_{j} }  (\widehat{X}_{j} ^0)  | \big] \rightarrow   0, \quad  \text{as}     \quad j  \to  \infty.
\end{align}

 In addition, from  \eqref{eq:inverse-lambda-stability} 
 \begin{equation}\label{eq:inverse-2-lambda-stability}  
  \mathbb{ E }  [   |   \lambda_{\rho_{j}}^{-2}( X_{j} ^0)  -    \lambda_{\widehat{\rho}_{j}}^{-2}(\widehat{X}_{j} ^0 ) | ] \to0,   \quad  \text{as}     \quad     j  \to  \infty. 
  \end{equation}

\end{remark}

\begin{remark}\label{rem:G-Lip-stability}
For each $j  \in  \mathbb{N}$, let
$
 q_j = \mathrm{Law}(X_j,Y_j),$
 $\widehat{q}_j= \mathrm{Law}(\widehat{X}_j,\widehat{Y}_j),
$
whose spatial marginals are denoted by $\rho_j,\widehat{\rho}_j$, respectively. Assume that  \textup{\textbf{(A1\textquotesingle)}-\textbf{(A2\textquotesingle)}}   hold, and  the four random variables are placed on the same coupling space, and 
$
\sup_j\mathbb{E} \big[  |X_j|^2+|\widehat X_j|^2+|Y_j|^4+|\widehat Y_j|^4\big]<\infty.
$
Then 
$$
 \mathbb{E}[| \mathcal{Q}(q_j;X_j,Y_j)- \mathcal{Q}(\widehat{q}_j;\widehat{X}_j,\widehat{Y}_j)|]
\le C (  
( \mathbb{E}[|X_j-\widehat{X}_j|^2])^{1/2}
+( \mathbb{E}[ |Y_j-\widehat{Y}_j|^2]  )^{1/2}) .
$$
Likewise, if $(X_j,\widehat X_j)$ are coupled with $\rho_j,  \widehat{\rho}_j$ such that
$
 \sup_j\mathbb{E} \big[ |X_j|^2+|\widehat X_j|^2\big]<\infty,
$
then
\begin{equation*} 
\mathbb{E} \big[ |\overline{\mathcal{Q}}_{\rho_j}(X_j)- \overline{\mathcal{Q}}_{\widehat{\rho}_j}(\widehat{X}_j)| \big] \le C    
( \mathbb{E}[|X_j-\widehat{X}_j|^2])^{1/2}.
\end{equation*}
 
\end{remark}

\medskip

\section{Global Lipschitz estimate under assumptions \textbf{(A1\textquotesingle)}-\textbf{(A2\textquotesingle)}}

\begin{lemma}\label{lem:limit-Lipschitz}
Under \textup{\textbf{(A1\textquotesingle)}-\textbf{(A2\textquotesingle)}}, for $\rho,\widehat{\rho} \in\mathcal{P}_1(\mathbb{R}^d)$ and $X ^0  , \widehat{X}^0\in\mathbb{R}^d$,
\begin{equation}\label{eq:full-coeff-Lip}
 |b_\rho(X ^0)-b_{\widehat{\rho}}(\widehat{X}^0)|
 + |\lambda_\rho^{-1}  (X ^0)  \sigma - \lambda_{\widehat{\rho}}^{-1}  (\widehat{X}^0)  \sigma |
 \le C\bigl(| X ^0   -\widehat{X}^0|+\mathbb{W}_1(\rho,\widehat{\rho})\bigr).
\end{equation}
\end{lemma}

\proof

Under Assumption \textbf{(A1\textquotesingle)} 
\begin{equation*} 
 |\lambda_\rho(X ^0)-\lambda_{\widehat{\rho}}(\widehat{X}^0)|
 \le L_\phi\bigl(| X ^0   -\widehat{X}^0|+\mathbb{W}_1(\rho,\widehat{\rho})\bigr),
\end{equation*}
and further it follows that
\begin{equation} \label{eq:inverse-lambda-Lipschitz}
|\lambda_\rho^{-1}  (X ^0)  \sigma - \lambda_{\widehat{\rho}}^{-1}  (\widehat{X}^0)  \sigma |
 \le C\bigl(| X ^0   -\widehat{X}^0|+\mathbb{W}_1(\rho,\widehat{\rho})\bigr).
\end{equation}
Under Assumption \textbf{(A1\textquotesingle)}, we have
$$
 |\phi(X ^0)-\phi(\widehat{X}^0)|\le L_\phi|X ^0-\widehat{X}^0|,
 \qquad
 |\nabla\phi(X ^0)-\nabla\phi(\widehat{X}^0)|
 \le  \mathrm{Lip}(\nabla\phi)|X ^0-\widehat{X}^0|,
$$
and further  
$$
 |\nabla\lambda_\rho(X ^0)-\nabla\lambda_{\widehat{\rho}}(\widehat{X}^0)|
 \le   \mathrm{Lip}(\nabla\phi)
      \bigl(| X ^0   -\widehat{X}^0|+\mathbb{W}_1(\rho,\widehat{\rho})\bigr).
$$

  From \eqref{eq:Qbar}, we have
 $$
  \overline{\mathcal{Q}}_{\rho }^{(N)}(X ^0)
 =\sum_{n=0}^N     \mathbb{E}_{1:n} \big[
J_{\rho }(X ^n)
 g_{n,n,\rho }(X ^0,X ^1, \ldots,X ^n)\big].
$$
where we adopt the convention that   the $n=0$ term is interpreted as $J_\rho(X ^0 )   \nabla R_\rho(X ^0 )$.  
Due   to \eqref{eq:inverse-lambda-Lipschitz},  we  get 
  \begin{align*} 
   \| J_{\rho }( X ^0)-
 J_{ \widehat{\rho}}(\widehat{X}^0)\|_{F}   \le C\bigl(| X ^0   -\widehat{X}^0|+\mathbb{W}_1(\rho,\widehat{\rho})\bigr).
\end{align*}
For $0\le k<n$, $0\le \ell\le n$, and coordinate $1\le a\le d$, denote
\begin{align*}
 \Phi_{k}&  :=\phi(X^k-X^{k+1}),&
 \widehat{\Phi}_{k}
 &:=\phi(\widehat{X}^k-\widehat{X}^{k+1}),\\
 \dot{\Phi}_{k}^a
 &:=\partial_a\phi(X^k-X^{k+1}),&
 \widehat{\dot{\Phi}}_{k}^a
 &:=\partial_a\phi(\widehat{X}^k-\widehat{X}^{k+1}),\\
 R_{\ell}&:=\lambda_{\rho}^{-1}(X^\ell),&
 \widehat{R}_{\ell}
 &:=\lambda_{\widehat{\rho}}^{-1}(\widehat{X}^\ell),\\
 H_{\ell}^a&:=  \partial_a\lambda_{\rho}(X^\ell)    \lambda_{\rho}^{-1}(X^\ell),&
 \widehat{H}_{\ell}^a
 &:=\partial_a\lambda_{\widehat{\rho}}(\widehat{X}^\ell)  \lambda_{\widehat{\rho}}^{-1}(\widehat{X}^\ell).
\end{align*}
Then we can obtain 
 \begin{align*}
| \Phi_{k} -  \widehat{\Phi}_{k}|\leq  {}&  C (| X^k  -\widehat{X}^k  |   +   | X^{k+1}  -\widehat{X}^{k+1}  |  ), \\
|\dot{\Phi}_{k}^a  -\widehat{\dot{\Phi}}_{k}^a  |\leq  {}& C (| X^k  -\widehat{X}^k  |   +   | X^{k+1}  -\widehat{X}^{k+1}  |  ),   \\
 |R_{\ell}  - \widehat{R}_{\ell}|\leq  {}&    C\bigl(| X ^\ell   -\widehat{X}^\ell|+\mathbb{W}_1(\rho,\widehat{\rho})\bigr)
\end{align*}
and 
 \begin{align*}
| H_{\ell}^a  - \widehat{H}_{\ell}^a | =  {}&    | \partial_a\lambda_{\rho }(X ^\ell)    \lambda_{\rho }^{-1}(X ^\ell)  - 
  \partial_a\lambda_{\widehat{\rho}}(\widehat{X}^\ell)  \lambda_{\widehat{\rho}}^{-1}(\widehat{X}^\ell)|   \\
\leq  {}&     | \lambda_{\rho}^{-1}(X^\ell)  (\partial_a\lambda_{\rho}(X^\ell)    
 -\partial_a\lambda_{\widehat{\rho}}(\widehat{X}^\ell) ) |   +    | \partial_a\lambda_{\widehat{\rho}}(\widehat{X}^\ell) (  \lambda_{\rho}^{-1}(X^\ell)   - \lambda_{\widehat{\rho}}^{-1}(\widehat{X}^\ell))| \\  
 \leq  {}& C\bigl(| X ^\ell   -\widehat{X}^\ell|+\mathbb{W}_1(\rho,\widehat{\rho})\bigr).
\end{align*}
Furthermore, from \eqref{eq:w-derivative-exact}, the $a$-th component of $g_{n,i,\rho}$ can be written as
\begin{align*}
 g_{n,i,\rho}^a(X^{0}, \dots,   X^{n})
 ={}&
 \boldsymbol{1}_{\{i<n\}}\dot\Phi_{i}^a
 \Bigg(\prod_{\substack{0\le k<n\\k\ne i}}\Phi_{k}\Bigg)
 \left(\prod_{\ell=0}^{n}R_{\ell}\right)\\
 &-\boldsymbol{1}_{\{i>0\}}\dot\Phi_{i-1}^a
 \Bigg(\prod_{\substack{0\le k<n\\k\ne i-1}}\Phi_{k}\Bigg)
 \left(\prod_{\ell=0}^{n}R_{\ell}\right)\\
 &-\left(\prod_{k=0}^{n-1}\Phi_{k}\right)
 \left(\prod_{\ell=0}^{n}R_{\ell}\right)H_{i}^a.
\end{align*}
  By similar arguments as in the proof of Lemma \ref{lem:Q-stability}, Assumption \textbf{(A1\textquotesingle)} and  \eqref{eq:finite-product-telescope},    we have
\begin{equation*} 
\begin{aligned}
 &\mathbb{E}_{1:n}\Big[|g_{n,i,\rho}(X^{0}, \dots,   X^{n})
 -g_{n,i,\widehat{\rho}_j}(\widehat{X}^{0}, \dots,  \widehat{ X}^{n})|  \Big]    \leq C (n+1)^2  r_{*}^{(n-2)_{+}} \bigl(| X ^0   -\widehat{X}^0|+\mathbb{W}_1(\rho,\widehat{\rho})\bigr) 
\end{aligned}
\end{equation*}
and
\begin{equation*} 
\begin{aligned}
   |g_{n,i,\rho}(X^{0}, \dots,   X^{n})  |       \leq C    r_{*}^{(n-2)_{+}},  
\end{aligned}
\end{equation*}
where    we enlarge it to the larger upper bound $C (n+1)^2  r_{*}^{(n-2)_{+}} \bigl(| X ^0   -\widehat{X}^0|+\mathbb{W}_1(\rho,\widehat{\rho})\bigr)$  in order to ensure that it holds for all  $n\ge 0$.
Consequently, we have
\begin{align*}
& \Big| \mathbb{E}_{1:n} \Big[    J_{\rho }(X ^n)   g_{n,n,\rho }(X^{0},  \ldots,X ^n)-   J_{\widehat{\rho} }(\widehat{X} ^n)  g_{n,n, \nu }(  \widehat{X}^{0} , \ldots,\widehat{X} ^n)\Big]\Big| \notag\\
\leq  {}&  \mathbb{E}_{1:n}\Big[\Big|    ( J_{\rho}(X^n)   -    J_{\widehat{\rho} }(\widehat{X}^n) ) g_{n,n,\rho}(X^0,\ldots,X^n)  \Big|\Big]    \notag    \\
  {}&   +   \mathbb{E}_{1:n}\Big[\Big|  J_{\widehat{\rho}}(\widehat{X}^n)    (g_{n,n,\rho}(X^0,\ldots,X^n)  -g_{n,n,\widehat{\rho}}(\widehat{X} ^0,\ldots,\widehat{X}^n)) \Big|\Big]\notag\\
 \leq  {}&   C (n+1)^2  r_{*}^{(n-2)_{+}} \bigl(| X ^0   -\widehat{X}^0|+\mathbb{W}_1(\rho,\widehat{\rho})\bigr). 
\end{align*}
  Then   we  get  
$$
   |\overline{\mathcal{Q}}_{\rho }^N(X^0 )
 -\overline{\mathcal{Q}}_{\widehat{\rho} }^N(\widehat{X}^0 )|  \leq C \bigl(| X ^0   -\widehat{X}^0|+\mathbb{W}_1(\rho,\widehat{\rho})\bigr)\sum_{n=0}^{N} (n+1)^2  r_{*}^{(n-2)_{+}},
$$ 
 which  implies 
$$
 |\overline{\mathcal{Q}}_{\rho }(X^0 )- \overline{\mathcal{Q}}_{\widehat{\rho} }(\widehat{X}^0)| \leq C \bigl(| X ^0   -\widehat{X}^0|+\mathbb{W}_1(\rho,\widehat{\rho})\bigr).
$$
  From Assumption \textbf{(A2\textquotesingle)}, we obtain
$$
 |G_{\rho}(X^0)-G_{\widehat{\rho}}(\widehat{X}^0)|
 \le  \mathrm{Lip}(\nabla K)
      \bigl(|X^0   -  \widehat{X}^0|+\mathbb{W}_1(\rho,\widehat{\rho})\bigr) 
$$
and 
$$
|G_{\rho}(X^0)|\leq  C .
$$
 By similar arguments as in    Remark \ref{rem:G-stability}, we  get 
$$
| F_{\rho  } (X  ^0)-  F_{\widehat{\rho}  }  (\widehat{X}  ^0)  |  \le   
     C  \bigl(|X^0   -  \widehat{X}^0|+\mathbb{W}_1(\rho,\widehat{\rho})\bigr) .
$$
Hence
\begin{equation*} 
 |b_\rho(X ^0)-b_{\widehat{\rho}}(\widehat{X}^0)|
 \le C\bigl(| X ^0   -\widehat{X}^0|+\mathbb{W}_1(\rho,\widehat{\rho})\bigr).
\end{equation*}

\endproof


\end{document}